\input ./style/arxiv-general.cfg
\documentclass[MSNbibl,number,citesort,dvips]{arxbj}
\makeatletter
   \@ifpackageloaded{graphicx}{}{\usepackage{graphicx}}
\makeatother
\usepackage{mathbh}
\usepackage{upgreek}

\aid{0}
\volume{22}
\issue{1}
\pubyear{2016}
\firstpage{143}
\lastpage{192}
\doi{10.3150/14-BEJ626} 

\makeatletter
\renewcommand{\pi}{\uppi}
\newtheorem{theorem}{Theorem}[section]
\newtheorem{proposition}[theorem]{Proposition}
\newtheorem{corollary}[theorem]{Corollary}
\newtheorem{lemma}[theorem]{Lemma}
\newproclaim{assumption}{Assumption}
\newremark{remark}[theorem]{Remark}
\newcommand{\E}{{\mathbb E}}
\newcommand{\F}{{\mathcal F}}
\newcommand{\R}{{\mathbb R}}
\newcommand{\N}{{\mathbb N}}
\newcommand{\Q}{{\mathbb Q}}
\newcommand{\supp}{\operatorname{supp}}
\newcommand{\argmin}{\mathop{\arg\min}}
\newcommand{\Var}{\operatorname{Var}}
\newcommand{\bull}{{\bolds{\cdot}}}
\renewcommand{\d}{{\mathrm{d}}}
\renewcommand{\phi}{\varphi}
\renewcommand{\theta}{\vartheta}
\renewcommand{\subset}{\subseteq}
\newcommand{\ind}{\mathbh{1}}
\newcommand{\B}{\mathcal B}
\makeatother

\begin{document}
\begin{frontmatter}

\title{Adaptive quantile estimation in deconvolution with unknown
error~distribution}
\runtitle{Quantile estimation in deconvolution}

\begin{aug}
\author[t1]{\inits{I.}\fnms{Itai}~\snm{Dattner}\corref{}\ead[label=e1]{idattner@stat.haifa.ac.il}\thanksref{t1}},
\author[t2]{\inits{M.}\fnms{Markus}~\snm{Rei\textup{\ss}}\ead[label=e2,mark]{mreiss@math.hu-berlin.de}\thanksref{t2,e2}} \and
\author[t2]{\inits{M.}\fnms{Mathias}~\snm{Trabs}\ead[label=e3,mark]{trabs@math.hu-berlin.de}\thanksref{t2,e3}}
\runauthor{I. Dattner, M. Rei{\ss} and M. Trabs} 
\address[t1]{Department of Statistics,
University of Haifa,
199 Abba Khoushy Ave,
Mount Carmel, Haifa 3498838, Israel.
\printead{e1}}
\address[t2]{Institut f\"ur Mathematik,
Humboldt-Universit\"at zu Berlin,
Unter den Linden 6,
10099 Berlin, Germany.
\printead{e2,e3}}
\end{aug}

\received{\smonth{5} \syear{2013}}
\revised{\smonth{4} \syear{2014}}

%
\begin{abstract}
Quantile estimation in deconvolution problems is studied
comprehensively. In particular, the more realistic setup of unknown
error distributions is covered. Our plug-in method is based on a
deconvolution density estimator and is minimax optimal under minimal
and natural conditions.
This closes an important gap in the literature. Optimal adaptive
estimation is obtained by a data-driven bandwidth choice. As a side
result, we obtain optimal rates for the plug-in estimation of
distribution functions with unknown error distributions. The method is
applied to a real data example.
\end{abstract}

%
\begin{keyword}
\kwd{adaptive estimation}
\kwd{deconvolution}
\kwd{distribution function}
\kwd{minimax convergence rates}
\kwd{plug-in estimator}
\kwd{quantile function}
\kwd{random Fourier multiplier}
\end{keyword}
\end{frontmatter}

\section{Introduction}\label{secintro}
Nonparametric deconvolution models are of high practical importance and
lead to challenging questions in statistical methodology. Let $X_1,
\ldots, X_n$ be independent random variables with a common Lebesgue
density $f\dvtx \R\to\R$. Suppose that we merely observe the random variables
\begin{eqnarray}
\nonumber
Y_j=X_j+\varepsilon_j,
\qquad j=1,\ldots,n,
\end{eqnarray}
that is the original $(X_j)$ corrupted by i.i.d. error variables
$\varepsilon_j$, independent of $(X_j)$ and with Lebesgue density
$f_\varepsilon$. For $\tau\in(0,1)$ the objective is to estimate the
$\tau$-quantile $q_\tau$ of the population $X$
from the observations $Y_1, \ldots, Y_n$. For practitioners estimated
quantiles are very relevant, but they depend in a nonlinear way on the
underlying density such that their estimation is not always obvious.
Abstractly, quantile estimation in deconvolution is an example of
nonlinear functional estimation in ill-posed inverse problems.

Two natural strategies may be pursued. Either a distribution function
estimator is inverted or an M-estimation paradigm is applied using a
density estimator of $f$. While the first possibility was studied by
Hall and Lahiri \cite{hall2008estimation}, the purpose of this paper
is the analysis of the second in a far more general setting.
Assuming that the distribution of the measurement error is completely
known, Carroll and Hall \cite{carroll1988optimal} have constructed a
kernel density estimator based on the empirical characteristic function
$\phi_{n}(u):=\frac{1}{n}\sum_{j=1}^n\mathrm{e}^{\mathrm{i}uY_j}, u\in\R$. In
practice, however, the distribution of the measurement error is usually
not known. Instead, we assume that we have at hand a sample from
$f_\varepsilon$ given by
\begin{eqnarray}
\nonumber
\varepsilon^*_1,\ldots,\varepsilon^*_m,
\qquad m\in\N.
\end{eqnarray}
Motivated from applications, we will not assume that the observations
$(\varepsilon^*_k)$ are independent from~$(Y_j)$. In particular, our
procedure applies to the experimental setup of repeated measurements,
as discussed below.

Let $\F g(u):=\int_{\R} \mathrm{e}^{\mathrm{i}ux}g(x)\,\mathrm{d}x$, $u\in\R$, denote the
Fourier transform of $g\in L^1(\R)\cup L^2(\R)$. Consequently, $\F
^{-1}[h(u)](x)=\frac{1}{2\pi}\int \mathrm{e}^{-\mathrm{i}ux}h(u)\,\d u,x\in\R$. Based
on the classical kernel estimator, Neumann \cite{neumann1997effect}
has proposed the following density estimator of $f$ for the case of
unknown error distributions:
%
\[
\widetilde f_b(x):=\F^{-1} \biggl[
\frac{\phi_n(u)\phi_K(bu)}{\phi
_{\varepsilon,m}(u)} {\ind_{\{|\phi_{\varepsilon,m}(u)|\geq
m^{-1/2}\}}} \biggr](x),\qquad x\in\R,
\]
where $\phi_K$ is the Fourier transform of a kernel $K$, $b>0$ is its
bandwidth and the characteristic function of the error distribution
$\phi_\varepsilon$ is estimated by its empirical counterpart $\phi
_{\varepsilon,m}(u):=\frac{1}{m}\sum_{k=1}^m\mathrm{e}^{\mathrm{i}u\varepsilon_k^*},
u\in\R$. Obviously, $\widetilde f_b$ depends on the sample sizes $n$
and $m$ which are suppressed in the notation. Applying a plug-in
approach, our estimator for the quantile $q_\tau$ is then given by the
minimum-contrast estimator
%
\begin{equation}
\label{eqDefEst} \widetilde q_{\tau,b}:=\argmin_{\eta\in
[-U_n,U_n]}\bigl|\widetilde
M_b(\eta)\bigr|\qquad\mbox{with } \widetilde M_b(\eta)=
\int_{-\infty
}^{\eta}\widetilde f_b(x)\,\d x-
\tau
\end{equation}
for some $U_n\to\infty$. 
We will show as the very first step that $\widetilde f_b$ is indeed
integrable with overwhelming probability and when not, we define
$\widetilde q_{\tau}$ to be the empirical $\tau$-quantile of the
observations $Y_j$'s. In this work we pursue the analysis for error
distributions whose characteristic function decays polynomially. As
shown by Fan \cite{fan1991}, these so-called ordinary smooth errors
lead to mildly ill-posed estimation problems. They are mathematically
more challenging than the so-called super-smooth errors, which we
discuss briefly in Section~\ref{SecDisc}.

Although the literature on deconvolution problems is extensive and very
broad, the problem of adaptive deconvolution with unknown measurement
errors was addressed only recently, see Comte and Lacour \cite
{comte2011data}, Johannes and Schwarz \cite{Johannes2010} and Kappus
\cite{kappus2012}
for adaptive density estimation with unknown error distributions in the
model selection framework. Minimax results and other properties for
nonadaptive methods are given by Neumann \cite
{neumann1997effect,neumann2007deconvolution}, Meister \cite
{meister2004effect}, Delaigle, Hall and Meister
\cite{delaigle2008deconvolution}, Johannes \cite
{johannes2009deconvolution} among others. To the best of our knowledge,
the problem of quantile estimation in deconvolution was considered only
in Hall and Lahiri \cite{hall2008estimation}. They have constructed a
quantile estimator for the case of known error distributions by
inverting the distribution function estimator, without proposing an
adaptive bandwidth choice. As we shall establish, the error of the quantile
estimator (\ref{eqDefEst}) is directly related to that of the
distribution function estimator (cf. the error representation (\ref
{eqtaylor}) below). Yet, the general analysis of the latter was not
clear before.

Fan \cite{fan1991} has proposed an estimator for the distribution
function by integrating the density deconvolution estimator. In order
to perform an exact analysis of its variance, a truncation of the
integral was required in the estimation procedure. This resulted in a
nonoptimal (in the minimax sense) estimation method for the case of
ordinary smooth errors and raised the conjecture that `plug-in does not
work optimally' for estimation of the distribution function in
deconvolution. Trying to circumvent this problem, Hall and Lahiri \cite
{hall2008estimation} as well as Dattner, Goldenshluger and Juditsky
\cite{dattnerEtAll2011} have constructed a distribution function
estimator based on a direct inversion formula.
Applying the Fourier multiplier approach by Nickl and Rei{\ss} \cite
{NicklReiss2012}, {S\"ohl} and Trabs \cite{soehlTrabs2012} have shown
that the integrated density estimator can indeed estimate the
distribution function with $\sqrt n$-rate under suitable conditions.
Since they prove a Donsker theorem, the imposed conditions are
restrictive. In particular, a global Sobolev regularity of $f$ is
assumed there which is not natural for pointwise loss.
So even with a known error distribution, it remained an open and
intriguing question whether the canonical plug-in estimator for
distribution or quantile function estimation yields asymptotically
optimal results under natural conditions.

In Section~\ref{secminimax}, we settle this question in the positive
under local H\"older regularity of $f$ by combining an exact analysis
like in Dattner, Goldenshluger and Juditsky \cite{dattnerEtAll2011}
together with abstract Fourier multiplier theory from {S\"ohl} and
Trabs \cite{soehlTrabs2012}. Moreover, we show that the optimal rates
continue to hold if the error distribution is unknown and has to be
estimated, which is mathematically nontrivial. Since the deconvolution
operator $\F^{-1}[1/\phi_\varepsilon]$ is not observable, we have to
study the estimated counterpart $\F^{-1}[ {\frac{\phi_K(bu)}{\phi
_{\varepsilon,m}(u)}\ind_{\{|\phi_{\varepsilon,m}(u)|\geq m^{-1/2}\}
}}]$. As a random Fourier multiplier, it preserves the mapping
properties of the deterministic $\F^{-1}[1/\phi_\varepsilon]$, but
its operator norm turns out to be (slightly) larger.

A lower bound result establishes that the rates under a local H\"older
condition are indeed minimax optimal. Surprisingly, the dependence of
the minimax rate on the error sample size $m$ is completely different
from the case of global Sobolev restrictions like in Neumann \cite
{neumann1997effect}.
The proof enlightens this interplay between the decay of one
characteristic function and estimation error in the other sample for
both, the $(Y_j)$ and the $(\varepsilon_j)$.

An adaptive (data-driven) bandwidth choice is developed in Section~\ref
{secadaptive}. To this end, a variant of Lepski's method is applied,
but because of the unknown and possibly dependent error distribution a
much more refined analysis is needed to establish that the resulting
adaptive quantile estimator is (up to log factors) still rate optimal.

In Section~\ref{secnumeric}, we implement our estimation procedure
and present simulation results which show a good performance of the
estimator. In a real data example, we consider multiple blood pressure
measurement data from different patients. Here, a measurement error is
clearly present, but of unknown distribution and we have to estimate it
by taking patient-wise differences. The completely data-driven method
yields reasonable quantile estimates which differ from the sample
quantiles of the directly measured $(Y_j)$. All proofs are postponed to
Section~\ref{secproofs}.

%
\section{Convergence rates}\label{secminimax}

\subsection{Setting and upper bounds}
Let us introduce some notation. Denoting $\langle\alpha\rangle$ as
the largest integer which is strictly smaller than $\alpha>0$, we
define for some function $g$ and any possibly unbounded interval
$I\subset\R$ the H\"older norm
%
\[
\|g\|_{C^{\alpha}(I)}:=\sum_{k=0}^{\langle\alpha\rangle}
\bigl\|g^{(k)}\bigr\| _{L^\infty(I)}+\sup_{x,y\in I\dvtx x\neq y}
\frac{|g^{\langle\alpha
\rangle}(x)-g^{\langle\alpha\rangle}(y)|}{|x-y|^{\alpha-\langle
\alpha\rangle}}.
\]
Let $C^0(I)$ denote the space of all continuous and bounded functions
on the interval $I$ and
\[
C^s(\R):=\bigcup_{R>0}C^s(
\R,R)\qquad\mbox{with } C^\alpha(I,R):= \bigl\{g\in C^0(I) |
\|g\|_{C^{\alpha}(I)}\leq R \bigr\}, R>0.
\]
In the sequel, we use the Landau notation $\mathcal{O}$ and $\mathcal{O}_P$. For two sequences $A_n(\theta), B_n(\theta)$ depending on a
parameter $\theta$, $A_n(\theta)=\mathcal{O}_P(B_n(\theta))$ holds
uniformly over a parameter set $\theta\in\Theta$ if there is for all
$c>0$ some $C>0$ such that $\sup_{\theta\in\Theta}P_\theta
(A_n(\theta)>CB_n(\theta))<c$. If $A_n(\theta)/B_n(\theta)$
converges in probability to zero, we write $A_n(\theta)=\mathrm{o}_P(B_n(\theta))$.

\renewcommand\theassumption{\Alph{assumption}}
%
\begin{assumption}\label{asKernel}
Let the kernel function $K\in L^1(\R)$ with Fourier transform $\phi
_K:=\F K$ satisfy
\begin{enumerate}[(ii)]
\item[(i)] $\supp{\phi_K}\subset[-1,1]$ and
\item[(ii)] $K$ has order $\ell\in\N$, i.e., for $k=0,\dots,
\ell$
\[
\int_{\R} \bigl|K(x)\bigr||x|^{\ell+1}\,\d x<\infty\quad
\mbox{and}\quad\int_{\R} x^kK(x)\,\d x =\cases{ 1,&
\quad if $k=0$,
\cr
0,&\quad otherwise.}
\]
\end{enumerate}
\end{assumption}

By construction the quantile estimator, $\widetilde q_{\tau,b}$ is the
approximated solution of the estimating equation
%
\begin{equation}
\label{eqestEq} 0=\widetilde M_b(\eta)=\int_{-\infty}^{\eta}
\widetilde f_b(x)\,\d x-\tau.
\end{equation}
If a solution exists, it does not have to be unique since $\widetilde
f_b$ is not necessarily nonnegative. Nevertheless, any choice\vspace*{1pt}
converges to the true quantile, assuming the latter is unique. Before,
integrability of $\widetilde f_b$ was an open problem, which we shall
settle now.

%
\begin{lemma}\label{lemWellDef}
Grant Assumption~\textup{\ref{asKernel}} with $\ell=0$. On the event
%
\begin{equation}
\label{eqBEps} B_\varepsilon(b):= \Bigl\{\inf_{u\in[-1/b,1/b]}\bigl|
\phi_{\varepsilon,m}(u)\bigr|\geq m^{-1/2}|\log b|^{3/2} \Bigr\}
\end{equation}
we have $\widetilde f_b\in L^1(\R)$ and estimating equation~(\ref
{eqestEq}) has a solution.
\end{lemma}

Therefore, a truncation of the integral as used by Fan \cite{fan1991}
is not necessary, implying that no tail condition on $f$ is required.
Although $\|\widetilde f_b\|_{L^1}$ is finite, it depends on the
observations as well as through $b$ on $n,m$. To quantify the behavior
of $\widetilde f_b$ more precisely, our analysis relies on the
following much stronger result.

\begin{lemma}\label{lemFourierMult}
Grant Assumption~\textup{\ref{asKernel}} with $\ell=0$. For some $\beta,
R>0$ suppose $\E[(\varepsilon_k^*)^4]\leq R$ and
\[
\bigl|\phi_\varepsilon(u)\bigr|^{-1}\leq R\bigl(1+|u|\bigr)^{\beta}\quad
\mbox{and}\quad\bigl|\phi_\varepsilon'(u)\bigr|\leq
R\bigl(1+|u|\bigr)^{-\beta-1}
\]
as well as $mb^{2\beta+1}\to\infty$. Then there exists a finite
random variable $\mathcal E_b$ which is $\mathcal{O}_P (1\vee\frac
{1}{m^{1/2}b^{\beta+1}} )$ with the constant depending only on
$\beta$ and $R$, such that for any $s>\beta^+>\beta$ on the event
$B_\varepsilon(b)$ from (\ref{eqBEps})
\begin{eqnarray*}
\biggl\|\F^{-1} \biggl[\frac{\phi_K(bu)}{\phi_{\varepsilon,m}(u)} \biggr]\ast
\psi\biggr\|_{C^{s-\beta^+}(\R)}
\leq\mathcal E_{b}\|\psi\|_{C^{s}(\R)}\qquad\mbox{for all }
\psi\in C^{s}(\R).
\end{eqnarray*}
\end{lemma}

The deterministic counterpart of this lemma was proved by {S\"ohl} and
Trabs \cite{soehlTrabs2012}. Here, we show that the random Fourier
multiplication operator $C^{s}(\R)\ni\psi\mapsto\F^{-1} [\frac
{\phi_K(bu)\F\psi(u)}{\phi_{\varepsilon,m}(u)} ]\in C^{s-\beta
^+}(\R)$ has a norm bound $\mathcal{O}_P (1\vee\frac
{1}{m^{1/2}b^{\beta+1}} )$ on the\vspace*{1pt} event $B_\varepsilon(b)$. The
condition on the derivative $\phi_\varepsilon'$ is natural in the
context of Fourier multipliers and is usually satisfied for
distributions with polynomial decaying characteristic functions, for
example, Gamma distributions with shape parameter $\beta>0$ satisfy it.

%
\begin{remark}
Depending only on the observations, condition (\ref{eqBEps}) can be
verified by the practitioner for a given bandwidth $b$. Under the
assumptions of Lemma~\ref{lemFourierMult} Talagrand's inequality %
yields $P(B_\varepsilon(b))\geq1- 2\mathrm{e}^{-mb^{2\beta+1}}$ (cf.
Lemma~\ref{lemPhiEps} and (\ref{eqTal2Phi}) below). Therefore, with
overwhelming probability $B_\varepsilon(b)$ holds true and the
estimating equation (\ref{eqestEq}) is rigorously defined.
\end{remark}

Before we start with the error analysis, let us describe the class of
densities we are interested in. Let $\mathcal Q(R)$ denote the set of
all probability densities on $\R$ which are uniformly bounded by
$R>0$. Following the minimax paradigm, we consider for $R,r,\zeta,U>0$
and the smoothness index $\alpha>0$ the classes
\begin{eqnarray*}
\mathcal C^\alpha(R,r,\zeta)&:=& {\bigcup_{U\in\N}
\mathcal C^\alpha(R,r,\zeta, U)}\quad\mbox{and}
\\
\mathcal C^\alpha(R,r,\zeta, {U})&:=& \bigl\{f\in\mathcal Q(R) | f\mbox{
has a }\tau\mbox{-quantile } {q_\tau\in[-U,U]}\mbox{ such that}
\nonumber
\\
&&\hspace*{5pt} f\in C^\alpha\bigl([q_\tau-
\zeta,q_\tau+\zeta],R\bigr)\mbox{ and }f(q_\tau)\geq r
\bigr\}.
\end{eqnarray*}
In contrast to Dattner, Goldenshluger and Juditsky \cite
{dattnerEtAll2011}, the smoothness is measured locally in a H\"older
scale and not globally by decay conditions of the Fourier transform of
$f$. The former is more natural since both, the distribution function
and the quantile function are estimated pointwise. Note that the
quantile $q_\tau$ is unique given the assumption $f(q_\tau)>0$.
Recalling that we write $\phi_\varepsilon:=\F f_\varepsilon$, the
conditions in Lemma~\ref{lemFourierMult} motivate the definition of
the class of error densities
\begin{eqnarray*}
\mathcal D^\beta(R,\gamma)&:=& \biggl\{f_\varepsilon\in\mathcal Q(
\infty) \Big|\frac{1}{R}\bigl(1+|u|\bigr)^{-\beta}\leq\bigl|\F f_\varepsilon(u)\bigr|
\leq R\bigl(1+|u|\bigr)^{-\beta},
\nonumber
\\
&&\hspace*{5pt} \bigl|(\F f_\varepsilon)'(u)\bigr|\leq R\bigl(1+|u|\bigr)^{-1-\beta},
\bigl\|x^{\gamma
}f_\varepsilon(x)\bigr\|_{L^1}\leq R \biggr
\}\label{eqClasserror}
\end{eqnarray*}
for some moment $\gamma\geq0$ and we use the same constant $R$
as above for convenience.

%
\begin{remark}
The upper and lower bounds for $|\phi_\varepsilon(u)|$ in $\mathcal
D^\beta(R,\gamma)$ are standard assumptions in deconvolution and are
used for deriving lower bounds for the estimation problem as well as
upper bounds for the risk of the estimators. Specifically, these bounds
correspond to ordinary smooth error distributions (Fan \cite
{fan1991}), cf. Section~\ref{SecDisc} below for the super-smooth case.

Applying the plug-in approach, we need to integrate the density
estimator over an unbounded interval. As mentioned above, additional
assumptions are necessary to control $\|\widetilde f_b\|_{L^1}$. We
apply Lemma~\ref{lemFourierMult} assuming $\gamma\geq4$, that
is $\E[(\varepsilon_1^*)^4]<\infty$, and a polynomial decay of
$|\phi_\varepsilon'|$. The latter is a natural Mihlin-type condition
in the context of Fourier multipliers. Note that $\phi_\varepsilon'$
exists if $f_\varepsilon$, the distribution of the measurement errors,
has a first moment. In view of the analysis by Neumann and Rei{\ss}
\cite{neumann2009nonparametric}, the moment assumption in particular
implies uniform convergence of $\phi_{\varepsilon,m}$.
\end{remark}

To control the estimation error of $\widetilde q_{\tau,b}$, we follow
the Z-estimator approach (cf. van~der Vaart \cite{vanderVaart1998}).
Let $M(\eta)$ be the deterministic counterpart of $\widetilde M_b(\eta
)$ defined in (\ref{eqDefEst}). The quantities $\widetilde q_{\tau,b}$
and $q_\tau$ are given by the (approximated) zeros of $\widetilde
M_{b}$ and $M$, respectively. From the Taylor expansion
$0\approx\widetilde M_b(\widetilde q_{\tau,b})=\widetilde M_b(q_\tau
)+(\widetilde q_{\tau,b}-q_\tau)\widetilde M_b^\prime(q_\tau^*)$
for some intermediate point $q_\tau^*$ between $q_\tau$ and
$\widetilde q_{\tau,b}$, we obtain
%
\begin{equation}
\label{eqtaylor} \widetilde q_{\tau,b}-q_\tau\approx-
\frac{\int_{-\infty}^{q_\tau
}(\widetilde f_b(x)-f(x))\,\d x}{\widetilde f_b(q_\tau^*)}.
\end{equation}
The following two propositions deal separately with the numerator and
the denominator in this representation. The results are intrinsic to
our analysis, but may also be of interest on their own. The first
proposition deals with the numerator in (\ref{eqtaylor}) and
establishes minimax rates of convergence for estimation of the
distribution function with unknown error distributions.
Note that the quotient in (\ref{eqtaylor}) might explode if
$\widetilde f_b(q_\tau^*)$ becomes very small for large stochastic
error. Excluding this event which has vanishing probability, we
establish convergence rates as $\mathcal{O}_P$-results.

\begin{proposition}\label{prDistUnknown}
Suppose that Assumption~\textup{\ref{asKernel}} holds with $\ell=\langle
\alpha\rangle+1$ and let $b^*_{n,m}=(n\wedge m)^{-1/(2\alpha+2(\beta
\vee1/2)+1)}$. Then for any $\alpha\geq1/2$, $\beta, R,r,\zeta
>0$ and $\gamma\geq4$ we have uniformly over $f\in\mathcal
C^\alpha(R,r,\zeta)$ and $f_\varepsilon\in\mathcal D^\beta(R,
\gamma)$ as $n\wedge m\to\infty$,
\[
\biggl|\int_{-\infty}^{q_\tau}\bigl(\widetilde
f_{b^*_{n,m}}(x)-f(x)\bigr)\,\d x \biggr| =\mathcal{O}_P \bigl(
\psi_{n\wedge m}(\alpha,\beta) \bigr),
\]
where for $k\geq1$
%
\begin{equation}
\label{eqrate} \psi_k(\alpha,\beta):=\cases{ k^{-1/2},&\quad
for $\beta\in(0,1/2)$,
\cr
(\log k/k)^{1/2},&\quad for $\beta=1/2$,
\vspace*{3pt}\cr
k^{-(\alpha+1)/(2\alpha+2\beta+1)},&\quad for $\beta>1/2$.}
\end{equation}
\end{proposition}

Since the techniques to obtain Proposition~\ref{prDistUnknown} differ
significantly from previous results for deconvolution with unknown
error distribution, let us briefly sketch the proof: we apply a smooth
truncation function $a_s$ to decompose the error into
%
\begin{eqnarray}\label{eqsketchProof}
&& \int_{-\infty}^{q_\tau}\bigl(\widetilde
f_b(x)-f(x)\bigr)\,\d x\nonumber
\\
&&\quad = \underbrace{\int_{-\infty}^{q_\tau}
\bigl(K_b\ast f(x)-f(x)\bigr)\,\d x}_{\mathrm{deterministic\ error}}
+\underbrace{\int
_{-\infty}^{q_\tau}a_s(x+q_\tau)
\bigl(\widetilde f_b(x)-K_b\ast f(x)\bigr)\,\d
x}_{\mathrm{singular\ part\ of\ stochastic\
error}}
\\
&&\qquad{}+\underbrace{\int_{-\infty}^{q_\tau}
\bigl(1-a_s(x+q_\tau)\bigr) \bigl(\widetilde
f_b(x)-K_b\ast f(x)\bigr)\,\d x}_{\mathrm{continuous\ part\ of\
stochastic\ error}}\nonumber
\end{eqnarray}
with the usual notation $K_b(\cdot)=b^{-1}K(\cdot/b)$. The function
$a_s$ can be chosen such that it has compact support and satisfies
$(\ind_{(-\infty,0]}-a_s)\in C^\infty(\R)$. Similar to the
classical bias-variance trade-off, the deterministic error and singular
part of the stochastic error will determine the rate. The continuous
part, however, corresponds to the estimation error of a smooth (but not
integrable) functional of the density. If the error distribution were
known, it would be of order $n^{-1/2}$.
For unknown errors we use Lemma~\ref{lemFourierMult}, where our
estimate of the operator norm of the random Fourier multiplier $\F
^{-1}[\phi_K(bu)/\phi_{\varepsilon,m}(u) {\ind_{\{|\phi
_{\varepsilon,m}(u)|\geq m^{-1/2}\}}}]$ is of order $\mathcal{O}_P(1\vee(m^{-1/2}b^{-\beta-1}))$. This might be larger than the
operator norm of the unknown deconvolution operator $\F^{-1}[1/\phi
_\varepsilon(u)]$ which is uniformly bounded. Yet, for $\alpha
\geq1/2$ the additional error that appears in the continuous part
of stochastic error in (\ref{eqsketchProof}) is negligible.

Next, we like to understand the denominator of (\ref{eqtaylor}).
Lounici and Nickl \cite{louniciNickl2011} have proved uniform risk
bounds for the deconvolution wavelet estimator on the whole real line
for a known error distribution. On a bounded interval, which is
sufficient for our purpose, uniform convergence of the deconvolution
estimator $\widetilde f_b$ can be proved more elementarily. With
$b_{n}=(\log n/n)^{1/(2\alpha+2\beta+1)}$ the following proposition
yields the minimax rate $(\log n/n)^{\alpha/(2\alpha+2\beta+1)}$ in
$L^\infty$-loss (at least if $\frac{n}{\log n}\leq m$).
%

\begin{proposition}\label{pruniform}
Grant Assumption~\textup{\ref{asKernel}} with $\ell=\langle\alpha\rangle$.
For any $\alpha,\beta,R,r,\zeta>0$ and $\gamma\geq0$ we have
uniformly over $f\in\mathcal C^\alpha(R,r,\zeta)$ and $f_\varepsilon
\in\mathcal D^\beta(R, \gamma)$ as $n\wedge m\to\infty$,
\[
\sup_{x\in(-\zeta,\zeta)}\bigl|\widetilde f_b(x+q_\tau)-f(x+q_\tau
)\bigr|=\mathcal{O}_P \biggl(b^\alpha+ \biggl(\frac{\log n}{n}
\vee\frac
{1}{m} \biggr)^{1/2}b^{-\beta-1/2} \biggr).
\]
In particular, if $b=b_{n,m}\to0$ and $(\frac{n}{\log n}\wedge
m)b_{n,m}^{2\beta+1}\to\infty$ as $n\wedge m\to\infty$,
$\widetilde f_{b_{n,m}}$ is a uniformly consistent estimator.
\end{proposition}

The two propositions above are the building blocks for the first main
result of this paper announced in the following theorem. The constant
preceding the rate depends only on the class parameters $\alpha,\beta
,\gamma,R,r,\zeta$. The location parameter $U_n$ can grow
logarithmically to infinity as $n\to\infty$.

%
\begin{theorem}\label{thRateUnk}
Let $\alpha\geq1/2$, $\beta,R,r,\zeta>0$ and $\gamma\geq
4$ and grant Assumption~\textup{\ref{asKernel}} with $\ell=\langle\alpha
\rangle+1$. Let $\widetilde q_{\tau,b^*_{n,m}}$ be the quantile
estimator defined in (\ref{eqDefEst}) associated with
$b^*_{n,m}=(n\wedge m)^{-1/(2\alpha+2(\beta\vee1/2)+1)}$ and with
$U_n\to\infty, U_n=\mathcal{O}(\log n)$. Then we have uniformly over
$f\in\mathcal C^\alpha(R,r,\zeta, {U_n})$ and $f_\varepsilon\in
\mathcal D^\beta(R, \gamma)$ as $n\wedge m\to\infty$,
\[
|\widetilde q_{\tau,b^*_{n,m}}-q_\tau|=\mathcal{O}_P
\bigl(\psi_{
{n\wedge m}}(\alpha,\beta) \bigr),
\]
where $\psi_\bull(\alpha,\beta)$ is given in (\ref{eqrate}).
\end{theorem}

Using the methods of the proof of Theorem~\ref{thRateUnk} and an
additional application of Bernstein's concentration inequality,
convergence rates for the uniform loss can be obtained, assuming
regularity in a neighborhood of some interval of quantiles. For $0<\tau
_1<\tau_2<1$ and $\alpha,R,r,\zeta,U_n>0$, define
\begin{eqnarray*}
&& \mathcal C^\alpha_{\infty}(\tau_1,
\tau_2,R,r,\zeta,U_n)
\\
&&\quad := \Bigl\{ f\in\mathcal Q(R) \big|
\mbox{for all }\tau\in(\tau_1,\tau_2)\dvt f\mbox{ has a }\tau
\mbox{-quantile } {q_\tau\in[-U_n,U_n]}\mbox{ and}
\\
&&\hspace*{28pt} f\in C^\alpha\bigl([q_{\tau_1}-\zeta,q_{\tau_2}+
\zeta],R\bigr),\inf_{\tau
\in(\tau_1,\tau_2)}f(q_\tau)\geq r \Bigr
\}.
\end{eqnarray*}
%

\begin{theorem}\label{thuniformLoss}
Let $\alpha\geq1/2$, $\beta,R,r,\zeta>0$ and $\gamma\geq
4$ and grant Assumption~\textup{\ref{asKernel}} with $\ell=\langle\alpha
\rangle+1$. For $0<\tau_1<\tau_2<1$ and $\tau\in(\tau_1,\tau_2)$
let $\widetilde q_{\tau,b^*_{n,m}}$ be the quantile estimator defined
in (\ref{eqDefEst}) associated with $b^*_{n,m}=(\frac{\log
n}{n}\vee\frac{1}{m})^{1/(2\alpha+2(\beta\vee1/2)+1)}$ and with
$U_n\to\infty, U_n=\mathcal{O}(\log n)$. Then we have uniformly over
$f\in\mathcal C^\alpha_{\infty}(\tau_1,\tau_2,R,r,\zeta,U_n)$ and
$f_\varepsilon\in\mathcal D^\beta(R, \gamma)$ as $n\wedge m\to
\infty$,
\[
\sup_{\tau\in(\tau_1,\tau_2)}|\widetilde q_{\tau,b^*_{n,m}}-q_\tau|=
\mathcal{O}_P \bigl(\psi_{(n/\log n)\wedge
m}(\alpha,\beta) \bigr),
\]
where $\psi_\bull(\alpha,\beta)$ is given in (\ref{eqrate}).
\end{theorem}

We finish this subsection by providing the minimax rates for estimating
the distribution function and the quantiles for the case of known error
distributions, restricting to pointwise loss. As above, the estimators
are given by plugging in the classical density estimator
%
\begin{equation}
\label{eqfhat} \widehat f_b(x):= \F^{-1} \biggl[
\frac{\phi_{n}(u)\phi_K(bu)}{\phi
_\varepsilon(u)} \biggr](x),\qquad x\in\R.
\end{equation}

%
\begin{corollary}\label{corRateknown}
Let $\alpha, \beta,R,r,\zeta>0$ and $\gamma\geq0$ and suppose
that the error distribution is known and $f_\varepsilon\in\mathcal
D^\beta(R, \gamma)$. Let Assumption~\textup{\ref{asKernel}} hold with $\ell
=\langle\alpha\rangle+1$. Let $\widehat q_{\tau,b}$ be the quantile
estimator based on the density deconvolution estimator (\ref
{eqfhat}) associated with $b^*_n=n^{-1/(2\alpha+2(\beta\vee
1/2)+1)}$ {and $U_n\to\infty, U_n=\mathcal{O}(\log n)$}. Then we
obtain uniformly over $f\in\mathcal C^\alpha(R,r,\zeta, {U_n})$ as
$n\to\infty$,
\begin{eqnarray*}
\biggl|\int_{-\infty}^{q_\tau}\bigl(\widehat
f_{b^*_n}(x)-f(x)\bigr)\,\d x \biggr| &=&\mathcal{O}_P \bigl(
\psi_n(\alpha,\beta) \bigr),
\\
|\widehat q_{\tau,b^*_n}-q_\tau|&=&\mathcal{O}_P
\bigl(\psi_n(\alpha,\beta) \bigr),
\end{eqnarray*}
where $\psi_\bull(\alpha,\beta)$ is given (\ref{eqrate}).
\end{corollary}

Here, we do not estimate the deconvolution operator and thus there is
no additional error in terms of $m$. Consequently, we do not need a
moment assumption on the error distribution and the convergence rates
hold true for all $\alpha>0$.

\subsection{Lower bounds}
In view of the lower bounds stated by Fan \cite{fan1991}, in case
$n\leq m$ the rates in Proposition~\ref{prDistUnknown} are optimal.
Using the error representation (\ref{eqtaylor}), the result for
distribution function estimation carries over to quantile estimation.
Therefore, we focus on the case $m<n$.
To provide a clear proof of the lower bound, we allow for a more
general class of distributions of $X_j$, assuming only local
assumptions. Using point measures, the estimation error of $\phi
_{\varepsilon}$ does not profit from the decay of the characteristic
function of $X_j$. One could also consider the case of bounded
densities $f$ and choose alternatives in the proof whose Fourier
transforms decay arbitrarily slowly, but this would require far more
technical arguments.


We define for $\alpha, R,r,\zeta>0$ and some interval $I\subset\R$
\begin{eqnarray*}
\widetilde{\mathcal C}^{\alpha+1}(R,r,I)&:=& \Bigl\{F \mbox{ c.d.f.} \big| F
\mbox{ has on }I\mbox{ a Lebesgue density }f\in C^\alpha(I,R)\mbox{ and
}\inf_{x\in I}f(x)\geq r \Bigr\},
\\
\widetilde{\mathcal C}^{\alpha+1}(R,r,\zeta)&:=& \bigl\{F \mbox{ c.d.f.}
| F
\mbox{ has a }\tau\mbox{-quantile } q_\tau\in\R\mbox{ and } F\in
\widetilde{\mathcal C}^{\alpha+1}\bigl(R,r,[q_\tau-
\zeta,q_\tau+\zeta]\bigr) \bigr\}.
\end{eqnarray*}

%
\begin{theorem}\label{thLowerBound}
Suppose that $Y_1,\dots,Y_n$ and $\varepsilon^*_1,\dots,\varepsilon
^*_m$ are independent. Let $q\in\R$ and $\alpha, \beta,R,r,\zeta
>0,\gamma\geq0$. Then for any $C>0$ there is some $\delta>0$
such that
\begin{eqnarray*}
&& \inf_{\bar F_{n,m}}\sup_{F\in\widetilde{\mathcal C}^{\alpha
+1}(R,r,[q-\zeta,q+\zeta])}\sup _{f_\varepsilon\in\mathcal D^\beta
(R,\gamma)}P \bigl(\bigl|\bar F_{n,m}(q)-F(q)\bigr|
\\
&&\hspace*{157pt} >C(n\wedge
m)^{-(\alpha
+1)/(2\alpha+(2\beta)\vee1+1)} \bigr)\geq\delta,
\\
&& \inf_{\bar q_{\tau,n,m}}\sup_{F\in\widetilde{\mathcal C}^{\alpha
+1}(R,r,\zeta)}\sup
_{f_\varepsilon\in\mathcal D^\beta(R,\gamma
)}P \bigl(\bigl|\bar q_{\tau,n,m}-q_\tau\bigr|>C(n
\wedge m)^{-(\alpha
+1)/(2\alpha+(2\beta)\vee1+1)} \bigr)\geq\delta,
\end{eqnarray*}
where the infima are taken over all estimators $\bar F_{n,m}$ and $\bar
q_{\tau,n,m}$, respectively.
\end{theorem}

This lower bound implies that the rates in Proposition~\ref
{prDistUnknown} and Theorem~\ref{thRateUnk} are minimax optimal,
except for the case $\beta=1/2$ where they deviate by a logarithmic factor.

\subsection{Discussion and extension}\label{SecDisc}

The previous results show that estimating the distribution function by
integrating a density deconvolution estimator is a minimax optimal
procedure and under the local H\"older condition the rates are
determined by $n\wedge m$.
In that point our results differ completely from previous studies.
Assuming $\alpha$-Sobolev regularity of $f$, the RMSE of the kernel
density estimator by Neumann \cite{neumann1997effect} is of order
$\mathcal{O}(n^{-\alpha/(2\alpha+\beta+1)}+m^{-((\alpha/\beta
)\wedge1)})$. Since the error in estimating $\phi_\varepsilon$ is
reduced by the decay of the characteristic function $\phi$ of $X_j$,
the risk is of much smaller order in $m$. Assuming local regularity on
$f$ only, ${\mathcal F}f$ can decay arbitrarily slowly such that this
reduction effect may not occur. Note that assuming global Sobolev
regularity would improve also the convergence rate of the plug-in estimator.

Interestingly, the dependence on $n$ and $m$ is not completely
symmetric. As an intrinsic property of the uniform loss, the
convergence rates are typically by a logarithmic factor slower than for
pointwise loss. Yet, in Proposition~\ref{pruniform} and Theorem~\ref
{thuniformLoss} this payment for uniform convergence affects only the
estimation of $\phi$ and thus the rate is determined by $\frac{\log
n}{n}\vee\frac{1}{m}$.

Although the focus of this paper is on ordinary smooth error
distributions, a generalization to supersmooth errors is worth
mentioning. Let us sketch this case of exponentially decaying $\phi
_\varepsilon$. Supposing $\E[|\varepsilon_k^*|^4]<\infty$ and
$|\phi_\varepsilon(u)|^{-1}=\mathcal{O} (\mathrm{e}^{\gamma_0|u|^\beta
} )$ as well as $|\phi_\varepsilon'(u)|=\mathcal{O}
(\mathrm{e}^{-\gamma_1|u|^\beta} ), u\in\R$, for some $\beta>0$ and
$\gamma_0\geq\gamma_1>0$, we obtain analogously to Lemma~\ref
{lemFourierMult} for sufficiently small $c,\gamma>0$ and for the
bandwidth $b_m^*=c(\log m)^{-1/\beta}$
\[
\biggl\|\F^{-1} \biggl[\frac{\phi_K(b_m^*u)}{\phi_{\varepsilon,m}(u)} \biggr
]\ast\psi\biggr\|_{C^{s}(\R)}
\ind_{B_\varepsilon(b_m^*)} \leq\mathcal E_{b_m^*}\|\psi\|
_{C^{s}(\R)}
\qquad\mbox{where } \mathcal E_b=\mathcal{O}_P \bigl(1
\vee \mathrm{e}^{\gamma b^{-\beta}} \bigr)
\]
for any $s\geq0$ and for any $\psi\in C^{s}(\R)$. In other
words, $\phi_K(bu)/\phi_{\varepsilon,m}(u)$ is a random Fourier
multiplier on H\"older spaces with exponentially increasing operator
norm {on the event $B_\varepsilon(b)$}. Following the lines of the
proof of Proposition~\ref{prDistUnknown}, one sees that the singular
as well as the continuous part of the stochastic error in (\ref
{eqsketchProof}) are of the order $\mathcal{O}_P((n\wedge
m)^{-1/2}\mathrm{e}^{\gamma b^{-\beta}})$. Combined with the estimate for the
deterministic error, the choice $b_{n,m}^*=c(\log(n\wedge
m))^{-1/\beta}$ yields for $f\in\mathcal C^\alpha(R,r,\zeta)$
\[
\biggl|\int_{-\infty}^{q_\tau}\bigl(\widetilde
f_{b^*_n}(x)-f(x)\bigr)\,\d x \biggr| =\mathcal{O}_P \bigl(\bigl(
\log(n\wedge m)\bigr)^{-(\alpha+1)/\beta} \bigr).
\]
Note that for $n\leq m$ this is the minimax rate for distribution
function estimation as given in Fan~\cite{fan1991}. Therefore, also
for supersmooth error distributions the integral domain does not need
to be truncated to estimate the distribution function via the plug-in approach.

%
\section{Adaptive estimation}\label{secadaptive}
The choice of the bandwidth $b$ is crucial in applications. Therefore,
we develop a fully data-driven procedure to determine a good bandwidth.
We follow the approach initiated by Lepski{\u\i} \cite{Lepski1990}.
More precisely, we use the version proposed in Goldenshluger and
Nemirovski \cite{goldenshluger1997spatially}. For simplicity, we
suppose $n=m$ and focus on the pointwise loss in this section.

Let us consider the family of estimators $\{\widetilde q_{\tau,b},b\in
\B_n\}$ where $\widetilde q_{\tau,b}$ is defined in (\ref
{eqDefEst}) and $\B_n$ is a finite set of bandwidths. In view of the
error representation~(\ref{eqtaylor}), it is important that
$\widetilde f_b(\widetilde q_{\tau,b})$ is a consistent estimator of
$f(q_\tau)$ for all $b\in\mathcal B_n$. Therefore, conditions on the
bandwidth as in Proposition~\ref{pruniform} are necessary for the
entire set $\mathcal B_n$. These depend on the true but unknown degree
of ill-posedness $\beta$ and on $\alpha$. We keep to the assumption
$\alpha>1/2$ such that the additional error due to bounding the random
Fourier multiplier is negligible. Note that the lower bound for the
bandwidth is not determined by the variance of the quantile estimator
itself but by the variance of the density estimator and the minimal
smoothing which results from $\alpha>1/2$.

Inspired by Comte and Lacour \cite{comte2011data}, we propose the
following construction of a feasible set $\B_n$: for some $L>1$ define
\[
b_{n,j}:=n^{-1}L^j\qquad\mbox{for }j=0,
\dots,N_n\mbox{ where } N_n\in\N\mbox{ satisfies
}n^{-1}L^{N_n}\sim(\log n)^{-3}.
\]
Choosing
%
\begin{equation}
\label{eqtildeJ} \widetilde j_{n}:=N_n\wedge\min\biggl
\{j=0,\dots,N_n-1\dvt \frac
{1}{2}\leq\biggl(
\frac{\log n}{n} \biggr)^{1/2}\int_{-1/b_{n,j}}^{1/b_{n,j}}
{\frac{\ind_{\{|\phi_{\varepsilon,m}(u)|\geq m^{-1/2}\}}}{|\phi
_{\varepsilon,m}(u)|}}\,\d u\leq1 \biggr\},\quad
\end{equation}
the bandwidth set is given by
%
\begin{equation}
\label{eqbandwidthset} \mathcal B_n:=\{b_{n,\widetilde j_{n}},
\dots,b_{n,N_n}\}.
\end{equation}
Note that by construction $\B_n$ is nonempty and it consists of a
monotone increasing sequence of bandwidths such that
$b_{n,j+1}/b_{n,j}$ is uniformly bounded in $j=\widetilde j_n,\dots
,N_n$ and $n\geq1$. Also, for $n\to\infty$ we have $N_n\lesssim
\log n$ and $(\log n)^2b_{n,N_n}\to0$. The following lemma establishes
two additional properties. The latter one ensures that for any $b\in
\mathcal B_n$ our estimators are consistent.

\begin{lemma}\label{lemBmin}
Let $(Y_j)$ and $(\varepsilon_k^*)$ be distributed according to $f\in
\mathcal C^\alpha(R,r,\zeta)$ and $f_\varepsilon\in\mathcal D^\beta
(R,\gamma)$ with $\alpha\geq1/2,\beta>0$. Then with probability
converging to one, $\widetilde j_n<N_n$ and the optimal bandwidth
$b^*_n=n^{-1/(2\alpha+2(\beta\vee1/2)+1)}$ is contained in the
interval $[b_{n,\widetilde j_n},b_{n,N_n}]$ as well as
\mbox{$nb_{n,\widetilde j_n}^{2\beta+2}\to\infty$}.
\end{lemma}

Given the bandwidth set, the adaptive estimator is obtained by
selection from the family of estimators $\{\widetilde q_{\tau,b},b\in
\B_n\}$. As proposed by Lepski{\u\i} \cite{Lepski1990} the adaptive
choice should mimic the trade-off between deterministic error and
stochastic error. The adaptive choice will be given by the largest
bandwidth such that the intersection of all confidence sets, which
corresponds to smaller bandwidths, is nonempty. As discussed above, it
is sufficient to consider the singular part of the stochastic error in
(\ref{eqsketchProof}) only. To estimate the variance of $\widetilde
q_{\tau, b}$ corresponding to the latter, we define for some $\delta>0$
%
\begin{eqnarray}
\label{eqTildeSigma} \hspace*{-25pt}\widetilde\Sigma_b&:=&\frac{(2\sqrt{2}+\delta)\sqrt
{\log\log
n}\max_{\mu\geq b}\widetilde\sigma_{\mu,X}+(\delta\log
n)^3\max_{\mu\geq b}\widetilde\sigma_{\mu,\varepsilon
}+(1+\delta)|\widetilde M_b(\widetilde q_{\tau,b})|}{|\widetilde
{f}_{b}(\widetilde q_{\tau,b})|},
\end{eqnarray}
with the truncation function $a_s$ from decomposition (\ref
{eqsketchProof}) and
%
\begin{eqnarray}
\label{eqTildeSigmaX} \widetilde\sigma^2_{b,X}&=&\frac{1}{n^2}
\sum_{j=1}^n \biggl(\int
_{-\infty}^0a_s(x) \F^{-1}
\biggl[\frac{\phi_K(bu)\mathrm{e}^{\mathrm{i}uY_j}}{\phi_{\varepsilon,m}(u)} \biggr
](x+\widetilde q_{\tau,b})\,\d x
\biggr)^2\quad\mbox{and}
\\
\label{eqTildeSigmaEPS} \widetilde\sigma^2_{b,\varepsilon}&=&
\frac{1}{4\pi^2m}\int_{-1/b}^{1/b}\bigl|
\phi_K(bu)\bigr| \biggl|\frac{\phi_n(u)}{\phi_{\varepsilon,m}(u)} \biggr|^2\,\d u\int
_{-1/b}^{1/b}\bigl|\phi_K(bu)\bigr| \biggl|
\frac{\F
a_s(u)}{\phi_{\varepsilon,m}(u)} \biggr|^2\,\d u.
\end{eqnarray}
The parameter $\delta$ has minor influence and should be chosen close
to zero. Note that we apply a monotonization in the numerator of
$\widetilde\Sigma_b$ by taking maxima of $\widetilde\sigma_{\mu,X}$ and
$\widetilde\sigma_{\mu,\varepsilon}$, respectively. The
correction term $|\widetilde M_b(\widetilde q_{\tau,b})|$ appears only
if $\widetilde q_{\tau,b}$ is not the exact solution of the estimating
equation~(\ref{eqestEq}). Define for any $b\in\B_n$
\begin{eqnarray}
\nonumber
\mathcal U_b:=[\widetilde q_{\tau,b}-
\widetilde{\Sigma}_b, \widetilde q_{\tau,b}+\widetilde{
\Sigma}_b].
\end{eqnarray}
The adaptive estimator is given by
%
\begin{eqnarray}
\label{eqadaptiveest} \widetilde q_\tau:=\widetilde q_{\tau,\widetilde
{b}^*_{n}}\qquad
\mbox{with } \widetilde{b}^*_{n}:=\max\biggl\{b\in\B_n \Big|
\bigcap_{\mu\leq
b,\mu\in\B_n} \mathcal U_\mu\neq\varnothing
\biggr\}.
\end{eqnarray}
Note that $\widetilde{b}^*_{n}$ is well defined since the intersection
in (\ref{eqadaptiveest}) is nonempty for $b=b_{n,1}$. The following
theorem shows that this estimator achieves the minimax rate up to a
logarithmic factor. The proof relies on a comparison with an
oracle-type choice of the bandwidth. All ingredients, though, have to
be estimated and the dependence between $Y_j$ and $\varepsilon^*_k$
requires special attention.

%
\begin{theorem}\label{thoracleinequality}
Let $n= m$ and $\alpha\geq1/2$, $\beta,R,r,\zeta>0,\gamma\geq
4$ and grant Assumption~\textup{\ref{asKernel}} with $\ell=\langle\alpha
\rangle+1$. Then the estimator $\widetilde q_\tau$ as defined in
(\ref{eqadaptiveest}) with $\B_n$ from (\ref{eqbandwidthset})
satisfies uniformly over $f\in\mathcal C^\alpha(R,r,\zeta, {U_n})$
and $f_\varepsilon\in\mathcal D^\beta(R, \gamma)$ as $n\to\infty$,
\begin{eqnarray*}
|\widetilde q_\tau-q_\tau|=\mathcal{O}_P
\bigl( {\psi_{n(\delta\log
n)^{-6}}(\alpha,\beta)} \bigr),
\end{eqnarray*}
where $\psi_\bull(\alpha,\beta)$ is given in (\ref{eqrate}).
\end{theorem}

As the theorem shows, the adaptive method achieves the minimax rate up
to a logarithmic factor. This additional loss is dominated by the
stochastic error which is due to the estimation of $\phi_\varepsilon
$. Since $Y_j$ and $\varepsilon^*_k$ are not independent, we have to
bound the stochastic error of $\widetilde q_{\tau,b}$ in a way that
separates the error terms coming from the estimation of $\phi$ and
$\phi_\varepsilon$, respectively. Estimating the remaining parts, we
lose the factor $(\delta\log n)^6$, which appears not to be optimal.
To improve the rate slightly, $\delta=\delta_n$ could be chosen as a
null sequence provided $\delta_n(\log n)^{1/2}\to\infty$. In the
case where the error density is known, we can achieve the better rate
$\psi_{n/\log\log n}(\alpha,\beta)$. The $\log\log n$-factor is
the additional payment for $\mathcal{O}_P$-adaptivity, which is known to
be unavoidable for a bounded loss function in standard regression, cf.
Spokoiny \cite{spokoiny1996}. For estimating the distribution
function, an analogous result can be obtained, but is omitted.

\section{Numerical results}\label{secnumeric}
\subsection{Simulation study}
We illustrate the implementation of the adaptive estimation procedure of
Section~\ref{secadaptive}. Our small simulation study serves as a
proof of viability of the proposed method.

We run $1000$ Monte Carlo simulations for four experimental setups. The
sample size is set to $n=1000$ and the external sample of the directly
observed error is set to $m=1000$ as well (here the external sample is
independent of the main one). We consider $\Gamma(1,1)$ and $\Gamma
(2,1)$ for the distribution of $X$ where $\Gamma(k,\eta)$ denotes the
gamma distribution with shape parameter $k$ and scale~$\eta$. 
Note that the shape $k$ of the gamma distribution determines the
Sobolev smoothness of the density while the density is smooth away from
the origin. For the error distribution, we consider $\Gamma(1,\sqrt
{2})$ centered around zero which corresponds to $\beta=1$ and the
standard Laplace distribution (scale equals $1$) corresponding to
$\beta=2$. In both cases, the variance of the error equals $2$.

%
\begin{table}
\tabcolsep=0pt
\caption{Empirical root mean square error (RMSE) of the adaptive
deconvolution estimator and the empirical quantiles of $(Y_j)$ (in
parenthesis) for estimating $q_\tau$ based on $1000$ Monte Carlo
simulations with $n=m=1000$} \label{tabsim1}
\begin{tabular*}{\tablewidth}{@{\extracolsep{\fill}}@{}lllll@{}}
\hline
RMSE &$k=1$, $\beta=1$&$k=2$, $\beta=1$&$k=1$, $\beta=2$&$k=2$, $\beta=2$\\
\hline
$\tau=0.1$ & 0.532 (0.886) & 0.252 (0.706) & 0.378 (1.029) & 0.191 (0.765)\\
$\tau=0.2$ & 0.265 (0.653) & 0.114 (0.508) & 0.175 (0.452) & 0.091 (0.349)\\
$\tau=0.3$ & 0.111 (0.461) & 0.070 (0.360) & 0.077 (0.178) & 0.090 (0.158)\\
$\tau=0.4$ & 0.067 (0.282) & 0.080 (0.212) & 0.112 (0.052) & 0.105 (0.064)\\
$\tau=0.5$ & 0.123 (0.110) & 0.092 (0.096) & 0.171 (0.175) & 0.116 (0.145)\\
$\tau=0.6$ & 0.162 (0.122) & 0.094 (0.123) & 0.200 (0.318) & 0.109 (0.255)\\
$\tau=0.7$ & 0.154 (0.326) & 0.098 (0.272) & 0.189 (0.462) & 0.098 (0.373)\\
$\tau=0.8$ & 0.107 (0.597) & 0.150 (0.481) & 0.115 (0.624) & 0.141 (0.506)\\
$\tau=0.9$ & 0.232 (1.015) & 0.312 (0.783) & 0.226 (0.849) & 0.293 (0.675)\\
\hline
\end{tabular*}
\end{table}

The target quantiles of interest are $q_\tau$ with $\tau
=0.1,0.2,\dots,0.9$. 
In the real data example in the next subsection we compare the adaptive
estimator to the ``naive'' quantile estimator given by the sample
quantiles of the observations $Y$. Therefore we have also applied the
naive estimator in the simulations. The results of this simulation
study are given in Table~\ref{tabsim1}. We can see that the results
support the theory -- the empirical root mean squared error (RMSE) is
higher in most cases for $\beta=2$ than for $\beta=1$. Also, we can
see that in most cases the RMSE is lower for $k=2$ than for $k=1$ since
the gamma distribution with larger shape parameter is smoother in our
context. At the tails, our estimation method is significantly better
than the naive estimator. Near the median the naive estimator behaves
nice when the distribution of the error is Laplace. This is not the
case under the gamma error distribution which may suggest that the
naive estimator profits from the symmetry of the error distribution.
Similar behavior was observed also in distribution deconvolution with
nonsymmetric error distributions, see Dattner and Reiser \cite
{dattner2012estimation}. 

%
\subsection{Real data example}

High blood pressure is a direct cause of serious cardiovascular disease
(Kannel \textit{et~al}. \cite{kannel1995framingham}) and determining reference
values for physicians is important. In particular, estimating
percentiles of systolic and diastolic blood pressure by sex, race or
ethnicity, age, etc. is of substantial interest. Blood pressure is
known to be measured with additional error which needs to be addressed
in its analysis (see e.g., Frese, Fick and Sadowsky \cite
{frese2011blood}). Therefore, measurement errors should be taken into
account, otherwise quantile estimates based on the observed blood
pressure measurements would be biased.

We illustrate our method using data from the Framingham Heart Study
(Carroll \textit{et~al}. \cite{carroll2006measurement}). This study
consists of
a series of exams taken two years apart where systolic blood pressure
(SBP) measurements of 1615 men aged 31--65 were taken. These data
were used as an illustration for density deconvolution by Stirnemann,
Comte and Samson \cite{stirnemann2012density} and for distribution
deconvolution by Dattner and Reiser \cite{dattner2012estimation}. We
denote by $Y_{j,1}$ and $Y_{j,2}$
the two repeated measures of SBP for each individual $j$ at two
different exams and denote by $X_j$ the long-term average SBP of
individual~$j$. Then we model that
\begin{eqnarray*}
Y_{j,1}=X_j+\varepsilon_{j,1},\qquad
Y_{j,2}=X_j+\varepsilon_{j,2},
\end{eqnarray*}
for individuals $j=1,\ldots,n$. Following Carroll \textit{et~al}. \cite
{carroll2006measurement}, we use the average of the two exams
$Y^{\prime}_j=(Y_{j,1}+Y_{j,2})/2$, so that the model in our case is
\begin{eqnarray}
\nonumber
Y^{\prime}_j=X_j+
\varepsilon^{\prime}_j,
\end{eqnarray}
where $\varepsilon^{\prime}_j=(\varepsilon_{j,1}+\varepsilon_{j,2})/2$.

%
\begin{figure}[t]

\includegraphics{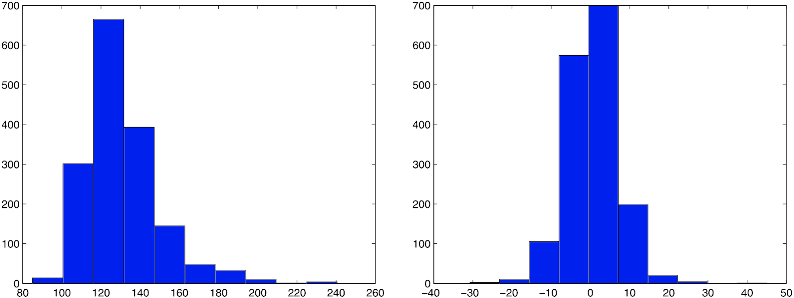}

\caption{Average systolic blood pressure $Y^\prime$ (\emph{left})
and the errors $\varepsilon^*$ (\emph{right}) over the two
measurements from the two visits of 1615 men aged 31--65 from the
Framingham Heart Study.}\label{figY}
\end{figure}

Taking advantage of the repeated measurements, we can avoid parametric
assumptions regarding the distribution of the errors. The only
assumption we will make is that the distribution of the measurement
error is symmetric around zero and does not vanish. We then set
$\varepsilon^*_j=(Y_{j,1}-Y_{j,2})/2$ and note that under the symmetry
assumption it is distributed as $\varepsilon^{\prime}_j$. We
emphasize the fact that our theoretical results do not require that the
sample $\varepsilon^*_j$ must be independent from that of the
$Y^{\prime}_j$.

%
\begin{figure}

\includegraphics{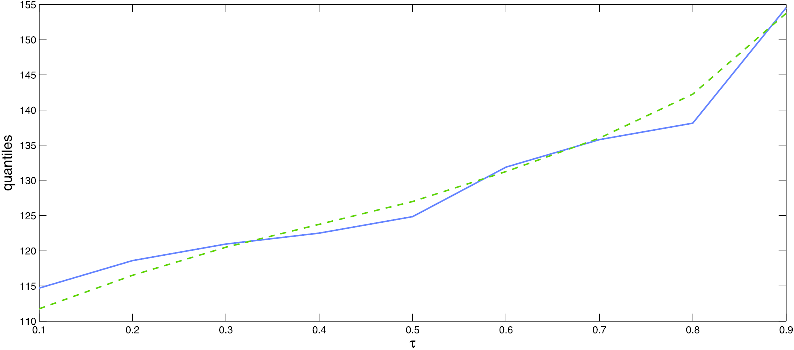}

\caption{Quantiles estimates for systolic blood pressure of 1615
men aged 31--65 from the Framingham Heart Study. Solid line for the
adaptive deconvolution estimator and dashed line for the empirical quantiles of $(Y_j)$.}\label{figqest}
\end{figure}

Histograms of $Y^\prime$ and $\varepsilon^*$ are presented in
Figure~\ref{figY}. Although Figure~\ref{figY} may suggest that the
error distribution does not entirely satisfy the symmetry assumption,
it serves as working hypothesis for our procedure and, indeed, it is
supposed in previous works on the same data set as well. The resulting
adaptive and naive quantiles estimates are displayed in Figure~\ref
{figqest}. 
We can see certain differences between the naive and adaptive estimates
which might result in important implications for medical research, but
here we do not aim at pursuing a more detailed statistical analysis.

%
\section{Proofs}\label{secproofs}
\subsection{\texorpdfstring{Proofs for Section~\protect\ref{secminimax}}{Proofs for Section 2}}
For convenience, we will write $A_n(\theta)\lesssim B_n(\theta)$ if
$A_n(\theta)=\mathcal{O}(B_n(\theta))$. For a better readability, we
assume throughout $\beta\neq1/2$. In the special case, $\beta=1/2$
the order of the stochastic error will be $(\log n/n)^{1/2}$ which can
be easily seen below in the bounds (\ref{eqETsX}) and (\ref
{eqETsEps}). The subscript $n$ at the bandwidth will be omitted.

Since $1/\phi_{\varepsilon,m}$ might explode for large stochastic
errors we need the following lemma.

\begin{lemma}\label{lemPhiEps}
Suppose $\E[|\varepsilon_k^*|^{\delta}]<\infty$ for some $\delta
>0$. Let $T_m\to\infty$ be an increasing sequence satisfying
$m^{1/2}\inf_{u\in[-T_m,T_m]}|\phi_\varepsilon(u)|\gtrsim(\log
T_m)^2$, then for any $p<2$
\[
P \Bigl(\inf_{u\in[-T_m,T_m]}\bigl|\phi_{\varepsilon,m}(u)\bigr|<m^{-1/2}(
\log T_m)^p \Bigr)=\mathrm{o}(1)\qquad\mbox{as } m\to\infty.
\]
\end{lemma}

\begin{pf}
The triangle inequality, the assumption on $T_m$ and Markov's
inequality yield for $m$ as well as $T_m$ large enough
\begin{eqnarray*}
&& P \Bigl(\inf_{u\in[-T_m,T_m]}\bigl|\phi_{\varepsilon,m}(u)\bigr|<m^{-1/2}(
\log T_m)^p \Bigr)
\\
&&\quad \leq P \Bigl(\sup_{u\in[-T_m,T_m]}\bigl|
\phi_\varepsilon(u)-\phi_{\varepsilon,m}(u)\bigr|> \inf_{u\in[-T_m,T_m]}\bigl|
\phi_\varepsilon(u)\bigr|-m^{-1/2}(\log T_m)^p
\Bigr)
\\
&&\quad \lesssim \frac{2}{(\log T_m)^2}\E\Bigl[\sup
_{u\in
[-T_m,T_m]}m^{1/2}\bigl|\phi_\varepsilon(u)-
\phi_{\varepsilon,m}(u)\bigr| \Bigr].
\end{eqnarray*}
Noting $\ind_{[-T_m,T_m]}(u)\leq w(u)/w(T_m)$ for $w(u):=(\log
(\mathrm{e}+|u|))^{-1/2-\eta}$ for some $\eta\in(0,1/2)$, the above display
can be bounded by
%
\begin{equation}
\label{eqboundNeumannReiss} \frac{2}{w(T_m)(\log T_m)^2}\E\Bigl[\sup
_{u\in\R}m^{1/2}w(u)
\bigl|\phi_\varepsilon(u)-\phi_{\varepsilon,m}(u) \bigr| \Bigr] \lesssim(\log
T_m)^{-3/2+\eta},
\end{equation}
where the expectation is bounded by applying Theorem 4.1 in Neumann and
Rei{\ss} \cite{neumann2009nonparametric}.
\end{pf}

To ensure consistency of the density estimator, the bandwidth satisfies
usually $(n\wedge m)b^{2\beta+1}\to\infty$ and is of polynomial
order in $n,m$. This implies $m^{1/2}\inf_{u\in[-1/b,1/b]}|\phi
_\varepsilon(u)|\gtrsim|\log b|^2$ for $f\in\mathcal D^\beta
(R,\gamma), \gamma>0$, and thus Lemma~\ref{lemPhiEps} can be
applied to $T_m=1/b$. Under this conditions on~$b$ the probability of
the event $B_\varepsilon(b)$, defined in (\ref{eqBEps}), tends to
one. In that case, it suffices to control terms on $B_\varepsilon(b)$,
a strategy that will follow in the sequel. For instance, the $\mathcal{O}_P$-convergence in Theorem~\ref{thRateUnk} is equivalent to $\lim_{C\to
\infty}\lim_{n,m\to\infty}P(|\widetilde q_{\tau,b^*}-q_\tau
|>C\psi_{n\wedge m}(\alpha,\beta))=0$ for which we have
\begin{eqnarray*}
&&\lim_{C\to\infty}\lim_{n,m\to\infty}P \bigl(|\widetilde
q_{\tau,b^*}-q_\tau|>C\psi_{n\wedge m}(\alpha,\beta) \bigr)
\\
&&\quad\leq\lim_{C\to\infty}\lim_{n,m\to\infty}P
\bigl(\bigl\{ |\widetilde q_{\tau,b^*}-q_\tau|>C
\psi_{n\wedge m}(\alpha,\beta)\bigr\} \cap B_\varepsilon\bigl(b^*\bigr)
\bigr)+\lim_{m\to\infty}P\bigl(B_\varepsilon\bigl(b^*
\bigr)^c\bigr),
\end{eqnarray*}
where the second term converges to zero by Lemma~\ref{lemPhiEps} and
it remains to bound the first one.

On $B_\varepsilon(b)$ the weaker estimate $|\phi_{\varepsilon
,m}(u)|\geq m^{-1/2}$ for $|u|\leq1/b$ will frequently be
enough, implying
\[
\frac{\phi_K(bu)}{\phi_{\varepsilon,m}(u)}\ind_{\{|\phi
_{\varepsilon,m}(u)|\geq m^{-1/2}\}}=\frac{\phi_K(bu)}{\phi
_{\varepsilon,m}(u)}\qquad\mbox{on
}B_\varepsilon(b).
\]

\subsubsection{\texorpdfstring{Proof of Lemma~\protect\ref{lemWellDef}}{Proof of Lemma 2.1}}
On $B_\varepsilon(b)$, we have by continuity of the characteristic
functions and the properties of the kernel that $g(u):=\frac{\phi
_n(u)\phi_K(bu)}{\phi_{\varepsilon,m}(u)}$ satisfies $g,g'\in L^2(\R
)$. Hence, $(1+x^2)^{1/2}\F^{-1}g(x)\in L^2(R)$ and the
Cauchy--Schwarz inequality yields $\|\widetilde f_b\|_{L^1}\leq\|
(1+x^2)^{-1/2}\|_{L^2}\|(1+x^2)^{1/2}\F^{-1}g(x)\|_{L^2}<\infty$ on
$B_\varepsilon(b)$. In particular, (\ref{eqestEq}) is well defined
on the event $B_\varepsilon(b)$.

On $B_\varepsilon(b)$, we have moreover $\lim_{\eta\to-\infty}\int
_{-\infty}^\eta\widetilde f_b(x)\,\d x=0$, by integrability of
$\widetilde f_b$, and $\int_{-\infty}^\infty\widetilde f_b(x)\,\d x=\F
[\widetilde f_b](0)=\phi_n(0)\phi_K(0)/\phi_{\varepsilon,m}(0)=1$.
Applying $\|\widetilde f_b\|_\infty\leq\|\phi_K(bu)/\break \phi
_{\varepsilon_m}(u)\|_{L^1}<\infty$, we conclude that $\eta\mapsto
\int_{-\infty}^\eta\widetilde f_b(x)\,\d x$ continuous and $[0,1]$ is
contained in its range. 

\subsubsection{\texorpdfstring{Proof of Lemma~\protect\ref{lemFourierMult}}{Proof of Lemma 2.2}}
Note that the assumption on $\phi_\varepsilon$ imply $|(\phi
_\varepsilon^{-1})'(u)|\lesssim(1+|u|)^{\beta-1}$ as well as $|\phi
_\varepsilon^{-1}(u)|\lesssim(1+|u|)^{\beta}, u\in\R$. We define
the random Fourier multiplier
\[
\psi(u):=(1+\mathrm{i}u)^{-\beta}\frac{\phi_K(bu)}{\phi_{\varepsilon
,m}(u)},\qquad u\in\R.
\]
On $B_\varepsilon(b)$, as defined in (\ref{eqBEps}), we will check
H\"ormander type conditions and derive an upper bound for the operator
norm of $\psi(u)$. Hence, we have to determine a suitable constant
$A_\psi>0$ satisfying
%
\begin{eqnarray}\label{eqFMProp}
\max_{l\in\{0,1\}} \biggl(
\int_{[-2,2]}\bigl|\psi^{(l)}(u)\bigr|^2\,\d u
\biggr)^{1/2}&\leq& A_\psi\quad\mbox{and}
\nonumber\\[-8pt]\\[-8pt]
\max_{l\in\{0,1\}}\sup_{T\in[1,\infty)}T^{l-1/2}
\biggl(\int_{T\leq|u|\leq4T}\bigl|\psi^{(l)}(u)\bigr|^2\,\d u
\biggr)^{1/2}&\leq& A_\psi.\nonumber
\end{eqnarray}
To find $A_\psi$, we note that
%
\begin{eqnarray}
\frac{1}{|\phi_{\varepsilon,m}(u)|^p} &\leq&\frac{p}{|\phi
_\varepsilon(u)|^p}+\frac{p|\phi
_{\varepsilon,m}(u)-\phi_\varepsilon(u)|^p}{|\phi_\varepsilon
(u)\phi_{\varepsilon,m}(u)|^p},\qquad
\mbox{for }p\in\{1,2\}\label
{eqPhiInv}
\end{eqnarray}
and thus on $B_\varepsilon(b)$
\begin{eqnarray*}
\frac{1}{|\phi_{\varepsilon,m}(u)|}\leq\frac{1+\Delta
_m(u)}{|\phi_\varepsilon(u)|},\qquad\Delta_m(u):=
\frac
{m^{1/2}}{|\log b|^{3/2}}\bigl|\phi_{\varepsilon,m}(u)-\phi_\varepsilon(u)\bigr|.
\end{eqnarray*}
By the assumptions on $\phi_\varepsilon$ and $K$, we conclude
%
\begin{equation}
\label{eqPsiBound} \bigl|\psi(u)\bigr|\leq\frac{|\phi_K(bu)|(1+\Delta
_m(u))}{(1+u^2)^{\beta/2}|\phi_\varepsilon(u)|}\lesssim\bigl(1+\Delta
_m(u)\bigr)\ind_{[-1/b,1/b]}(u).
\end{equation}
Concerning the derivative, we estimate $b\leq2(1+|u|)^{-1}$ for
$|u|\leq1/b$ and $b<1/2$ and consequently by $|\phi_\varepsilon
'(u)/\phi_\varepsilon(u)|\lesssim(1+|u|)^{-1}$
%
\begin{eqnarray}\label{eqPsiPBound}
\bigl|\psi'(u)\bigr|&\leq&(\beta+1) \bigl(1+u^2
\bigr)^{-(\beta+1)/2} \biggl|\frac
{\phi_K(bu)}{\phi_{\varepsilon,m}(u)} \biggr| +b\bigl(1+u^2
\bigr)^{-\beta/2} \biggl|\frac{\phi_K'(bu)}{\phi_{\varepsilon,m}(u)} \biggr|
\nonumber
\\
&&{}+\bigl(1+u^2\bigr)^{-\beta/2} \biggl|\frac{\phi_{\varepsilon
,m}'(u)}{\phi_{\varepsilon,m}(u)}
\frac{\phi_K(bu)}{\phi
_{\varepsilon,m}(u)} \biggr|
\nonumber
\\
&\lesssim&\frac{|\psi(u)|}{1+|u|}+\bigl|\psi(u)\bigr| \biggl|\frac{\phi
_{\varepsilon,m}'(u)}{\phi_{\varepsilon,m}(u)} \biggr|
\nonumber
\\
&\lesssim&\bigl(1+\Delta_m(u)\bigr) \biggl(\frac{1}{1+|u|}+
\bigl(1+\Delta_m(u)\bigr) \biggl|\frac{\phi_{\varepsilon,m}'(u)}{\phi
_{\varepsilon}(u)} \biggr| \biggr)
\ind_{[-1/b,1/b]}(u)
\\
&\lesssim&\bigl(1+\Delta_m(u)\bigr) \biggl(\frac{2+\Delta_m(u)}{1+|u|}+
\bigl(1+\Delta_m(u)\bigr) \bigl(1+|u|\bigr)^{\beta}\bigl|
\phi_{\varepsilon,m}'(u)-\phi_\varepsilon'(u)\bigr|
\biggr)\nonumber
\\
&&{}\times \ind_{[-1/b,1/b]}(u)\nonumber
\\
&\lesssim&\frac{(1+\Delta_m(u))^2}{1+|u|} \bigl(1+\bigl(1+|u|\bigr)^{\beta
+1}\bigl|\phi_{\varepsilon,m}'(u)-
\phi_\varepsilon'(u)\bigr| \bigr)\ind_{[-1/b,1/b]}(u).\nonumber
\end{eqnarray}
With these bounds at hand, we can show (\ref{eqFMProp}). For $l=0$,
the estimate (\ref{eqPsiBound}) and $1/T\lesssim(1+|u|)^{-1}$ for
$|u|\leq4T$ yield
\begin{eqnarray*}
\int_{-2}^2\bigl|\psi(u)\bigr|^2\,\d u &\lesssim&\int_{-2}^2\bigl(1+\Delta^2_m(u)
\bigr)\ind_{[-1/b,1/b]}(u)\,\d u,
\\
\frac{1}{T}\int_{T\leq|u|\leq4T}\bigl|\psi(u)\bigr|^2\,\d u &\lesssim&\frac{1}{T}\int_{T\leq|u|\leq4T}\bigl(1+\Delta
^2_m(u)\bigr)\ind_{[-1/b,1/b]}(u)\,\d u
\\
&\lesssim&1+\int_{-1/b}^{1/b}\bigl(1+|u|\bigr)^{-1}
\Delta^2_m(u)\,\d u,
\end{eqnarray*}
for $b$ small enough. Hence, the conditions (\ref{eqFMProp}) for
$l=0$ are satisfied for $A_\psi$ of the order
$
(1+\int_{-1/b}^{1/b}(1+|u|)^{-1}\Delta^2_m(u)\,\d u)^{1/2}$.
For $l=1$, we verify by (\ref{eqPsiPBound}) and $T\leq(1+|u|)$
for $|u|>T$
\begin{eqnarray*}
&& \int_{-2}^2\bigl|\psi'(u)\bigr|^2
\,\d u
\lesssim \int_{-2}^2 \bigl(1+
\Delta_m^4(u) \bigr) \bigl(1+\bigl(1+|u|\bigr)^{2\beta+2}\bigl|
\phi_{\varepsilon,m}'(u)-\phi_\varepsilon'(u)\bigr|^2
\bigr)\,\d u\quad\mbox{and}
\\
&& T\int_{T\leq|u|\leq4T}\bigl|\psi'(u)\bigr|^2\,\d u
\\
&&\quad
\lesssim\int_{T\leq|u|\leq4T}\frac{T\,\d
u}{(1+|u|)^{2}}
\\
&&\qquad{} +\int
_{-1/b}^{1/b} \biggl(\frac{\Delta
_m^4(u)}{1+|u|}+\bigl(1+
\Delta_m^4(u)\bigr) \bigl(1+|u|\bigr)^{2\beta+1}\bigl|\phi
_{\varepsilon,m}'(u)-\phi_\varepsilon'(u)\bigr|^2
\biggr)\,\d u
\\
&&\quad \lesssim 1+\int_{-1/b}^{1/b} \biggl(
\frac{\Delta
_m^4(u)}{1+|u|}+\bigl(1+\Delta_m^4(u)\bigr)
\bigl(1+|u|\bigr)^{2\beta+1}\bigl|\phi_{\varepsilon,m}'(u)-
\phi_\varepsilon'(u)\bigr|^2 \biggr)\,\d u.
\end{eqnarray*}
Therefore, we find a constant $A'>0$, depending only on $R,\beta$,
such that (\ref{eqFMProp}) holds for
%
\begin{eqnarray}\label{eqOpNorm}
A_\psi&:=&A' \biggl(1+\int_{-1/b}^{1/b}
\biggl(\frac{\Delta
_m^2(u)+\Delta_m^4(u)}{1+|u|}
\nonumber\\[-8pt]\\[-8pt]
&&\hspace*{62pt}{} +\bigl(1+\Delta_m^4(u)\bigr)
\bigl(1+|u|\bigr)^{2\beta+1}\bigl|\phi_{\varepsilon,m}'(u)-\phi
_\varepsilon'(u)\bigr|^2 \biggr)\,\d u
\biggr)^{1/2}.\nonumber
\end{eqnarray}
The conditions (\ref{eqFMProp}) imply that $\psi$ is indeed a
Fourier multiplier on $B_\varepsilon(b)$ and thus by Theorem~4.8 and
Corollary 4.13 by Girardi and Weis \cite{girardi2003operator} with
$p=2$, $l=1$ there is a universal constant $C>0$ such that for all
$\eta>0$ and $f\in C^{s+\beta+\eta}(\R)$
\begin{eqnarray*}
\biggl\|\F^{-1} \biggl[\frac{\phi_K(bu)}{\phi_{\varepsilon,m}(u)} \biggr]\ast
f \biggr\|_{C^{s}}
&=& \biggl\|\F^{-1} \biggl[\frac{\phi_K(bu)}{\phi_{\varepsilon,m}(u)}\F f
\biggr] \biggr\|_{C^s}
\leq CA_\psi\bigl\|\F^{-1} \bigl[(1+\mathrm{i}u)^\beta\F f
\bigr] \bigr\| _{C^{s+\eta}}.
\end{eqnarray*}
Choosing $\eta>0$ such that $s+\beta+\eta,s+\eta\notin\N$, the
Fourier multiplier $(1+\mathrm{i}u)^\beta$ induces an isomorphism from
$C^{s+\beta+\eta}(\R)$ onto $C^{s+\eta}(\R)$
(Triebel \cite{triebel2010}, Thm. 2.3.8). Hence, there is another
universal constant $C'>0$ such that the second assertion of the lemma follows:
\begin{eqnarray*}
\biggl\|\F^{-1} \biggl[\frac{\phi_K(bu)}{\phi_{\varepsilon,m}(u)}\F f \biggr]
\biggr\|_{C^{s}}
\leq\mathcal E_b\|f\|_{C^{s+\beta+\eta
}}\qquad\mbox{with }\mathcal
E_b:=C'A_\psi.
\end{eqnarray*}
To bound $\mathcal E_b$, we apply Markov's inequality on $A_\psi$ from
(\ref{eqOpNorm}). The inequality by Rosenthal \cite{rosenthal1970} yields
\[
{\sup_{u\in\R}}\E\bigl[m^{p/2} \bigl|\phi_{\varepsilon,m}^{(l)}(u)-
\phi_{\varepsilon}^{(l)}(u) \bigr|^{p} \bigr]<\infty
\]
for $l=0$ and $p\in\N$ as well as $l=1$ and $p\in\{1,\dots,4\}$.
Combined with the Markov inequality and Cauchy--Schwarz inequality, we obtain
%
\begin{eqnarray}\label{eqBoundEb}
&& P \biggl(B_\varepsilon(b)\cap\biggl\{\mathcal E_b>
\frac{c^{1/2}}{m^{1/2}b^{\beta+1}\wedge1} \biggr\} \biggr)\nonumber 
\\
&&\quad \leq
c^{-1}\bigl(mb^{2\beta+2}\wedge1\bigr)\E\bigl[\mathcal
E_b^2\ind_{B_\varepsilon(b)} \bigr]
\nonumber
\\
&&\quad \lesssim \frac{1}{c}\bigl(mb^{2\beta+2}\wedge1\bigr) \biggl(1+
\int_{-1/b}^{1/b} \bigl(\bigl(1+|u|\bigr)^{-1}\E\bigl[
\Delta_m^2(u)+\Delta_m^4(u)\bigr]
\\
&&\hspace*{112pt}\qquad{}+\E\bigl[\bigl(1+\Delta_m^4(u)\bigr)
\bigl(1+|u|\bigr)^{2\beta+1}\bigl|\phi_{\varepsilon,m}'(u)-
\phi_\varepsilon'(u)\bigr|^2 \bigr] \bigr)\,\d u \biggr)
\nonumber
\\
&&\quad \lesssim\frac{mb^{2\beta+2}\wedge1}{c} \biggl(1+\frac
{1}{|\log b|^{3}}\int
_{-1/b}^{1/b}\frac{\d u}{1+|u|}+\frac{1}{m}
\int_{-1/b}^{1/b}\bigl(1+|u|\bigr)^{2\beta+1}\,\d u \biggr)
\lesssim\frac{1}{c},\nonumber
\end{eqnarray}
which shows $\mathcal E_b=\mathcal{O}_P(m^{-1/2}b^{-\beta-1}\vee1)$.

\subsubsection{\texorpdfstring{Proof of Proposition~\protect\ref{prDistUnknown}}{Proof of Proposition 2.5}}\label{secDistUnknown}
The following lemma establishes a bound for the bias term of the
estimator for the distribution function.

%
\begin{lemma}\label{lembias}
Let Assumption~\textup{\ref{asKernel}} hold with $\ell=\langle\alpha\rangle
+1$, $\alpha>0$ and $f(\bull+q_\tau)\in C^{\alpha}([-\zeta,\zeta
],R)$. Then we have
\[
\sup_{f(\bull+q_\tau)\in C^\alpha([-\zeta,\zeta],R)} \biggl|\int_{-\infty
}^{q_\tau}K_b
\ast f(x)\,\d x-\int_{-\infty}^{q_\tau}f(x)\,\d x \biggr|\leq D
b^{\alpha+1},
\]
where $D=(R/(\langle\alpha\rangle+1)!+2\zeta^{-\alpha-1})\|
K(x)x^{\alpha+1}\|_{L^1}$.
\end{lemma}

\begin{pf}
Let $F(x):=\int_{-\infty}^xf(y)\,\d y$.
Fubini's theorem yields
\begin{eqnarray}\label{eqDistConv}
\nonumber
\int_{-\infty}^{q_\tau}K_b\ast f(x)
\,\d x 
&=&\int_{-\infty}^\infty
K_b(x)F(q_\tau-x)\,\d x,
\end{eqnarray}
where $K_b(x):=b^{-1}K(x/b),x\in\R$. Therefore, the bias depends only
locally on $f$. Note that $F(\bull+q_\tau)\in C^{\alpha+1}([-\zeta,\zeta
])$ by assumption. A Taylor expansion of $F$ around $q_\tau$
yields for \mbox{$|bz|<\zeta$}
\begin{eqnarray*}
F(q_\tau-bz)-F(q_\tau)=-bzF^\prime(q_\tau)+
\cdots+(-bz)^{\langle\alpha\rangle+1}\frac{F^{(\langle\alpha
\rangle+1)}(q_\tau-\kappa bz)}{(\langle\alpha\rangle+1)!},
\end{eqnarray*}
where $0\leq\kappa\leq1$. Using the fact that $\int
x^kK(x)\,\d x=0$ for $k=1,\ldots,\langle\alpha\rangle+1$ and the
properties of the class, we obtain
\begin{eqnarray*}
&&\biggl|\int_{-\infty}^{q_\tau} \bigl(K_b\ast
f(x)-f(x) \bigr)\,\d x \biggr|
\\
&&\quad = \biggl|\int_{-\infty}^\infty K(z)
\bigl(F(q_\tau-bz)-F(q_\tau) \bigr)\,\d z \biggr|
\\
&&\quad \leq \biggl|\int_{|z|<\zeta/b} K(z) (-bz)^{\langle\alpha\rangle
+1}
\frac{F^{(\langle\alpha\rangle+1)}(q_\tau-\kappa
bz)-F^{(\langle\alpha\rangle+1)}(q_\tau)}{(\langle\alpha\rangle
+1)!}\,\d z \biggr|
\\
&&\qquad{}+\int_{|z|\geq\zeta/b}\bigl|K(z)\bigr|\bigl|F(q_\tau-bz)-F(q_\tau)\bigr|
\,\d z
\\
&&\quad \leq\frac{b^{\langle\alpha\rangle+1}R}{(\langle\alpha
\rangle+1)!}\int_{-\infty}^\infty\bigl|K(z)\bigr||z|^{\langle\alpha\rangle
+1}|\kappa bz|^{\alpha+1-(\langle\alpha\rangle+1)}\,\d z+2\int_{|z|\geq
\zeta/b}\bigl|K_b(z)\bigr|
\,\d z
\\
&&\quad \leq \biggl(\frac{b^{\alpha+1}R}{(\langle\alpha\rangle
+1)!}+2 \biggl(\frac{b}{\zeta}
\biggr)^{\alpha+1} \biggr)\int_{-\infty
}^\infty\bigl|K(z)\bigr||z|^{\alpha+1}\,\d z,
\end{eqnarray*}
and the statement follows.
\end{pf}

\begin{pf*}{Proof of Proposition~\ref{prDistUnknown}}
We will show uniformly over $f\in\mathcal C^\alpha(R,r,\zeta)$ and
$f_\varepsilon\in\mathcal D^\beta(R,\gamma)$ for any $b$ such that
$(n\wedge m)b^{2\beta+1}\to\infty$
\begin{eqnarray*}
&& \biggl|\int_{-\infty}^{q_\tau}\bigl(\widetilde
f_b(x)-f(x)\bigr)\,\d x \biggr|
\\
&&\quad =\mathcal{O}_P
\biggl(b^{\alpha+1}+\frac{1}{\sqrt{(n\wedge
m)(b^{2\beta-1}\wedge1)}}+\frac{1}{\sqrt{(n\wedge m)(mb^{2\beta
+2}\wedge1)}} \biggr).
\end{eqnarray*}
The third term on the right-hand side is of smaller or of the same
order than the second one if and only if $(mb^{1\wedge2\beta
+2})^{-1}\lesssim1$. Hence, when $\alpha\geq1/2$ the
asymptotically optimal choice $b=(n\wedge m)^{-1/(2\alpha+2(\beta\vee
1/2)+1)}$ yields
\[
\biggl|\int_{-\infty}^{q_\tau}\bigl(\widetilde
f_b(x)-f(x)\bigr)\,\d x \biggr| =\mathcal{O}_P \bigl((n\wedge
m)^{-(\alpha+1)/(2\alpha+2\beta+1)}\vee(n\wedge m)^{-1/2} \bigr).
\]

\textit{Step 1}:
As usual, we decompose the error into a deterministic error term and a
stochastic error term, writing $\phi_X=\F f$,
\begin{eqnarray*}
&& \biggl|\int_{-\infty}^{q_\tau}\bigl(\widetilde
f_b(x)-f(x)\bigr)\,\d x \biggr|
\\[-1.5pt]
&&\quad \leq \biggl|\int_{-\infty}^{q_\tau}
\bigl(K_b\ast f(x)-f(x)\bigr)\,\d x \biggr|
+ \biggl|\int_{-\infty}^{q_\tau}\F^{-1} \biggl[
\frac{\phi_n(u)\phi
_K(bu)}{\phi_{\varepsilon,m}(u)}-\phi_K(bu)\phi_X(u) \biggr](x)\,\d x
\biggr|.
\end{eqnarray*}
The bias is of order $\mathcal{O}(b^{\alpha+1})$ by Lemma~\ref
{lembias}. As discussed above, we decompose the stochastic error into
a singular part and a continuous one using a smooth truncation
function. Let $a_c\in C^\infty(\R)$ satisfy $a_c(x)=1$ for
$x\leq-1$ and $a_c(x)=0$ for $x\geq0$ and define
$a_s(x):=\mathbf\ind_{(-\infty,0]}(x)-a_c(x)$. Then\vspace*{-3pt}
%
\begin{eqnarray}\label{eqTsTc}
&&\int_{-\infty}^{q_\tau}\F^{-1} \biggl[
\phi_K(bu) \biggl(\frac{\phi
_n(u)}{\phi_{\varepsilon,m}(u)}-\phi_X(u) \biggr)
\biggr](x)\,\d x
\nonumber
\\[-1.5pt]
&&\quad =\int_{\R}a_s(x)\F^{-1} \biggl[
\phi_K(bu) \biggl(\frac{\phi
_n(u)}{\phi_{\varepsilon,m}(u)}-\phi_X(u) \biggr)
\biggr](x+{q_\tau})\,\d x
\nonumber\\[-9pt]\\[-9pt]
&&\qquad{} +\int_{\R}a_c(x)\F^{-1}
\biggl[\phi_K(bu) \biggl(\frac{\phi
_n(u)}{\phi_{\varepsilon,m}(u)}-\phi_X(u)
\biggr) \biggr](x+{q_\tau})\,\d x\nonumber
\\[-1.5pt]
&&\quad =:T_s+T_c.\nonumber
\end{eqnarray}
The singular term $T_s$ will be treated in the next step while we bound
the continuous, but not integrable term $T_c$ in Step~3.

\textit{Step 2}:
Lemma~\ref{lemPhiEps} shows that the probability of the complement
$B_\varepsilon(b)^c$ of $B_\varepsilon(b)$ from (\ref{eqBEps})
converges to zero. We obtain for any $c>0$ with Markov's inequality
\begin{eqnarray*}
&& P \biggl( |T_s |>\frac{c}{\sqrt{(n\wedge m)(b^{2\beta-1}\vee
1)}} \biggr)
\\[-1.5pt]
&&\quad \leq  P
\biggl(B_\varepsilon(b)\cap\biggl\{ |T_s |>\frac
{c}{\sqrt{(n\wedge m)(b^{2\beta-1}\vee1)}}
\biggr\} \biggr)+P \bigl(B_\varepsilon(b)^c \bigr)
\\[-1.5pt]
&&\quad \leq\frac{1}{c}\sqrt{(n\wedge m) \bigl(b^{2\beta-1}\vee1\bigr)}
\E\bigl[ |T_s |\ind_{B_\varepsilon(b)} \bigr]+\mathrm{o}(1).
\end{eqnarray*}
To bound $\E[|T_s|\ind_{B_\varepsilon(b)}]$, we first note by
Plancherel's identity
%
\begin{eqnarray}\label{eqtstx}
T_s&=&\frac{1}{2\pi}\int_{\R}\F
a_s(u)\mathrm{e}^{-\mathrm{i}uq_\tau}\phi_K(bu) \biggl(
\frac{\phi_n(u)}{\phi_{\varepsilon,m}(u)}-\phi_X(u) \biggr)\,\d u
\nonumber
\\[-1.5pt]
&=&\frac{1}{2\pi}\int_{\R}\F a_s(u)\mathrm{e}^{-\mathrm{i}uq_\tau}
\phi_K(bu) \biggl(\frac{\phi_n(u)}{\phi_\varepsilon(u)}-\phi_X(u) \biggr)
\,\d u
\nonumber\\[-9pt]\\[-9pt]
&&{}+\frac{1}{2\pi}\int_{\R}\F a_s(u)\mathrm{e}^{-\mathrm{i}uq_\tau}
\frac{\phi
_K(bu)\phi_n(u)}{\phi_{\varepsilon}(u)} \biggl(\frac{\phi
_{\varepsilon}(u)}{\phi_{\varepsilon,m}(u)}-1 \biggr)\,\d u \nonumber
\\[-1.5pt]
&=:&\frac{1}{2\pi}
(T_{s,x}+T_{s,\varepsilon} ).\nonumber
\end{eqnarray}
The first term, $T_{s,x}$ corresponds to the error due to the unknown
density $f$ while $T_{s,\varepsilon}$ is dominated by the error of the
estimator $\phi_{\varepsilon,m}$.
Since $a_s$ is of bounded variation and has compact support, there is a
constant $A_s\in(0,\infty)$ such that $|\F a_s(u)|\leq
A_s(1+|u|)^{-1}$. Plancherel's identity yields
%
\begin{eqnarray}\label{eqBoundTsX}
\Var(T_{s,x})&=&\E\bigl[|T_{s,x}|^2\bigr]
\leq\frac{1}{n}\E\biggl[ \biggl|\int_{\R}\F
a_s(u)\mathrm{e}^{-\mathrm{i}uq_\tau
}\frac{\phi_K(bu)}{\phi_\varepsilon(u)}\mathrm{e}^{\mathrm{i}uY_1}\,\d u
\biggr|^2 \biggr]\nonumber
\\
&\leq&\frac{4\pi^2}{n}\|f_Y\|_\infty\biggl\|
\F^{-1} \biggl[\F a_s(u)\frac{\phi_K(bu)}{\phi_\varepsilon(u)} \biggr]
\biggr\|_{L^2}^2
\nonumber\\[-8pt]\\[-8pt]
&\leq&\frac{4\pi^2}{n}\|K\|_{L^1}^2
\|f_Y\|_\infty\int_{-1/b}^{1/b}
\frac{|\F a_s(u)|^2}{|\phi_\varepsilon(u)|^2}\,\d u
\nonumber
\\
&\leq&\frac{4\pi^2}{n}\|K\|_{L^1}^2A_s^2
\|f_Y\|_\infty\int_{-1/b}^{1/b}
\frac{1}{(1+|u|)^2|\phi_\varepsilon(u)|^2}\,\d u.
\nonumber
\end{eqnarray}
Using the assumption $\|f\|_\infty<R$ and $f_\varepsilon\in\mathcal
D^\beta(R,\gamma)$, we get
%
\begin{equation}
\label{eqETsX} \E\bigl[|T_{s,x}|^2\bigr]\lesssim
\frac{1}{n}\int_{-1/b}^{1/b}\bigl(1+|u|\bigr)^{2\beta
-2}
\,\d u\lesssim\frac{1}{nb^{2\beta-1}}\vee\frac{1}{n}.
\end{equation}
To bound $T_{s,\varepsilon}$, we will use the following version of a
lemma by Neumann \cite{neumann1997effect}: by the definition~(\ref
{eqBEps}) of $B_\varepsilon(b)$ and applying (\ref{eqPhiInv})
it holds
%
\begin{eqnarray}\label{eqNeumann}
&& \E\biggl[ \biggl|\frac{\phi_{\varepsilon}(u)}{\phi_{\varepsilon,m}(u)}-1
\biggr|^{2}\ind_{B_\varepsilon(b)} \biggr]\nonumber
\\
&&\quad \leq 2\E\biggl[\frac{|\phi_{\varepsilon,m}(u)-\phi
_\varepsilon(u)|^2}{|\phi_\varepsilon(u)|^2} \biggr]+2 \E\biggl[\frac
{|\phi_{\varepsilon,m}(u)-\phi_\varepsilon
(u)|^4}{|\phi_\varepsilon(u)\phi_{\varepsilon,m}(u)|^2}\ind
_{B_\varepsilon(b)} \biggr]
\nonumber\\[-8pt]\\[-8pt]
&&\quad \leq \frac{2\E[|\phi_{\varepsilon,m}(u)-\phi_\varepsilon
(u)|^2]}{|\phi_\varepsilon(u)|^2} +\frac{2m\E[|\phi_{\varepsilon
,m}(u)-\phi_\varepsilon
(u)|^4]}{|\phi_\varepsilon(u)|^2}\nonumber
\\
&&\quad \leq \frac{18}{m|\phi_\varepsilon(u)|^2}.\nonumber
\end{eqnarray}
We estimate with the Cauchy--Schwarz inequality
\begin{eqnarray*}
T_{s,\varepsilon}^2&\leq&\|K\|_{L^1}^2\int
_{-1/b}^{1/b} \biggl|\frac{\phi_n(u)}{\phi_{\varepsilon}(u)} \biggr|^2\,\d u
\int_{-1/b}^{1/b}\bigl|\F a_s(u)\bigr|^2
\biggl|\frac{\phi_{\varepsilon}(u)}{\phi
_{\varepsilon,m}(u)}-1 \biggr|^2\,\d u
\nonumber
\\
&\leq&2\|K\|_{L^1}^2 \biggl(\|\phi_X
\|_{L^2}^2+\int_{-1/b}^{1/b}
\frac{|\phi_n(u)-\phi_Y(u)|^2}{|\phi_{\varepsilon
}(u)|^2}\,\d u \biggr)
\\
&&{}\times \int_{-1/b}^{1/b}\bigl|\F
a_s(u)\bigr|^2 \biggl|\frac{\phi
_{\varepsilon}(u)}{\phi_{\varepsilon,m}(u)}-1 \biggr|^2\,\d u.
\nonumber
%
\end{eqnarray*}
Applying again the Cauchy--Schwarz inequality, Fubini's theorem, the
decay of $\F a_s$ and (\ref{eqNeumann}), we obtain
%
\begin{eqnarray}\label{eqETsEps}
&& \E\bigl[|T_{s,\varepsilon}|\ind_{B_\varepsilon(b)}\bigr]\nonumber
\\
&&\quad \leq \sqrt{2}\|K
\|_{L^1} \biggl(\|\phi_X\|_{L^2}^2 +
\int_{-1/b}^{1/b}\frac{\E[|\phi_n(u)-\phi_Y(u)|^2]}{|\phi
_{\varepsilon}(u)|^2}\,\d u
\biggr)^{1/2}
\nonumber\\[-8pt]\\[-8pt]
&&\qquad{}\times\biggl(\int_{-1/b}^{1/b}\frac{A_s^2}{(1+|u|)^2}\E
\biggl[ \biggl|\frac{\phi_{\varepsilon}(u)}{\phi_{\varepsilon,m}(u)}-1 \biggr|^2\ind
_{B_\varepsilon(b)} \biggr]\,\d u
\biggr)^{1/2}\nonumber
\\
&&\quad \leq \frac{\sqrt{36}\|K\|_{L^1}A_s}{\sqrt{m}} \biggl(\|\phi_X\|
_{L^2}^2+\int_{-1/b}^{1/b}
\frac{\d u}{n|\phi_{\varepsilon
}(u)|^2} \biggr)^{1/2} \biggl(\int_{-1/b}^{1/b}
\frac{\d u}{(1+|u|)^2|\phi_{\varepsilon
}(u)|^2} \biggr)^{1/2}.\nonumber
\end{eqnarray}
The assumptions $\|f\|_\infty\lesssim1, |\phi_\varepsilon
(u)|\lesssim(1+|u|)^{-\beta}$ and $n^{-1}b^{-2\beta-1}\to0$ for the
optimal $b=b^*$ yield
\[
\E\bigl[|T_{s,\varepsilon}|\ind_{B_\varepsilon(b)}\bigr]\lesssim\biggl(1+\frac
{1}{nb^{2\beta+1}}
\biggr)^{1/2} \biggl(\frac{1}{\sqrt mb^{\beta
-1/2}}\vee\frac{1}{\sqrt m} \biggr)
\lesssim\frac{1}{\sqrt mb^{\beta
-1/2}}\vee\frac{1}{\sqrt m}.
\]
Together with (\ref{eqETsX}) and (\ref{eqtstx}) this implies the
optimal order
\[
\E\bigl[|T_s|\ind_{B_\varepsilon(b)}\bigr]\lesssim\bigl((n\wedge m)
\bigl(b^{2\beta
-1}\wedge1\bigr) \bigr)^{-1/2}.
\]

\textit{Step 3}: The empirical measures of $(Y_j)$ and $(\varepsilon
_k)$ are given by $\mu_{Y,n}:=\frac{1}{n}\sum_{j=1}^n\delta_{Y_j}$
and $\mu_{\varepsilon,m}:=\frac{1}{m}\sum_{k=1}^m\delta
_{\varepsilon_k}$, respectively, with Dirac measure $\delta_x$ in
$x\in\R$. We can write
\begin{eqnarray*}
T_c&=&\int_{\R}a_c(x)
\F^{-1} \biggl[\frac{\phi_K(bu)}{\phi
_{\varepsilon,m}(u)} \bigl(\phi_n(u)-
\phi_{\varepsilon,m}(u)\phi_X(u) \bigr) \biggr](x+q_\tau)
\,\d x
\\
&=&\F^{-1} \biggl[\frac{\phi_K(-bu)}{\phi_{\varepsilon,m}(-u)} \bigl(\phi_n(-u)-
\phi_{\varepsilon,m}(-u)\phi_X(-u) \bigr) \biggr]\ast
a_c(-q_\tau)
\\
&=&\F^{-1} \biggl[\frac{\phi_K(bu)}{\phi_{\varepsilon,m}(u)} \biggr]\ast
\bigl(
\mu_{Y,n}\ast a_c(-\bull)-\mu_{\varepsilon,m}\ast f\ast
a_c(-\bull) \bigr) (q_\tau).
\end{eqnarray*}
Applying Lemma~\ref{lemFourierMult}, we obtain on $B_\varepsilon(b)$
for any integer $s>\beta$
\begin{eqnarray*}
|T_c|&\leq&\biggl\|\F^{-1} \biggl[\frac{\phi_K(bu)}{\phi
_{\varepsilon,m}(u)}
\biggr]\ast\bigl(\mu_{Y,n}\ast a_c(-\bull)-\mu
_{\varepsilon,m}\ast f\ast a_c(-\bull) \bigr) \biggr\|_\infty
\\
&\leq&\mathcal E_b \bigl\|\mu_{Y,n}\ast
a_c(-\bull)-\mu_{\varepsilon,m}\ast f\ast a_c(-\bull)
\bigr\|_{C^s}
\\
&\lesssim& \mathcal E_b\sum
_{l=0}^s \bigl\|\mu_{Y,n}\ast
a_c^{(l)}(-\bull)-\mu_{\varepsilon,m}\ast f\ast
a_c^{(l)}(-\bull) \bigr\|_\infty.
\end{eqnarray*}
Therefore,
\begin{eqnarray*}
&&P \biggl(B_\varepsilon(b)\cap\biggl\{|T_c|>
\frac{c}{\sqrt{(n\wedge
m)}(\sqrt mb^{\beta+1}\wedge1)} \biggr\} \biggr)
\\[-3pt]
&&\quad\leq P \biggl(B_\varepsilon(b)
\cap\biggl\{\mathcal E_b> \biggl(\frac{c}{mb^{2\beta+2}\wedge1}
\biggr)^{1/2} \biggr\} \biggr)
\\[-3pt]
&&\qquad{}+P \Biggl(\sum_{l=0}^s \bigl\|
\mu_{Y,n}\ast a_c^{(l)}-\mu_{\varepsilon,m}\ast f
\ast a_c^{(l)} \bigr\|_\infty> \biggl(
\frac
{c}{n\wedge m} \biggr)^{1/2} \Biggr)
\\[-3pt]
&&\quad =:P_1+P_2.
\end{eqnarray*}
By Lemma~\ref{lemFourierMult}, more precisely estimate (\ref
{eqBoundEb}), the first probability is of the order $1/c$. To bound~$P_2$, it suffices to show $\|\mu_{Y,n}\ast a_c^{(l)}-\mu
_{\varepsilon,m}\ast f\ast a_c^{(l)} \|_\infty=\mathcal{O}_P((n\wedge m)^{-1/2})$ for all $l=0,\dots,s$. Denoting the density
of $Y_j$ as $f_Y=f\ast f_\varepsilon$, we decompose
\begin{eqnarray*}
&& \bigl\|\mu_{Y,n}\ast\bigl(a_c^{(l)}(-\bull)\bigr)-
\mu_{\varepsilon,m}\ast f\ast\bigl(a_c^{(l)}(-\bull)\bigr)
\bigr\|_\infty
\\[-2pt]
&&\quad \leq \bigl\|\mu_{Y,n}\ast\bigl(a_c^{(l)}(-
\bull)\bigr)-f_Y\ast\bigl(a_c^{(l)}(-\bull)
\bigr) \bigr\|_\infty
\\[-2pt]
&&\qquad{} + \bigl\|f_\varepsilon\ast\bigl(f\ast
\bigl(a_c^{(l)}(-\bull)\bigr)\bigr)-\mu_{\varepsilon,m}\ast
\bigl(f\ast\bigl(a_c^{(l)}(-\bull)\bigr)\bigr)
\bigr\|_\infty
\\[-2pt]
&&\quad \leq \biggl\|\int a_c^{(l)}(y-\bull)\mu_{Y,n}(\d
y) -\E\bigl[a_c^{(l)}(Y_1-\bull)\bigr]
\biggr\|_\infty
\\[-3pt]
&&\qquad{} + \biggl\|\E\bigl[\bigl(f\ast a_c^{(l)}\bigr) (
\varepsilon_1-\bull)\bigr]-\int\bigl(f\ast a_c^{(l)}
\bigr) (z-\bull)\mu_{\varepsilon,m}(\d z) \biggr\|_\infty.
\end{eqnarray*}
By construction all $a_c^{(l)}, l\geq1$, have compact support and
are bounded. Therefore, $\|a_c^{(l)}\|_{L^1}<\infty,\|(a_c\ast
f)^{(l)}\|_{L^1}\leq\|a_c^{(l)}\|_{L^1}\|f\|_{L^1}<\infty$ and
thus $a_c^{(l)}(\bull-t)$ and $a_c^{(l)}\ast f(\bull-t)$, $l\geq0$,
are of bounded variation for all $t\in\R$. Since the set of functions
with bounded variation is a Donsker class (cf. Theorem 2.1 by Dudley
\cite{dudley1992}), the two terms in the previous display converge in
probability to a tight limit with $\sqrt{n}$-rate and $\sqrt
{m}$-rate, respectively. Consequently,
\[
\sqrt{n\wedge m} \bigl\|\mu_{Y,n}\ast\bigl(a_c^{(l)}(-
\bull)\bigr)-\mu_{\varepsilon,m}\ast f\ast\bigl(a_c^{(l)}(-
\bull)\bigr) \bigr\|_\infty=\mathcal{O}_P(1)
\]
for all $\ell=0,\dots,s$ and $P_2$ is arbitrary small for $c$ large.
\end{pf*}

For the adaptive estimator, we will later need the following uniform
version of Proposition~\ref{prDistUnknown}.

\begin{corollary}\label{coruniformQ}
Suppose Assumption~\textup{\ref{asKernel}} holds with $l=\langle\alpha
\rangle+1$ and let the set $\mathcal B=\mathcal B_n$ be given by~(\ref{eqbandwidthset}).
For critical values $(\delta_b)_{b\in\mathcal B}$ satisfying $\delta
_b>3Db^{\alpha+1}$ and for any sequence $(x_n)_n$ with $x_n\to\infty
$ arbitrarily slowly we obtain uniformly in $\mathcal C^\alpha
(R,r,\zeta)$ and $\mathcal D^\beta(R,\gamma)$
\begin{eqnarray*}
&&P \biggl(\exists b\in\mathcal B\dvt  \biggl|\int_{-\infty}^{q_\tau}
\bigl(\widetilde f_b(x)-f(x) \bigr)\,\d x \biggr|>\delta_b
\biggr)
\\[-3pt]
&&\quad=\mathcal{O} \biggl(\sum_{b\in\mathcal B} \biggl(
\frac{1}{\delta
_b} \bigl((n\wedge m) \bigl(b^{2\beta-1}\wedge1\bigr)
\bigr)^{-1/2}+\frac
{1}{\delta_b^2}\frac{x_n}{(n\wedge m)(mb^{2\beta+2}\wedge1)} \biggr)
\biggr)+\mathrm{o}(1).
\end{eqnarray*}
In particular, if $|\mathcal B|\lesssim\log n, \max_{b\in\B}b\to0$
and $\min_{b\in\mathcal B}(n\wedge m)b^{2\beta+1}\to\infty$, then
\[
\sup_{b\in\mathcal B} \biggl|\int_{-\infty}^{q_\tau}
\bigl(\widetilde f_b(x)-f(x)\bigr)\,\d x \biggr|\stackrel{P} {\rightarrow}0.
\]
\end{corollary}

\begin{pf}
With the notation of the proof of Proposition~\ref{prDistUnknown} and
applying Lemma~\ref{lembias}, we obtain
\begin{eqnarray*}
\biggl|\int_{-\infty}^{q_\tau}\bigl(\widetilde
f_b(x)-f(x)\bigr)\,\d x \biggr| &\leq& \biggl|\int_{-\infty}^{q_\tau}
\bigl(K_b\ast f(x)-f(x) \bigr)\,\d x \biggr|+|T_s|+|T_c|
\\
&\leq& Db^{\alpha+1}+|T_s|+|T_c|,
\end{eqnarray*}
where $T_s$ and $T_c$ are the stochastic errors of the singular part
and of the continuous part, respectively, as defined in (\ref
{eqTsTc}). Since both terms depend on $b$ let us write $T_s(b)$ and
$T_c(b)$. By definition $b_1\leq b$ implies $B_\varepsilon
(b_1)\subset B_\varepsilon(b)$. Then, Step~2 in the previous proof shows
\begin{eqnarray*}
P (\exists b\in\mathcal B\dvt T_s>\delta_b/3 )
&\leq& \biggl(\sum_{b\in\mathcal B}P
\bigl(\bigl\{T_s(b)>\delta_b/3\bigr\} \cap
B_\varepsilon(b_1) \bigr) \biggr)+\mathrm{o}(1)
\\
&\leq& \biggl(\sum_{b\in\mathcal B}\delta_b^{-1}
\E\bigl[\bigl|T_s(b)\bigr|\ind_{B_\varepsilon(b_1)}\bigr] \biggr)+\mathrm{o}(1)
\\
&\lesssim& \biggl(\sum_{b\in\mathcal B}\delta_b^{-1}
\bigl((n\wedge m) \bigl(b^{2\beta-1}\wedge1\bigr) \bigr)^{-1/2}
\biggr)+\mathrm{o}(1).
\end{eqnarray*}
Following Step~3 in the previous proof, we obtain with the random
operator norm $\mathcal E_b$, for some integer $s>\beta$ and for a
diverging sequence $(x_{(n\wedge m)})$
\begin{eqnarray*}
&& P (\exists b\in\mathcal B\dvt T_c>\delta_b/3 )
\\
&&\quad \leq
P \bigl( \bigl\{\exists b\in\mathcal B\dvt \mathcal E_b>\delta
_b (n\wedge m)^{1/2}/\bigl(3(x_{(n\wedge m)})^{1/2}
\bigr) \bigr\}\cap B_\varepsilon(b_1) \bigr)+P
\bigl(B_\varepsilon(b_1)^c\bigr)
\\
&&\qquad{}+P \Biggl( \Biggl\{\sum_{l=0}^s \bigl\|
\mu_{Y,n}\ast a_c^{(l)}-\mu_{\varepsilon,m}\ast f
\ast a_c^{(l)} \bigr\|_\infty> \biggl(
\frac{x_{(n\wedge m)}}{(n\wedge m)} \biggr)^{1/2} \Biggr\} \Biggr)
\\
&&\quad \leq \biggl(\sum_{b\in\mathcal B}P \bigl(\bigl\{\mathcal
E_b> \delta_b (n\wedge m)^{1/2}/
\bigl(3(x_{(n\wedge m)})^{1/2}\bigr)\bigr\}\cap B_\varepsilon
(b_1) \bigr) \biggr)+\mathrm{o}(1)
\\
&&\quad \lesssim \biggl(\sum_{b\in\mathcal B}\frac{x_{(n\wedge m)}}{\delta
_b^2(n\wedge m)(mb^{2\beta+2}\wedge1)}
\biggr)+\mathrm{o}(1),
\end{eqnarray*}
where we have used (\ref{eqBoundEb}) in the last estimate.
\end{pf}

\subsubsection{\texorpdfstring{Proof of Proposition~\protect\ref{pruniform}}{Proof of Proposition 2.6}}
Without loss of generality, we set $q_\tau=0$. Recall definition~(\ref
{eqfhat}) of the pseudo-estimator $\widehat f_b$ which knows the
error distribution. We estimate
\begin{eqnarray*}
\sup_{x\in(-\zeta,\zeta)}\bigl|\widetilde f_b(x)-f(x)\bigr| &\leq&
\sup_{x\in(-\zeta,\zeta)}\bigl|\widehat f_b(x)-f(x)\bigr| +\|\widetilde
f_b-\widehat f_b\|_\infty
\\
&\leq&\sup_{x\in(-\zeta,\zeta)}\bigl|\widehat f_b(x)-f(x)\bigr| + \biggl\|
\frac{\phi_K(bu)\phi_n(u)}{\phi_\varepsilon(u)} \biggl(\frac
{\phi_\varepsilon(u)}{\phi_{\varepsilon,m}(u)}-1 \biggr) \biggr\|_{L^1}.
\end{eqnarray*}
The analysis of the first term is very classical. However, we are not
aware of any reference in the given setup. Both terms will be treated
separately in the following two steps. All estimates will be uniform in
$f\in\mathcal C^\alpha(R,r,\zeta)$ and $f_\varepsilon\in\mathcal
D^\beta(R,\gamma)$.

\textit{Step 1}: Let $b\in(0,1)$. We will show that there are
constants $d,D>0$ such that for any $t>d(b^\alpha+(nb^{2\beta+1})^{-2})$
%
\begin{equation}
\label{eqConcentrationDensity} P \Bigl(\sup_{x\in(-\zeta,\zeta
)}\bigl|\widehat
f_b(x)-f(x)\bigr|>t \Bigr)\leq2\exp\bigl(2\log
n-Dnb^{(2\beta+1)}\bigl(t\wedge t^2\bigr) \bigr).
\end{equation}
Then the result follows by choosing $t\sim b^\alpha+ (\frac{\log
n}{nb^{2\beta+1}} )^{1/2}$. Let us define $x_k:=-\zeta+kn^{-2}$
for $k=1,\dots,\lfloor2\zeta n^2\rfloor=:M$ as well as
\begin{eqnarray*}
\chi_{j}(x)&:=&\F^{-1} \biggl[\frac{\phi_K(bu)}{\phi_\varepsilon
(u)}\mathrm{e}^{\mathrm{i}uY_j}
\biggr](x)-\E\biggl[\F^{-1} \biggl[\frac{\phi_K(bu)}{\phi
_\varepsilon(u)}\mathrm{e}^{\mathrm{i}uY_j}
\biggr](x) \biggr]
\\
&=&K_b\ast\F^{-1} \biggl[\ind_{[-b^{-1},b^{-1}]}(u)
\frac
{\mathrm{e}^{\mathrm{i}uY_j}}{\phi_\varepsilon(u)} \biggr](x)-K_b\ast f(x),\qquad x\in\R.
\end{eqnarray*}
Therefore, $\widehat f_b(x)-\E[\widehat f_b(x)]=\frac{1}{n}\sum
_{j=1}^n\chi_j(x)$ and thus
\begin{eqnarray*}
\sup_{|x|<\zeta}\bigl|\widehat f_b(x)-f(x)\bigr| &\leq&\sup
_{|x|<\zeta}\bigl|\E\bigl[\widehat f_b(x)\bigr]-f(x)\bigr|+\sup
_{|x|<\zeta}\bigl|\widehat f_b(x)-\E\bigl[\widehat
f_b(x)\bigr]\bigr|
\\
&\leq&\sup_{|x|<\zeta}\bigl|\E\bigl[\widehat f_b(x)
\bigr]-f(x)\bigr|+ \sup_{|x|<\zeta}\min_{k=1,\dots,M} \Biggl|
\frac{1}{n}\sum_{j=1}^n \bigl(
\chi_j(x)-\chi_j(x_k) \bigr) \Biggr|
\\
&&{} +\max
_{k=1,\dots,M} \Biggl|\frac{1}{n}\sum_{j=1}^n
\chi_j(x_k) \Biggr|
\\
&=:& B+V_1+V_2.
\end{eqnarray*}
The bias term $B$ can be bounded as in the classical density estimation
setup (cf. also Fan \cite{fan1991}, Thms. 1 and 2), noting that the
constant does not depend on $x\in(-\zeta,\zeta)$. Hence,
$|B|\lesssim b^{\alpha}$. Using a continuity argument and the
properties of $f_\varepsilon\in\mathcal D^\beta(R,\gamma)$, the
term $V_1$ can be bounded by
\begin{eqnarray*}
|V_1|&\leq&\frac{1}{n^2} \Biggl\|\frac{1}{n}\sum
_{j=1}^n\chi_j'
\Biggr\|_\infty
\\
&=&\frac{1}{n^3} \Biggl\|\sum_{j=1}^n
\bigl(K_b'\bigr)\ast\biggl(\F^{-1} \biggl[
\ind_{[-b^{-1},b^{-1}]}(u)\frac{\mathrm{e}^{\mathrm{i}uY_j}}{\phi_\varepsilon(u)} \biggr
]-f \biggr) \Biggr\|_\infty
\\
&\leq&\frac{1}{n^2b}\bigl\|K'\bigr\|_{L^1} \bigl(\bigl\|\ind
_{[-b^{-1},b^{-1}]}\phi_\varepsilon^{-1}\bigr\|_{L^1}+\|f
\|_\infty\bigr) \lesssim n^{-2}b^{-(\beta+2)}
\lesssim
\bigl(nb^{2\beta+1}\bigr)^{-2}.
\end{eqnarray*}
Therefore, $|B+V_1|\leq D_1(b^\alpha+(nb^{2\beta+1})^{-2})$ for
some constant $D_1>0$. We obtain for all $t>d(b^\alpha+(nb^{2\beta
+1})^{-2})$ with $d:=2D_1$
\begin{eqnarray*}
P \Bigl(\sup_{|x|<\zeta} \bigl|\widehat f_b(x)-f(x) \bigr|>t
\Bigr) &\leq& P \Biggl(\max_{k=1,\dots,M} \Biggl|\frac{1}{n}\sum
_{j=1}^n\chi_j(x_k)
\Biggr|>\frac{t}{2} \Biggr)
\\
& \leq& \sum_{k=1}^MP
\Biggl( \Biggl|\frac{1}{n}\sum_{j=1}^n
\chi_j(x_k) \Biggr|>\frac{t}{2} \Biggr).
\end{eqnarray*}
Finally, we will apply Bernstein's inequality. To this end, we estimate
\[
\max_{j,k}\bigl|\chi_{j}(x_k)\bigr|\leq2
\|K_b\|_{L^1}\bigl\|\ind_{[-b^{-1},b^{-1}]}\phi_\varepsilon^{-1}
\bigr\|_{L^1}\leq D_2b^{-(\beta+1)},
\]
with some constant $D_2>0$. Using Plancherel's identity, the variance
can be estimated by
\begin{eqnarray*}
\Var\bigl(\chi_{j}(x_k)\bigr)&=&\E\biggl[
\F^{-1} \biggl[\frac{\phi_K(bu)}{\phi
_\varepsilon(u)}\mathrm{e}^{\mathrm{i}uY_j} \biggr]^2(x_k)
\biggr]-(K_b\ast f)^2(x_k)
\\
&\leq&\frac{1}{2\pi}\|f\|_\infty\biggl\|\frac{\phi
_K(-bu)}{\phi_\varepsilon(-u)}
\biggr\|^2_{L^2} \lesssim D_3b^{-(2\beta+1)},
\end{eqnarray*}
for some $D_3>0$. Then Bernstein's inequality yields
\begin{eqnarray*}
P \Bigl(\sup_{x\in(-\zeta,\zeta)} \bigl|\widehat f_b(x)-f(x) \bigr|>t
\Bigr) &\leq&\sum_{k=1}^MP \Biggl( \Biggl|
\sum_{j=1}^n\chi_{j}(x_k)
\Biggr|>nt/2 \Biggr)
\\
&\leq&2\exp\biggl(\log M-\frac{nb^{(2\beta
+1)}t^2}{8(D_3+D_2t/3)} \biggr)
\\
&\leq& 2\exp\bigl(2\log n-Dnb^{(2\beta+1)}\bigl(t\wedge t^2
\bigr) \bigr),
\end{eqnarray*}
with some constant $D>0$.

\textit{Step 2}: By the Cauchy--Schwarz inequality, we have
\begin{eqnarray*}
&&\E\biggl[ \biggl\|\frac{\phi_K(bu)\phi_n(u)}{\phi_\varepsilon
(u)} \biggl(\frac{\phi_\varepsilon(u)}{\phi_{\varepsilon,m}(u)}-1 \biggr)
\biggr\|_{L^1}\ind_{B_\varepsilon(b)} \biggr]
\\
&&\quad \lesssim \biggl(\E\biggl[ \biggl\|\frac{\phi_n(u)}{\phi_\varepsilon
(u)}\ind_{[-1/b,1/b]}(u)
\biggr\|^2_{L^2} \biggr]\E\biggl[ \biggl\| \biggl(\frac{\phi_{\varepsilon}(u)}{\phi
_{\varepsilon,m}(u)}-1
\biggr)\ind_{[-1/b,1/b]}(u) \biggr\|^2_{L^2}
\ind_{B_\varepsilon(b)} \biggr] \biggr)^{1/2}
\\
&&\quad \leq \biggl(\|\phi_X\|_{L^2}
+ \biggl(\int
_{-1/b}^{1/b}\frac{\E
[|\phi_n(u)-\phi_Y(u)|^2]}{|\phi_{\varepsilon}(u)|^2}\,\d u
\biggr)^{1/2} \biggr)
\\
&&\qquad{}\times \biggl(\int_{-1/b}^{1/b}
\E\biggl[ \biggl|\frac{\phi
_{\varepsilon}(u)}{\phi_{\varepsilon,m}(u)}-1 \biggr|^2\ind_{B_\varepsilon
(b)} \biggr]\,\d u
\biggr)^{1/2}
\\
&&\quad \lesssim \biggl(\|\phi_X\|_{L^2}+ \biggl(
\frac{1}{nb^{2\beta+1}} \biggr)^{1/2} \biggr) \biggl(\frac{1}{mb^{2\beta+1}}
\biggr)^{1/2},
\end{eqnarray*}
where we have used (\ref{eqNeumann}) for the last step.
Therefore, the additional error due to the unknown error distribution
satisfies for any $\delta>0$ by Markov's inequality and by Lemma~\ref
{lemPhiEps}
%
\begin{eqnarray}\label{eqconcBoundEps}
&& P \biggl( \biggl\|\frac{\phi_K(bu)\phi_n(u)}{\phi_\varepsilon(u)} \biggl(\frac
{\phi_\varepsilon(u)}{\phi_{\varepsilon,m}(u)}-1 \biggr) \biggr\|
_{L^1}>\delta\biggr)\nonumber
\\
&&\quad \leq \frac{1}{\delta}\E\biggl[ \biggl\|
\frac
{\phi_K(bu)\phi_n(u)}{\phi_\varepsilon(u)} \biggl(\frac{\phi
_\varepsilon(u)}{\phi_{\varepsilon,m}(u)}-1 \biggr) \biggr\|_{L^1}\ind
_{B_\varepsilon(b)} \biggr]
+P \Bigl(\inf_{|u|\leq1/b}\bigl|\phi_{\varepsilon,m}(u)\bigr|<
m^{-1/2} \Bigr)\qquad
\\
&&\quad \lesssim \frac{1}{\delta} \biggl(\frac{1}{mb^{2\beta+1}} \biggr
)^{1/2}+\mathrm{o}(1)\nonumber
\end{eqnarray}
and thus $\|\widetilde f_b-\widehat f_b\|_\infty=\mathcal{O}_P((mb^{2\beta+1})^{-1/2})$. Note that the second term does not
depend on $\delta$ and thus $\mathrm{o}(1)$ is sufficient.

\subsubsection{\texorpdfstring{Proof of Theorems~\protect\ref{thRateUnk} and \protect\ref{thuniformLoss}}{Proof of Theorems 2.7 and 2.8}}
We start with a lemma that establishes consistency of the quantile
estimator and then prove the theorems. To apply this lemma also for the
adaptive result, we prove convergence uniformly over a set of bandwidths.

%
\begin{lemma}\label{lemConsitency}
Grant Assumption~\textup{\ref{asKernel}} with $\ell=1$. Let $\mathcal B$ be a
set of bandwidths satisfying $|\mathcal B|\lesssim\log n, \max
\mathcal B \to0$ and $\min_{b\in\mathcal B}(\log n)^2/((n\wedge
m)b^{2\beta+1})\to0$. Then
\[
\sup_{f\in\mathcal C^\alpha(R,r,\zeta, {U_n})}\sup_{f_\varepsilon
\in\mathcal D^\beta(R,\gamma)}P \Bigl(\sup
_{b\in\mathcal
B}|\widetilde q_{\tau,b}-q_\tau|>\delta
\Bigr)\to0\qquad\mbox{for all }\delta>0.
\]
\end{lemma}

\begin{pf}
{We follow the general strategy of the proof of Theorem 5.7 by van~der
Vaart \cite{vanderVaart1998} in the classical M-estimation setting.
Recall the definition of $\widetilde M_b$ given in (\ref{eqestEq})
and its deterministic counterpart $M(\eta)=\int_{-\infty}^\eta
f(x)\,\d x-\tau$. To this end, we first claim that
%
\begin{equation}
\label{eqZEst} \sup_{b\in\mathcal B}\widetilde M_b(
\widetilde q_{\tau,b})=\mathrm{o}_P(1),
\end{equation}
Since $\widetilde q_{\tau,b}$ minimizes $\widetilde M_b$ on the
interval $[-U_n,U_n]$ for $U_n\lesssim\log n$ and $M(q_\tau)=0$ with
$q_\tau\in[-U_n,U_n]$, Corollary~\ref{coruniformQ} implies for any
$\delta>0$
%
\begin{eqnarray}\label{eqApprox0}
P \Bigl(\sup_{b\in\mathcal B}\bigl|\widetilde M_b(\widetilde
q_{\tau,b})\bigr|>\delta\Bigr) &\leq& P \Bigl(\sup_{b\in\mathcal B}\bigl|
\widetilde M_b(q_\tau)-M(q_\tau)\bigr|>\delta
\Bigr)
\nonumber\\[-8pt]\\[-8pt]
&=&P \biggl(\sup_{b\in\mathcal B}\biggl|\int_{-\infty}^{q_\tau}
\bigl(\widetilde f_b(x)-f(x)\bigr)\,\d x\biggr|>\delta\biggr)\to0,
\nonumber
\end{eqnarray}
which gives (\ref{eqZEst}).}

Now, we show that $f$ satisfies the uniqueness condition
%
\begin{equation}
\label{eqUnique} \inf_{\eta\dvtx |\eta-q_\tau|\geq\delta}\bigl|M(\eta
)\bigr|>0\qquad\mbox{for any }
\delta>0.
\end{equation}
By the H\"older regularity $M'(\eta)=f(\eta)\geq f(q_\tau
)-|f(q_\tau)-f(\eta)|\geq r-R|q_\tau-\eta|^{1\wedge\alpha
}\geq r/2$ for $|q_\tau-\eta|\leq(\frac{r}{2R})^{1\vee
\alpha^{-1}}$. Without loss of generality, we can assume $\delta
\leq(\frac{r}{2R})^{1\vee\alpha^{-1}}$, otherwise consider
$\delta\wedge(\frac{r}{2R})^{1\vee\alpha^{-1}}$. Recall that
$q_\tau$ is given by the root of $M$ and that $M$ is increasing.
Hence, we obtain
\begin{eqnarray*}
\inf_{\eta\dvtx |\eta-q_\tau|\geq\delta}\bigl|M(\eta)\bigr| &=&\inf_{\eta\in\{
-\delta,\delta\}}
\bigl|M(q_\tau-\eta)-M(q_\tau) \bigr| \geq\delta\inf
_{\eta\dvtx |\eta-q_\tau|\geq\delta}M'(\eta) \geq\frac{\delta r}{2}.
\end{eqnarray*}
Applying (\ref{eqZEst}) and (\ref{eqUnique}) yield
%
\begin{eqnarray}\label{eqConsistencyBound}
P \Bigl(\sup_{b\in\mathcal B}|\widetilde q_{\tau,b}-q_\tau|>
\delta\Bigr) &\leq& P \Bigl(\sup_{b\in\mathcal B}\bigl|M(\widetilde
q_{\tau,b})\bigr|\geq\delta r/2 \Bigr)
\nonumber
\\
&=& {P \Bigl(\sup_{b\in\mathcal B}\bigl|M(\widetilde q_{\tau,b})-
\widetilde M_b(\widetilde q_{\tau,b})\bigr|\geq\delta r/3
\Bigr)+\mathrm{o}(1)}
\nonumber\\[-8pt]\\[-8pt]
&\leq& P \Bigl(\sup_{b\in\mathcal B}\sup_{\eta\in
[-U_n,U_n]}\bigl|M(
\eta)-\widetilde M_b(\eta)\bigr|\geq\delta r/3 \Bigr)+\mathrm{o}(1)
\nonumber
\\
&=&P \biggl(\sup_{b\in\mathcal B}\sup_{\eta\in[-U_n,U_n]} \biggl|\int
_{-\infty}^\eta\bigl(\widetilde f_b(x)-f(x)
\bigr)\,\d x \biggr|\geq\delta r/3 \biggr)+\mathrm{o}(1).\nonumber
\end{eqnarray}
Hence, it remains to show uniform consistency of $\int_{-\infty}^\eta
\widetilde f_b(x)\,\d x$. Write
\begin{eqnarray*}
\biggl|\int_{-\infty}^\eta\bigl(\widetilde
f_b(x)-f(x)\bigr)\,\d x \biggr| &\leq& \biggl|\int_{-\infty}^\eta
\bigl(K_b\ast f(x)-f(x)\bigr)\,\d x \biggr|+ \biggl|\int_{-\infty}^\eta
\bigl(\widetilde f_b(x)-K_b\ast f(x)\bigr)\,\d x \biggr|
\\
&=&\bigl|K_b\ast F(\eta)-F(\eta)\bigr|+ \biggl|\int_{-\infty}^\eta
\bigl(\widetilde f_b(x)-K_b\ast f(x)\bigr)\,\d x \biggr|.
\end{eqnarray*}
We have $|K_b\ast F(\eta)-F(\eta)|=|\int K_b(z)(F(\eta-z)-F(\eta
))\,\d z|\leq b\|f\|_\infty\|zK(z)\|_{L^1}$ by the boundedness of
$f$. Further note for $\eta\in[-U_n,U_n]$
\begin{eqnarray*}
&& \biggl|\int_{-\infty}^\eta\bigl(\widetilde
f_b(x)-K_b\ast f(x)\bigr)\,\d x \biggr|
\\
&&\quad \leq \biggl|\int_{-\infty}^{q_\tau}\bigl(\widetilde
f_b(x)-K_b\ast f(x)\bigr)\,\d x \biggr| + \biggl|\int
_{q_\tau\wedge\eta}^{q_\tau\vee\eta}\bigl(\widetilde f_b(x)-K_b
\ast f(x)\bigr)\,\d x \biggr|
\\
&&\quad \leq \biggl|\int_{-\infty}^{q_\tau}\bigl(\widetilde
f_b(x)-K_b\ast f(x)\bigr)\,\d x \biggr|+\sqrt{2U_n}
\biggl(\int_{-\infty}^\infty\bigl(\widetilde
f_b(x)-K_b\ast f(x)\bigr)^2\,\d x
\biggr)^{1/2},
\end{eqnarray*}
where we have used the Cauchy--Schwarz inequality for the last step.
Hence, together with (\ref{eqConsistencyBound}) we obtain for all
$\delta>6\|f\|_\infty\|zK(z)\|_{L^1}/r \sup_{b\in\mathcal B}b$
\begin{eqnarray*}
&& P \Bigl(\sup_{b\in\mathcal B}|\widetilde q_{\tau,b}-q_\tau|>
\delta\Bigr)
\\
&&\quad \leq P \biggl(\sup_{b\in\mathcal B}\sup
_{\eta\in[-U_n,U_n]} \biggl|\int_{-\infty}^\eta\bigl(
\widetilde f_b(x)-f(x)\bigr)\,\d x \biggr|\geq\delta r/3 \biggr)+\mathrm{o}(1)
\\
&&\quad \leq P \biggl(\sup_{b\in\mathcal B} \biggl|\int_{-\infty
}^{q_\tau}
\bigl(\widetilde f_b(x)-K_b\ast f(x)\bigr)\,\d x \biggr|
\geq\frac
{\delta r}{9} \biggr)
\\
&&\qquad{} +P \biggl(\sup_{b\in\mathcal B}\int
_{\R}\bigl(\widetilde f_b(x)-K_b
\ast f(x)\bigr)^2\,\d x\geq\frac{\delta^2 r^2}{162U_n} \biggr).
\end{eqnarray*}
Corollary~\ref{coruniformQ} shows under the conditions on $\mathcal
B$ that
\[
P \biggl(\sup_{b\in\mathcal B} \biggl|\int_{-\infty}^{q_\tau
}
\bigl(\widetilde f_b(x)-K_b\ast f(x)\bigr)\,\d x \biggr|>
\delta r/9 \biggr)\to0.
\]
Hence, it remains to show
%
\begin{equation}
\label{eqDenUnk} P \biggl(\sup_{b\in\mathcal B}\int_{\R}
\bigl(\widetilde f_b(x)-K_b\ast f(x)
\bigr)^2\,\d x>\delta^2 r^2/(162U_n)
\biggr)\to0.
\end{equation}
On the event $B_\varepsilon(b)$, (\ref{eqDenUnk}) follows basically
from the work of Neumann \cite{neumann1997effect}. More precisely,
Plancherel's equality, (\ref{eqNeumann}) and the Cauchy--Schwarz
inequality yield for any $b\in\mathcal B$
\begin{eqnarray*}
&&\E\biggl[\int_{\R}\bigl(\widetilde f_b(x)-K_b
\ast f(x)\bigr)^2\,\d x\ind_{B_\varepsilon(b)} \biggr]
\\
&&\quad =\frac{1}{2\pi}\int_{\R}\bigl|\phi_K(bu)\bigr|^2
\E\biggl[ \biggl|\frac{\phi
_n(u)}{\phi_{\varepsilon,m}(u)}-\frac{\phi_Y(u)}{\phi_\varepsilon
(u)} \biggr|^2
\ind_{B_\varepsilon(b)} \biggr]\,\d u
\\
&&\quad \lesssim\int_{-1/b}^{1/b} \biggl(\E\biggl[
\frac{|\phi_n(u)-\phi
_Y(u)|^2}{|\phi_{\varepsilon,m}(u)|^2}\ind_{B_\varepsilon(b)} \biggr
]+\bigl|\phi_Y(u)\bigr|^2
\E\biggl[ \biggl|\frac{1}{\phi_{\varepsilon,m}(u)}-\frac{1}{\phi_{\varepsilon
}(u)} \biggr|^2
\ind_{B_\varepsilon
(b)} \biggr] \biggr)\,\d u
\\
&&\quad \lesssim\int_{-1/b}^{1/b} \biggl(\E\biggl[
\frac{|\phi_n(u)-\phi
_Y(u)|^2}{|\phi_{\varepsilon}(u)|^2} \bigl(1+m\bigl|\phi_{\varepsilon
,m}(u)-\phi_\varepsilon(u)\bigr|^2
\bigr) \biggr]+\frac{|\phi_Y(u)|^2}{m|\phi_{\varepsilon}(u)|^4} \biggr
)\,\d u
\\
&&\quad \leq\int_{-1/b}^{1/b}\frac{1}{|\phi_\varepsilon(u)|^{2}}
\\
&&\hspace*{42pt}{}\times \biggl( \bigl(\E\bigl[\bigl|\phi_n(u)-\phi_Y(u)\bigr|^4
\bigr]\E\bigl[2+2m^2\bigl|\phi_{\varepsilon,m}(u)-\phi_\varepsilon(u)\bigr|^4
\bigr] \bigr)^{1/2}
+\frac
{|\phi_X(u)|^2}{m} \biggr)\,\d u
\\
&&\quad \lesssim\int_{-1/b}^{1/b}\bigl|\phi_\varepsilon
(u)\bigr|^{-2}\bigl(n^{-1}+m^{-1}\bigr)\,\d u \lesssim
\frac{1}{(n\wedge m)b^{2\beta+1}}.
\end{eqnarray*}
Using $B_\varepsilon(\min\mathcal B)\subset B_\varepsilon(b)$ and
Lemma~\ref{lemPhiEps}, (\ref{eqDenUnk}) follows from Markov's inequality
\begin{eqnarray*}
&& P \biggl(\sup_{b\in\mathcal B}\int_{\R}\bigl(
\widetilde f_b(x)-K_b\ast f(x)\bigr)^2\,\d
x>\delta^2 r^2/(162U_n) \biggr)
\\
&&\quad \lesssim\frac{U_n}{\delta^2}\sum_{b\in\mathcal B}\E\biggl[
\int_{\R}\bigl(\widetilde f_b(x)-K_b
\ast f(x)\bigr)^2\,\d x\ind_{B_\varepsilon(\min
\mathcal B)} \biggr]+P \bigl(
\bigl(B_\varepsilon(\min\mathcal B)\bigr)^c \bigr)
\\
&&\quad \lesssim
\frac{(\log n)^2}{\delta^2(n\wedge m)b^{2\beta+1}}+\mathrm{o}(1).
\end{eqnarray*}\upqed
\end{pf}

\begin{pf*}{Proof of Theorem~\ref{thRateUnk}}
A Taylor expansion yields
%
\begin{eqnarray}\label{eqtaylor2}
\widetilde q_{\tau,b}-q_\tau&=&\frac{ {\widetilde M_b(\widetilde q_{\tau
,b})}-\widetilde M_b(q_\tau
)}{\widetilde M_b^{\prime}(q_\tau^*)} =
\frac{ {\widetilde M_b(\widetilde q_{\tau,b})}-\int_{-\infty
}^{q_\tau}\widetilde f_b(x)\,\d x+\tau}{\widetilde f_b(q_\tau
^*)}
\nonumber\\[-8pt]\\[-8pt]
&=&\frac{ {\widetilde M_b(\widetilde q_{\tau,b})}-\int_{-\infty
}^{q_\tau}(\widetilde f_b(x)-f(x))\,\d x}{\widetilde f_b(q_\tau
^*)},\nonumber
\end{eqnarray}
for some intermediate point $q_\tau^*$ between $q_\tau$ and
$\widetilde q_{\tau,b}$. By Proposition~\ref{prDistUnknown} {and
(\ref{eqApprox0})}, the numerator in the above display is of order
$\mathcal{O}_P(n^{-(\alpha+1)/(2\alpha+2\beta+1)})$ for the optimal
bandwidth $b^*$. For the denominator, we will show $\widetilde
f_b(q_\tau^*)=f(q_\tau)+\mathrm{o}_p(1)$ which completes the proof. Since
$f(\bull+q_\tau)\in C^\alpha([-\zeta,\zeta],R)$, we obtain
$|f(x+q_\tau)-f(q_\tau)|<t/2$ for all $|x|\leq(\frac
{t}{2R})^{1\vee\alpha^{-1}}\wedge\zeta=:\delta$ for any $t>0$. Therefore,
%
\begin{eqnarray}\label{eqUniformCons}
&& P\bigl(\bigl|\widetilde f_b\bigl(q_\tau^*
\bigr)-f(q_\tau)\bigr|>t\bigr)\nonumber
\\
&&\quad \leq P \Bigl(\sup_{x\in[-\delta,\delta]}\bigl|
\widetilde f_b(x+q_\tau)-f(q_\tau)\bigr|>t
\Bigr)+P\bigl(|\widetilde q_{\tau,b}-q_\tau|>\delta\bigr)
\\
&&\quad \leq P \Bigl(\sup_{x\in[-\delta,\delta]}\bigl|\widetilde f_b(x+q_\tau
)-f(x+q_\tau)\bigr|>t/2
\Bigr) +P\bigl(|\widetilde q_{\tau,b}-q_\tau|>\delta\bigr).\nonumber
\end{eqnarray}
Checking that the bandwidth satisfies $b\to0$ and $\log
(n)/(nb^{2\beta+1})\to0$ for $n\to\infty$, the first term on the
right-hand side above converges to zero by the uniform consistency
proved in Proposition~\ref{pruniform}. The second one vanishes
asymptotically by Lemma~\ref{lemConsitency}.
\end{pf*}

\begin{pf*}{Proof of Theorem~\ref{thuniformLoss}}
Under the smoothness condition the interval $(\tau_1,\tau_2)$
coincides with a bounded interval of quantiles $(q_{\tau_1},q_{\tau
_2})$. Noting that all our estimates are independent of the quantile,
Theorem~\ref{thuniformLoss} can be proved along the same lines as
Theorem~\ref{thRateUnk} with only minor adaptation to $\sup_{\tau
\in(\tau_1,\tau_2)}$ given a uniform version of Proposition~\ref
{prDistUnknown}:
uniformly over $f$ in the class defined in the theorem and
$f_\varepsilon\in\mathcal D^\beta(R,\gamma)$ for any $b$ such that
$(n\wedge m)b^{2\beta+1}\to\infty$ it holds
%
\begin{eqnarray}\label{equniformDist}
&& \sup_{\tau\in(\tau_1,\tau_2)} \biggl|\int_{-\infty
}^{q_\tau
}
\bigl(\widetilde f_b(x)-f(x)\bigr)\,\d x \biggr|
\nonumber\\[-8pt]\\[-8pt]
&&\quad =\mathcal{O}_P \biggl(b^{\alpha+1}+\biggl(\frac{\log n}{n}\vee
\frac
{1}{m}\biggr)^{1/2}\bigl(b^{-\beta+1/2}\vee1\bigr)
 +\biggl(\frac{1}{n}\vee\frac
{1}{m}\biggr)^{1/2}
\bigl(m^{-1/2}b^{-\beta-1}\vee1\bigr) \biggr).\quad
\nonumber
\end{eqnarray}
Hence, when $\alpha\geq1/2$ the asymptotically optimal choice
$b=(\frac{\log n}{n}\wedge\frac{1}{m})^{1/(2\alpha+2(\beta\vee
1/2)+1)}$ yields
\[
\biggl|\int_{-\infty}^{q_\tau}\bigl(\widetilde
f_b(x)-f(x)\bigr)\,\d x \biggr| =\mathcal{O}_P \biggl(\biggl(
\frac{\log n}{n}\vee\frac{1}{m}\biggr)^{(\alpha
+1)/(2\alpha+2\beta+1)}\vee\biggl(
\frac{\log n}{n}\vee\frac
{1}{m}\biggr)^{1/2} \biggr).
\]

The result (\ref{equniformDist}) can be obtained as Proposition~\ref
{prDistUnknown} except for the term $T_{s,x}=T_{s,x}(q_\tau)$,
defined in (\ref{eqtstx}), which will be treated in the following.
Defining the grid $\tau_1=\sigma_0\leq\cdots\leq\sigma
_M=\tau_2$ such that $q_{\sigma_{k+1}}-q_{\sigma_k}\leq
(q_{\tau_2}-q_{\tau_1})/M$ for $k=1,\dots,M$ and $M\in\N$, we
decompose for any $c>0$
%
\begin{eqnarray}\label{equniformLossDecomp}
P \Bigl(\sup_{\tau\in(\tau_1,\tau_2)}\bigl|T_{s,x}(q_\tau)\bigr|>c
\Bigr) &\leq& P \Bigl(\max_{k=1,\dots,M}\bigl|T_{s,x}(q_{\sigma_k})\bigr|>c/2
\Bigr)
\nonumber\\[-8pt]\\[-8pt]
&&{}+P \Bigl(\mathop{\sup_{q_1,q_2\in(q_{\tau_1},q_{\tau
_2})\dvtx }}_{|q_1-q_2|\leq(q_{\tau_2}-q_{\tau
_1})/(2M)}\bigl|T_{s,x}(q_1)-T_{s,x}(q_2)\bigr|>c/2
\Bigr).
\nonumber
\end{eqnarray}
For the first term, we deduce a concentration inequality. We write
\[
\frac{1}{2\pi}T_{s,x}=\frac{1}{2\pi}T_{s,x}(q_\tau)=
\frac
{1}{n}\sum_{j=1}^n\bigl(
\xi_{j,b}(q_\tau)-\E\bigl[\xi_{j,b}(q_\tau)
\bigr]\bigr)
\]
with
\[
\xi_{j,b}(q_\tau)=\int_{-\infty}^0
a_s(x)\F^{-1} \biggl[\frac{\phi
_K(bu)\mathrm{e}^{\mathrm{i}uY_j}}{\phi_\varepsilon(u)}
\biggr](x+q_\tau)\,\d x=\F^{-1} \biggl[\F a_s(-u)
\frac{\phi_K(bu)\mathrm{e}^{\mathrm{i}uY_j}}{\phi_\varepsilon
(u)} \biggr](q_\tau).
\]
Uniformly in $q_\tau$ we have the deterministic bound
%
\begin{eqnarray}
\label{eqdetBoundXi} \bigl|\xi_{j,b}(q_\tau)\bigr|\leq
\frac{1}{2\pi}\int_{-1/b}^{1/b}\bigl|\F
a_s(-u)\bigr| \biggl|\frac{\phi_K(bu)}{\phi_\varepsilon(u)} \biggr|\,\d u \lesssim\int_{-1/b}^{1/b}
\frac{1}{(1+|u|)|\phi_\varepsilon(u)|}\,\d u \lesssim b^{-\beta}.
\end{eqnarray}
Hence, $|\xi_{j,b}(q_\tau)-\E[\xi_{j,b}(q_\tau)]|\lesssim
b^{-\beta}$. Since the variance of $T_{s,x}(q_\tau)$ is bounded by
(\ref{eqBoundTsX}), Bernstein's inequality (e.g., Massart \cite
{massart2007}, Prop. 2.9) yields for some constant $C>0$ independent of~$q_\tau$
\begin{eqnarray*}
P \bigl(\bigl|T_{s,x}(q_\tau)\bigr|\geq\kappa
\bigl(n^{-1/2}b^{-\beta+1/2}\vee n^{-1/2}\bigr) \bigr)
\leq2\exp\biggl(-\frac{C\kappa^2}{1 + \kappa
(nb)^{-1/2}} \biggr).
\end{eqnarray*}
For the second term on the right-hand side of (\ref
{equniformLossDecomp}), we estimate
\begin{eqnarray*}
\bigl|T_{s,x}(q_1)-T_{s,x}(q_2)\bigr|&
\leq&\biggl\| \biggl(\F^{-1} \biggl[\F a_s(u)
\frac{\phi_K(bu)}{\phi_\varepsilon(u)} \bigl(\phi_n(u)-\phi_Y(u) \bigr)
\biggr] \biggr)' \biggr\|_\infty|q_1-q_2|
\\
&\leq&\frac{|q_1-q_2|}{2\pi}\int_{\R}|u|\bigl|\F
a_s(u)\bigr|\frac
{|\phi_K(bu)|}{|\phi_\varepsilon(u)|} \bigl|\phi_n(u)-
\phi_Y(u) \bigr|\,\d u
\\
&\lesssim&|q_1-q_2|\int_{-1/b}^{1/b}\bigl(1+|u|\bigr)^{\beta}
\bigl|\phi_n(u)-\phi_Y(u) \bigr|\,\d u.
\end{eqnarray*}
Using Markov's inequality, we thus estimate (\ref
{equniformLossDecomp}) by
\begin{eqnarray*}
&&P \Bigl(\sup_{\tau\in(\tau_1,\tau_2)}\bigl|T_{s,x}(q_\tau)\bigr|>
\kappa\bigl(n^{-1/2}b^{-\beta+1/2}\vee n^{-1/2}\bigr) \Bigr)
\\
&&\quad \lesssim M\exp\biggl(-\frac{C\kappa^2}{4 + 2\kappa(nb)^{-1/2}} \biggr)
\\
&&\qquad{} +\frac{(q_{\tau_2}-q_{\tau_1})n^{1/2}(b^{\beta-1/2}\wedge
1)}{M\kappa}\E
\biggl[\int_{-1/b}^{1/b}\bigl(1+|u|\bigr)^{\beta} \bigl|\phi
_n(u)-\phi_Y(u) \bigr|\,\d u \biggr]
\\
&&\quad \lesssim M\exp\biggl(-\frac{C\kappa^2}{4 + 2\kappa(nb)^{-1/2}} \biggr
)+\frac{(q_{\tau_2}-q_{\tau_1})(b^{-3/2}\wedge(b^{-\beta
-1}))}{M\kappa}.
\end{eqnarray*}
Choosing $M=n^{2}$ and $\kappa=(\frac{9}{C}\log n)^{1/2}$, we have
$\kappa(nb)^{-1/2}=\mathrm{o}(1)$ and the previous display converges to zero. Hence,
\[
\sup_{\tau\in(\tau_1,\tau_2)}\bigl|T_{s,x}(q_\tau)\bigr|=\mathcal{O}_P \biggl( \biggl(\frac{\log n}{n} \biggr)^{1/2}b^{-\beta+1/2}
\vee\biggl(\frac
{\log n}{n}\biggr)^{1/2} \biggr).
\]\upqed
\end{pf*}

%
\subsubsection{\texorpdfstring{Proof of Theorem~\protect\ref{thLowerBound}}{Proof of Theorem 2.10}}
To prove the lower bound for the estimation of the distribution
function, we can assume without loss of generality $q=0$. For
$n\leq m$ the estimation error of $\bar F_{n,m}(0)$ is bounded
from below by the estimation error with known error distribution. A
lower bound for the latter is proved by Fan \cite{fan1991} whose
construction can be used in our setting, too.

To\vspace*{1pt} prove the lower bound for $m<n$, we will apply Theorem~2.1 in
Tsybakov \cite{tsybakov2009}. To this end, we construct two
alternatives $(F_i,f_{\varepsilon,i})\in\widetilde{\mathcal
C}^{\alpha+1}(R,r,[-\zeta,\zeta])\times\mathcal D^\beta(R,\gamma
), i=1,2$, such that the $\chi^2$-distance of the corresponding laws
of $(Y_1,\dots,Y_n,\varepsilon^*_1,\dots,\varepsilon^*_m)$ is
bounded by some small constant and such that $|F_1(0)-F_2(0)|$ is
bounded from below with the right rate. Recall that the convolution of
a c.d.f. $F$ with a function $g$ is defined as $F\ast g(x)=\int
g(x-y)\,\d F(y)$. Following the idea by Neumann \cite{neumann1997effect}
our construction will satisfy $F_1\ast f_{\varepsilon,1}=F_2\ast
f_{\varepsilon,2}$ and is thus independent of $n$.

\textit{Step 1}:
For the construction of the alternatives, we need the following: let
$f_0$ be a bounded density whose corresponding distribution is in
$\mathcal C^{\alpha+1}(R,r,\zeta)$ satisfying $q_\tau=0$. Let
$f_{\varepsilon,0}$ be an inner point of $\mathcal D^\beta(R,\gamma
)$ with
%
\begin{eqnarray}
\label{eqPropFe} f_{\varepsilon,0}(x)\gtrsim\bigl(1+|x|\bigr)^{-\gamma-2}, \qquad
\bigl|(\F
f_{\varepsilon,0})^{(k)}(u)\bigr|\lesssim\bigl(1+ |u|\bigr)^{-\beta},\qquad k=0,\dots,K
\end{eqnarray}
for $x,u\in\R$ and an integer $K>\gamma/2+1$. Let the perturbation
$g\in C^\infty(\R)\cap L^1(\R)$ satisfy
\begin{eqnarray*}
\int_{\R}g(x)\,\d x&=&0,\qquad \int_{-\infty}^0
g(x)\,\d x\neq0,
\\
\bigl\|\bigl(1\vee x^{\gamma\vee1} \bigr)g(x)\bigr\|_{L^1}&<&
\infty,\qquad \supp\F g\subset[-2,-1]\cup[1,2].
\end{eqnarray*}
Define $g_b:=b^{-1}g(\bull/b)$ for $b>0$ and for some $a\in(0,1),c>0$
%
\begin{eqnarray}\label{eqalternatives}
F_1(x)&:=&a\int_{-\infty}^x
f_0(y)\,\d y+(1-a)\ind_{[2\zeta,\infty
)}(x),\nonumber
\\
f_{\varepsilon,1}(x)&:=&f_{\varepsilon,0}+cb^{\alpha
+1}\bigl(f_{\varepsilon,0}\ast g_b(\bull+2\zeta)\bigr) (x),
\nonumber\\[-8pt]\\[-8pt]
F_2(x)&:=&F_1(x)+cb^{\alpha+1}\int
_{-\infty}^xg_b(\bull+2\zeta)\ast
F_1(y)\,\d y,\nonumber
\\
f_{\varepsilon,2}(x)&:=&f_{\varepsilon,0}(x).\nonumber
\end{eqnarray}
Owing to $\int g_b=0$, $F_i$ are distribution functions admitting
Lebesgue densities on $[-\zeta,\zeta]$ which are at least $\alpha
$-H\"older continuous.\vspace*{2pt} Estimating $\|f_0\ast g_b\|_{C^\alpha(\R
)}\lesssim\|f_0\|_{L^1}\|g_b\|_{C^\alpha(\R)}\lesssim b^{-\alpha
-1}$, we infer that $\,\d F_2$ is contained in a closed H\"older ball.
Hence, $F_i\in\widetilde{\mathcal C}^{\alpha+1}(R,r,[-\zeta,\zeta
])$ for $c>0$ sufficiently small. $f_{\varepsilon_i}\in\mathcal
D^\beta(R,\gamma)$ can be verified, using $\int g=0,\|\F g\|_\infty
\leq\|g\|_{L^1}$ and $\|(\F g)'(u)(1+|u|)\|_\infty<\infty$.

\textit{Step 2}: To bound the distance $|F_1(0)-F_2(0)|$ from below we
note, using Fubini's theorem, $\int g=0$ and $\|f_0\|_\infty<\infty$,
%
\begin{eqnarray}
\label{eqdistance} F_2(0)-F_1(0)&=&b^{\alpha+1} \biggl(ac
\int_{\R}\int_{2\zeta
}^{-y+2\zeta}f_0(x)
g_b(y)\,\d x\,\d y+(1-a)c\int_{-\infty}^{0}g_b(x)
\,\d x \biggr)
\nonumber
\\
&=&b^{\alpha+1} \biggl((1-a)c\int_{-\infty}^{0}g(x)
\,\d x+\mathcal{O} \bigl(\bigl\|yg_b(y)\bigr\|_{L^1} \bigr) \biggr)
\\
&=&b^{\alpha+1} \biggl((1-a)c\int_{-\infty}^0g(x)
\,\d x+\mathcal{O}(b) \biggr),\nonumber
\end{eqnarray}
for $b$ small enough. Therefore, $|F_1(0)-F_2(0)|\gtrsim b^{\alpha+1}$.

\textit{Step 3}: Using the independence of the observations, the
sample $(Y_1,\dots,Y_n,\varepsilon^*_1,\dots,\varepsilon^*_m)$ is
distributed according to $(F_i\ast f_{\varepsilon,i})^{\otimes
n}\otimes f_{\varepsilon,i}^{\otimes m}$ under the hypotheses $i=1,2$.
By construction $F_1\ast f_{\varepsilon,1}=F_2\ast f_{\varepsilon,2}$
such that the $\chi^2$-distance of the laws of the observations equals
%
\begin{equation}
\label{eqChi2dist} \chi^2\bigl(f_{\varepsilon,1}^{\otimes
m},f_{\varepsilon,2}^{\otimes
m}
\bigr)= \biggl(1+\int_{\R}\frac{(f_{\varepsilon,1}-f_{\varepsilon
,2})^2(x)}{f_{\varepsilon,2}(x)}\,\d x
\biggr)^m-1.
\end{equation}
We decompose
\begin{eqnarray*}
&&\int_{\R}\frac{(f_{\varepsilon,1}-f_{\varepsilon
,2})^2(x)}{f_{\varepsilon,2}(x)}\,\d x
\\
&&\quad =c^2b^{2\alpha+2} \biggl(\int_{|x|\leq1}
\frac{(f_{\varepsilon,0}\ast g_b(\bull+2\zeta))^2(x)}{f_{\varepsilon
,0}(x)}\,\d x+\int_{|x|>1}\frac{(f_{\varepsilon,0}\ast g_b(\bull+2\zeta
))^2(x)}{f_{\varepsilon,0}(x)}\,\d x
\biggr)
\\
&&\quad =:c^2b^{2\alpha+2} (I_1+I_2 ).
\end{eqnarray*}
For the first integral, we use $\inf_{|x|\leq1}f_{\varepsilon
,0}(x)>0$, Plancherel's identity, $f_{\varepsilon,0}\in\mathcal
D^\beta(R,\gamma)$ and the support of $\F g$ to estimate
\begin{eqnarray*}
|I_1|\lesssim\int_{\R}\bigl|\F f_{\varepsilon,0}(u)
\F g(bu)\mathrm{e}^{-\mathrm{i}2\zeta
u}\bigr|^2\,\d u \lesssim\int_{1/b\leq|u|\leq2/b}\bigl(1+|u|\bigr)^{-2\beta}
\,\d u\lesssim b^{2\beta-1}.
\end{eqnarray*}
Using (\ref{eqPropFe}), $I_2$ can be estimated similarly
\begin{eqnarray*}
|I_2|&\lesssim&\int_{|x|>1}\bigl(1+|x|\bigr)^{\gamma+2}|x|^{-2K}\bigl|
\F^{-1}\bigl[\bigl(\F f_{\varepsilon,0}\F g_b\mathrm{e}^{-\mathrm{i}2\zeta\bull}
\bigr)^{(K)}\bigr]\bigr|^2(x)\,\d x
\\
&\sim&\int_{1/b\leq|u|\leq2/b}\bigl|\bigl(\F f_{\varepsilon,0}(u)\F
g(bu)\mathrm{e}^{-\mathrm{i}2\zeta u}\bigr)^{(K)}\bigr|^2\,\d u \lesssim
b^{2\beta-1}.
\end{eqnarray*}
We conclude from (\ref{eqChi2dist}) for some constant $C>0$ that
\[
\chi^2\bigl(f_{\varepsilon,1}^{\otimes m},f_{\varepsilon,2}^{\otimes m}
\bigr) \leq\bigl(1+Cc^{2}b^{2\alpha+2\beta+1}\bigr)^m-1
\leq\exp\bigl(Cc^{2}mb^{2\alpha+2\beta+1}\bigr)-1,
\]
which can be bounded by an arbitrarily small constant if $c$ is chosen
sufficiently small and $b=m^{-1/(2\alpha+2\beta+1)}$. We obtain from
Step~2 that $|F_1(0)-F_2(0)|\geq C m^{-(\alpha+1)/(2\alpha+2\beta
+1)}$, for some positive constant $C$.

\textit{Step 4}: Replacing in (\ref{eqalternatives}) the factor
$b^{\alpha+1}$ in $F_2$ and $f_{\varepsilon,1}$ by $cm^{-1/2}$ for
some sufficiently small constant $c>0$ and choosing $b=1$, the previous
steps yield the lower bound $m^{-1/2}$.

Let us finally conclude the lower bound for the estimation error of the
quantiles. We use the construction from Step~1, denoting the $\tau
$-quantile of $F_i$ by $q_{\tau,i}$. We note $|q_{\tau,1}|<\delta$
for any $\delta>0$ if we choose $a$ close enough to one and thus $F_1$
is regular in an interval around $q_{\tau,1}$. Moreover, it holds
\begin{eqnarray*}
\|F_1-F_2\|_\infty&\leq& c
\bigl(m^{-1/2}\vee b^{\alpha+1}\bigr) \bigl\|\bigl(a f_0+(1-a)
\delta_{-2\zeta}\bigr)\ast g_b(\bull+2\zeta)
\bigr\|_{L^1}
\\
&\leq& c\bigl(m^{-1/2}\vee b^{\alpha+1}\bigr)\|g
\|_{L^1}\to0.
\end{eqnarray*}
We infer analogously to (\ref{eqConsistencyBound}) that $|q_{\tau
,1}-q_{\tau,2}|<\delta$ for any $\delta>0$ and $m$ sufficiently
large implying $F_i\in\widetilde{\mathcal C}^{\alpha+1}(R,r,\zeta)$.
Applying a Taylor expansion similar to (\ref{eqtaylor}), we obtain
\[
q_{\tau,2}-q_{\tau,1}=-\frac{F_2(q_{\tau,1})-F_1(q_{\tau
,1})}{F_2'(q_\tau^*)}
\]
for some intermediate point between $q_{\tau,1}$ and $q_{\tau,2}$.
The denominator $F_2'(q_\tau^*)$ is bounded from above and below owing
to $\sup_{|x|\leq\zeta}|F_2'(x)-af_0(x)|\to0$, $|q_{\tau
,2}|\leq|q_{\tau,2}-q_{\tau,1}|+|q_{\tau,1}|<2\delta$ and
$f_0(0)>0$. (\ref{eqdistance}) yields $|q_{\tau,2}-q_{\tau,1}|\gtrsim
m^{-1/2}\vee b^{\alpha+1}$. The assertion follows from
Steps~3 and 4 above.

\subsection{\texorpdfstring{Proofs for Section~\protect\ref{secadaptive}}{Proofs for Section 3}}
We start with Lemma~\ref{lemBmin} concerning the bandwidth set $\B
_n$ from (\ref{eqbandwidthset}).
\subsubsection{\texorpdfstring{Proof of Lemma~\protect\ref{lemBmin}}{Proof of Lemma 3.1}}

By Lemma~\ref{lemPhiEps}, we can argue on the event $B_\varepsilon
(b)$ from (\ref{eqBEps}). The deterministic counterpart of
$\widetilde j_n$, defined in (\ref{eqtildeJ}), is given by
%
\begin{equation}
\label{eqj0} j_{0,n}:=\min\biggl\{j=0,\dots,N_n\dvt 2
\leq\biggl(\frac{\log
n}{n} \biggr)^{1/2}\int
_{-1/b_{j}}^{1/b_{j}}\bigl|\phi_{\varepsilon
}(u)\bigr|^{-1}
\,\d u\leq4 \biggr\}.
\end{equation}
Noting that for $f_\varepsilon\in\mathcal D^\beta(R,\gamma)$
\[
4\geq\biggl(\frac{\log n}{n} \biggr)^{1/2}\int
_{-1/b_{j_{0,n}}}^{1/b_{j_{0,n}}}\bigl|\phi_{\varepsilon}(u)\bigr|^{-1}
\,\d u \gtrsim\biggl(\frac{\log n}{nb_{j_{0,n}}^{2\beta+2}} \biggr)^{1/2}
\]
we obtain $nb_{j_{0,n}}^{2\beta+2}\to\infty$ and thus it is
sufficient to prove
%
\begin{equation}
\label{eqBminProb} \inf_{f\in\mathcal C^\alpha(R,r,\zeta)}\inf
_{f_\varepsilon\in
\mathcal D^\beta(R,\gamma)} P
\bigl(\bigl\{b_{j_{0,n}}\leq b_{\widetilde j_{n}}\leq b^*\bigr
\}\cap
B_\varepsilon(b_{j_{0,n}}) \bigr)\to1\qquad\mbox{as } n\to\infty,
\end{equation}
for the optimal bandwidth $b^*=n^{-1/(2\alpha+2(\beta\vee1/2)+1)}$.
For convenience, we define
\[
I_n(b):= \biggl(\frac{\log n}{n} \biggr)^{1/2}\int
_{-1/b}^{1/b}\frac{\d
u}{|\phi_{\varepsilon}(u)|},\qquad\widetilde
I_n(b):= \biggl(\frac{\log n}{n} \biggr)^{1/2}\int
_{-1/b}^{1/b}\frac{\d u}{|\phi_{\varepsilon,m}(u)|}.
\]
Assume $b_{\widetilde j_{n}}<b_{j_{0,n}}$, then monotonicity implies
$\widetilde I_n(b_{j_{0,n}})\leq\widetilde I_n(b_{\widetilde
j_n})\leq1$. Combined with $I_n(j_{0,n})\geq2$, we obtain
$I_n(b_{j_{0,n}})-\widetilde I_n(b_{j_{0,n}})\geq1$. Hence,
%
\begin{equation}
\label{eqbmin1} \{b_{\widetilde j_{n}}<b_{j_{0,n}} \}\subset\bigl\{
\bigl|I_n(b_{j_{0,n}})-\widetilde I_n(b_{j_{0,n}})\bigr|
\geq1 \bigr\}.
\end{equation}
On the other hand, if $b^*<b_{\widetilde j_{n}}$, we get $\widetilde
I_n(b^*)\geq\widetilde I_n(b_{\widetilde j_n})\geq1/2$.
Since $I_n(b^*)\lesssim(\frac{\log n}{n(b^*)^{2\beta+2}})^{1/2}$
converges to zero, $I_n(b^*)\leq1/4$ for $n$ large enough. Thus,
%
\begin{equation}
\label{eqbmin2} \bigl\{b_{\widetilde j_n}>b^* \bigr\}\subset\bigl
\{\bigl|I_n\bigl(b^*\bigr)-\widetilde I_n\bigl(b^*\bigr)\bigr|
\geq1/4 \bigr\}.
\end{equation}
To show that the probabilities of the right-hand sides of (\ref
{eqbmin1}) and (\ref{eqbmin2}) converge to zero, we first apply the
Cauchy--Schwarz inequality
\begin{eqnarray*}
\bigl|I_n(b)-\widetilde I_n(b) \bigr|^2 &\leq&
\frac{\log n}{n}\int_{-1/b}^{1/b}\frac{\d u}{|\phi
_\varepsilon(u)|^2}
\int_{-1/b}^{1/b} \biggl|\frac{\phi_\varepsilon
(u)}{\phi_{\varepsilon,m}(u)}-1
\biggr|^2\,\d u
\\
&\lesssim& \frac{\log n}{nb^{2\beta+1}}\int_{-1/b}^{1/b}
\biggl|\frac
{\phi_\varepsilon(u)}{\phi_{\varepsilon,m}(u)}-1 \biggr|^2\,\d u.
\end{eqnarray*}
Markov's inequality and (\ref{eqNeumann}) yield for $b\in\{
b_{\min},b^*\}$
\begin{eqnarray*}
P \biggl( \biggl\{ \bigl|I_n(b)-\widetilde I_n(b) \bigr|\geq
\frac
{1}{4} \biggr\}\cap B_\varepsilon(b_{j_{0,n}}) \biggr)
& \lesssim& \frac{\log n}{nb^{2\beta+1}}\int_{-1/b}^{1/b}\E\biggl[
\biggl|\frac{\phi_\varepsilon(u)}{\phi_{\varepsilon,m}(u)}-1 \biggr|^2\ind
_{B_\varepsilon(b_{j_{0,n}})} \biggr]\,\d u
\\
&\lesssim&
\frac{\log n}{nmb^{4\beta+2}}
\end{eqnarray*}
which converges to zero. Therefore, (\ref{eqBminProb}) holds
true.

\subsubsection{Preparations to the Proof of Theorem \texorpdfstring{\protect\ref{thoracleinequality}}{3.2}}
Before we can prove Theorem~\ref{thoracleinequality}, some
preparations are needed. By Lemma~\ref{lembias} there is a constant
$D>0$ such that the bias can be bounded by $B_b:=Db^{\alpha+1}$. By
the error representation~(\ref{eqtaylor2}), we have for any $b\in\B$
%
\begin{eqnarray}
\label{eqerrordecompadapt} |\widetilde q_{\tau,b}-q_\tau|&=& \biggl|
\frac{\int_{-\infty}^{q_\tau
}(\widetilde f_b(x)-f(x))\,\d x-\widetilde M_{b}(\widetilde q_{\tau
,b})}{\widetilde f_b(\widetilde q^*)} \biggr|
\nonumber\\[-8pt]\\[-8pt]
& \leq& \frac
{B_b+|V_{b,X}+V_{b,\varepsilon}+V_{b,c}| +|\widetilde
M_{b}(\widetilde q_{\tau,b})|}{|\widetilde{f}_b(q^*)|}\nonumber
\end{eqnarray}
with some $q^*\in[(q_\tau\wedge\widetilde q_{\tau,b}), (q_\tau\vee
\widetilde q_{\tau,b})]$ and where the stochastic error is decomposed in
\begin{eqnarray}
\label{eqVbX} V_{b,X}&:=&\frac{1}{n}\sum
_{j=1}^n \bigl(\xi_j(b)-\E\bigl[
\xi_j(b)\bigr] \bigr)\qquad\mbox{with}\nonumber
\\
\xi_j(b)&:=&\int_{-\infty}^0a_s(x)
\F^{-1} \biggl[\frac{\phi
_K(bu)\mathrm{e}^{\mathrm{i}uY_j}}{\phi_\varepsilon(u)} \biggr](x+q_{\tau})\,\d x,
\nonumber\\[-8pt]\\[-8pt]
V_{b,\varepsilon}&:=&\int_{-\infty}^0a_s(x)
\F^{-1} \biggl[\frac{\phi
_K(bu)\phi_n(u)}{\phi_\varepsilon(u)} \biggl(\frac{\phi_{\varepsilon
}(u)}{\phi_{\varepsilon,m}(u)}-1 \biggr)
\biggr](x+q_{\tau})\,\d x,
\nonumber
\\
V_{b,c}&:=&\int_{-\infty}^0
a_c(x)\F^{-1} \biggl[\phi_K(bu) \biggl(
\frac
{\phi_n(u)}{\phi_{\varepsilon,m}(u)}-\phi_X(u) \biggr) \biggr](x+q_\tau)
\,\d x.
\nonumber
\end{eqnarray}
In view of the analysis in Section~\ref{secDistUnknown}, the part of
the stochastic error which is due to the continuous part $a_c$ will be
negligible. Hence, we concentrate on $V_{b,X}$ and $V_{b,\varepsilon
}$. By independence of $(\xi_j(b))_j$, we obtain
%
\begin{eqnarray}
\label{eqVxVar} \Var(V_{b,X})&\leq&\frac{1}{n}\E\bigl[
\xi_j(b)^2\bigr] =\frac{1}{n}\E\biggl[ \biggl(
\int_{-\infty}^0a_s(x)\F^{-1}
\biggl[\frac{\phi_K(bu)\mathrm{e}^{\mathrm{i}uY_j}}{\phi
_\varepsilon(u)} \biggr](x+q_\tau)\,\d x \biggr)^2
\biggr]
\nonumber\\[-8pt]\\[-8pt]
&=:&\sigma^2_{b,X}.\nonumber
\end{eqnarray}
We will determine the variance of $V_{b,\varepsilon}$ on the event
$B_\varepsilon(b)$, defined in (\ref{eqBEps}). We apply Plancherel's
identity and the Cauchy--Schwarz inequality to separate $Y_i$ and
$\varepsilon_i$ from each other:
%
\begin{eqnarray}\label{eqVbEpsDecomp}
&&\E\bigl[|V_{b,\varepsilon}|\ind_{B_\varepsilon(b)}\bigr]\nonumber
\\
&&\quad =\frac{1}{2\pi}\E\biggl
[ \biggl|\int_{\R} \F a_s(-u)\mathrm{e}^{-\mathrm{i}uq_\tau
}\frac{\phi_K(bu)\phi_n(u)}{\phi_\varepsilon(u)} \biggl(\frac{\phi
_\varepsilon(u)}{\phi_{\varepsilon,m}(u)}-1 \biggr)\,\d u \biggr|\ind
_{B_\varepsilon(b)} \biggr]
\nonumber
\\
&&\quad \leq \frac{1}{2\pi}\E\biggl[ \biggl(\int_{\R} \bigl|
\phi_K(bu)\bigr| \biggl|\frac{\phi_n(u)}{\phi_\varepsilon(u)} \biggr|^2\,\d u
\biggr)^{1/2}
\nonumber\\[-8pt]\\[-8pt]
&& \hspace*{44pt}{}\times  \biggl(\int_{\R}\bigl|\phi_K(bu)\bigr|\bigl|
\F a_s(-u)\bigr|^2 \biggl|\frac{\phi_\varepsilon
(u)}{\phi_{\varepsilon,m}(u)}-1 \biggr|^2\,\d u
\biggr)^{1/2}\ind_{B_\varepsilon(b)} \biggr]
\nonumber
\\
&&\quad \leq \frac{1}{2\pi}\E\biggl[ \biggl(\int_{\R}\bigl|
\phi_K(bu)\bigr| \biggl|\frac{\phi_n(u)}{\phi_\varepsilon(u)} \biggr|^2\,\d u
\biggr)^{1/2} \biggl(\int_{\R}\bigl|\phi_K(bu)\bigr|
\biggl|\frac{\F a_s(-u)}{\phi_{\varepsilon,m}(u)} \biggr|^2\,\d u\ind_{B_\varepsilon(b)}
\biggr)^{1/2}
\nonumber
\\
&&\hspace*{44pt}{} \times\sup_{|u|\leq1/b}\bigl|\phi_{\varepsilon,m}(u)-\phi
_{\varepsilon}(u)\bigr| \biggr].\nonumber
\end{eqnarray}
Let us define
%
\begin{eqnarray}
\label{eqsigmaEps} \sigma_{b,\varepsilon}&:=&\frac{1}{2\pi}m^{-1/2}\sigma
_{b,\varepsilon,1}\sigma_{b,\varepsilon,2}
\end{eqnarray}
with
\begin{eqnarray*}
\sigma_{b,\varepsilon,1}&:=&\E\biggl[ \biggl(\int_{\R}\bigl|\phi
_K(bu)\bigr| \biggl|\frac{\phi_n(u)}{\phi_\varepsilon(u)} \biggr|^2\,\d u
\biggr)^{1/2} \biggr],
\\
\sigma_{b,\varepsilon,2}&:=&\E\biggl[ \biggl(
\int_{\R}\bigl|\phi_K(bu)\bigr| \biggl|\frac{\F a_s(-u)}{\phi_{\varepsilon,m}(u)}
\biggr|^2\,\d u \biggr)^{1/2}\ind_{B_\varepsilon(b)} \biggr].
\end{eqnarray*}
With the bounds $\sigma_{b,X}$ and $\sigma_{b,\varepsilon}$ at hand,
we obtain the following concentration results.

\begin{lemma}\label{lemconcIneq}
Let $\mathcal B$ be a set satisfying $|\mathcal B|\lesssim\log n,
(\log\log n)/nb_1\to0$ for $b_1=\min\mathcal B$ as well as $|\log
b_1|\lesssim\log n$. Then we obtain uniformly over $f\in\mathcal
C^\alpha(R,r,\zeta)$ and $f_\varepsilon\in\mathcal D^\beta
(R,\gamma)$ for any $\delta>0$:
\begin{enumerate}[(iii)]
\item[(i)]$
P (\exists b\in\mathcal B\dvt |V_{b,X}|\geq(1+\delta)\sqrt
{\log\log n} (\sqrt{2}\sigma_{b,X}+\mathrm{o}(n^{-1/2}(b^{-\beta
+1/2}\vee1)) ) )\to0$.\vspace*{3pt}
\item[(ii)]$
P (\exists b\in\mathcal B\dvt  |V_{b,\varepsilon}|\geq
\delta(\log n)^3\sigma_{b,\varepsilon} )\to0$.\vspace*{3pt}
\item[(iii)] Assuming further $mb_1^{(2\beta\wedge1)+2}\gtrsim1$,\\[3pt]
$
P (\exists b\in\mathcal B\dvt  |V_{b,c}|\geq(\log
n)^{3/2}n^{-1/2}(b^{-\beta+1/2}\vee1) )\to0$.
\end{enumerate}
\end{lemma}

\begin{pf}
\textup{(i})
Using the deterministic bound (\ref{eqdetBoundXi}), we obtain $|\xi
_j(b)-\E[\xi_j(b)]|\leq Cb^{-\beta}$ for some constant $C>0$.
Since the variance is bounded by (\ref{eqVxVar}), Bernstein's
inequality (e.g., Massart \cite{massart2007}, Prop. 2.9) yields for
any positive $\kappa_n=\mathrm{o}(nb)$
\begin{eqnarray*}
P \biggl(|V_{b,X}|\geq\sqrt{2\sigma_{b,X}^2
\kappa_n}+\frac
{C\kappa_n}{3nb^\beta} \biggr) &\leq&2\mathrm{e}^{-\kappa_n}.
\end{eqnarray*}
Hence, $\sqrt{\kappa_n}(nb^\beta)^{-1}\lesssim(n(b^{2\beta
-1}\wedge1))^{-1/2} (\kappa_n/(nb) )^{1/2}$ yields uniformly
in $\mathcal C^\alpha(R,r,\zeta)$ and $\mathcal D^\beta(R,\gamma)$
\[
P \bigl(|V_{b,X}|\geq\sqrt{\kappa_n} \bigl(\sqrt{2}
\sigma_{b,X}+\mathrm{o}\bigl(n^{-1/2}\bigl(b^{-\beta+1/2}\vee1\bigr)
\bigr) \bigr) \bigr)\leq2\mathrm{e}^{-\kappa_n}.
\]
The result follows from choosing $\kappa=(1+\delta)^2\log\log n$ and
using $|\mathcal B|\lesssim\log n$.

(ii) Using an estimate as in (\ref{eqVbEpsDecomp}), we obtain
\begin{eqnarray*}
|V_{b,\varepsilon}| &\leq&\frac{1}{2\pi} \biggl(\int
_{\R}\bigl|\phi_K(bu)\bigr| \biggl|\frac
{\phi_n(u)}{\phi_\varepsilon(u)}
\biggr|^2\,\d u \biggr)^{1/2} \biggl(\int_{\R}\bigl|
\phi_K(bu)\bigr|\bigl|\F a_s(-u)\bigr|^2 \biggl|
\frac{\phi_\varepsilon
(u)-\phi_{\varepsilon,m}(u)}{\phi_{\varepsilon,m}(u)} \biggr|^2\,\d u \biggr)^{1/2}
\\
&\leq&\frac{1}{2\pi}\underbrace{ \biggl(\int_{\R}\bigl|
\phi_K(bu)\bigr| \biggl|\frac{\phi_n(u)}{\phi_\varepsilon(u)} \biggr|^2\,\d u
\biggr)^{1/2}}_{=:V_{b,\varepsilon,1}}
\\
&&{}\times \underbrace{ \biggl(\int
_{\R}\bigl|\phi_K(bu)\bigr| \biggl|\frac{\F a_s(-u)}{\phi_{\varepsilon,m}(u)}
\biggr|^2\,\d u \biggr)^{1/2}}_{=:V_{b,\varepsilon,2}}\sup
_{|u|\leq1/b} \bigl|\phi_\varepsilon(u)-\phi_{\varepsilon,m}(u) \bigr|.
\end{eqnarray*}
Hence, for any $c\in(0,1/4)$
\begin{eqnarray*}
&&P \bigl( \bigl\{ |V_{b,\varepsilon}|\geq\delta(\log n)^3
\sigma_{b,\varepsilon} \bigr\}\cap B_\varepsilon(b) \bigr)
\\
&&\quad\leq P \bigl(|V_{b,\varepsilon,1}|\geq(\log n)^{1+c}
\sigma_{b,\varepsilon,1} \bigr) +P \bigl( \bigl\{|V_{b,\varepsilon
,2}|\geq(\log
n)^{1+c}\sigma_{b,\varepsilon,2} \bigr\}\cap B_\varepsilon(b) \bigr)
\\
&&\qquad{}+P \Bigl(\sup_{|u|\leq1/b} \bigl|\phi_\varepsilon(u)-\phi
_{\varepsilon,m}(u) \bigr|\geq\delta(\log n)^{1-2c} m^{-1/2}
\Bigr)
\\
&&\quad =:P_{b,1}+P_{b,2}+P_{b,3}.
\end{eqnarray*}
The first two probabilities can be bounded by Markov's inequality:
\begin{eqnarray*}
P_{b,1}&\leq&(\log n)^{-1-c}\sigma_{b,\varepsilon,1}^{-1}
\E[V_{b,\varepsilon,1}]=(\log n)^{-1-c},
\\
P_{b,2}&\leq&(\log n)^{-1-c}\sigma_{b,\varepsilon,2}^{-1}
\E[V_{b,\varepsilon,2}\ind_{B_\varepsilon(b)}]=(\log n)^{-1-c}.
\end{eqnarray*}
For $P_{b,3}$ we will apply the following version of Talagrand's
inequality (cf. Massart \cite{massart2007}, (5.50)): let $T$ be a
countable index and for all $t\in T$ let $Z_{1,t},\dots, Z_{n,t}$ be
an i.i.d. sample of centered, complex valued random variables
satisfying $\|Z_{k,t}\|_\infty\leq b$, for all $t\in T,k=1,\dots
,n$, as well as $\sup_{t\in T}\Var(\sum_{k=1}^nZ_{k,t})\leq
v<\infty$. Then for all $\kappa>0$
%
\begin{eqnarray}
\label{eqTalagrand} P \Biggl(\sup_{t\in T} \Biggl|\sum
_{k=1}^n Z_{k,t} \Biggr|\geq4\E\Biggl[
\sup_{t\in T} \Biggl|\sum_{k=1}^n
Z_{k,t} \Biggr| \Biggr]+\sqrt{2v\kappa}+\frac{2}{3}b\kappa\Biggr)
\leq2\mathrm{e}^{-\kappa}.
\end{eqnarray}
Choosing the rational numbers $T=\Q\cap[-\frac{1}{b},\frac{1}{b}]$
and $Z_{k,t}:=\mathrm{e}^{\mathrm{i}t\varepsilon^*_k}-\phi_\varepsilon(t)$,
Talagrand's inequality applies with $b=2$ and $v=n$. As in (\ref
{eqboundNeumannReiss}), we use Theorem 4.1 by Neumann and Rei{\ss}
\cite{neumann2009nonparametric} to obtain for any $\eta\in(0,1/2)$
\[
m^{1/2}\E\Bigl[{\sup_{|u|\leq1/b}}\bigl|\phi_{\varepsilon,m}(t)-
\phi_\varepsilon(t)\bigr| \Bigr]\lesssim|\log b|^{1/2+\eta}.
\]
Therefore on the assumptions $\kappa_n^{-1}(\log n)^{1+2\eta}\to0$
and $\kappa_n/m\to0$
\begin{eqnarray*}
4\E\Bigl[{\sup_{|u|\leq1/b,u\in\Q}}\bigl|\phi_{\varepsilon,m}(u)-
\phi_\varepsilon(u)\bigr| \Bigr]+\sqrt{\frac{2\kappa_n}{m}}+\frac
{4}{3m}
\kappa_n=\sqrt{\frac{\kappa_n}{m}} \bigl(\sqrt{2}+\mathrm{o}(1) \bigr)
\end{eqnarray*}
and thus continuity of $\phi_{\varepsilon,m}$ and (\ref
{eqTalagrand}) yield
%
\begin{eqnarray}
P_{b,3}=P \Bigl(\sup_{|u|\leq1/b,u\in\Q} \bigl|\phi
_{\varepsilon,m}(u)-\phi_{\varepsilon}(u) \bigr|\geq\bigl(\sqrt
2+\mathrm{o}(1)\bigr)
\sqrt{\kappa_n/m} \Bigr) \leq2\mathrm{e}^{-\kappa_n}\label{eqTal2Phi}.
\end{eqnarray}
With $\kappa_n=\frac{\delta}{2}(\log n)^{2-4c}$ for $c<1/4-\eta/2$,
we obtain $P_3\leq2n^{-\delta/2}$. Using $b_1=\min\mathcal B,
|B|\lesssim\log n$ and Lemma~\ref{lemPhiEps}, we finally get
\begin{eqnarray*}
P \Bigl(\sup_{b\in\mathcal B} |V_{b,\varepsilon}|\geq(\sqrt2+
\delta) (\log n)^3\sigma_{b,\varepsilon} \Bigr) &\leq&\sum
_{b\in\mathcal B} (P_{b,1}+P_{b,2}+P_{b,3}
)+P\bigl(B_\varepsilon(b_1)^c\bigr)=\mathrm{o}(1).
\end{eqnarray*}

(iii) Corollary~\ref{coruniformQ} shows for $\delta_b>0$
and for any sequence $(x_n)_n$ that tends to infinity
\[
P \bigl(\exists b\in\mathcal B\dvt  |V_{b,c}|\geq\delta_b \bigr)
\lesssim\sum_{b\in\mathcal B}\frac{x_n}{\delta_b^2n(mb^{2\beta
+2}\wedge1)}+\mathrm{o}(1).
\]
Choosing $\delta_b=(\log n)^{3/2}n^{-1/2}(b^{-\beta+1/2}\vee1)$ and
$x_n=\mathrm{o}((\log n)^{1/2})$ yields
\begin{eqnarray*}
&&P \bigl(\exists b\in\mathcal B\dvt  |V_{b,c}|\geq(\log
n)^{3/2}n^{-1/2}\bigl(b^{-\beta+1/2}\vee1\bigr) \bigr)
\\
&&\quad\lesssim\sum_{b\in\mathcal B}\frac{x_n}{(\log
n)^{3}(mb^{(2\beta\wedge1)+2}\wedge1)}+\mathrm{o}(1)
\lesssim\frac{x_n}{(\log n)^{2}(mb^{(2\beta\wedge1)+2}\wedge1)}+\mathrm{o}(1)=\mathrm{o}(1).
\end{eqnarray*}\upqed
\end{pf}

For the denominator in the error representation~(\ref
{eqerrordecompadapt}) we need uniform consistency. A uniform result
on the error $|\widetilde q_{\tau,b}-q_\tau|$ follows immediately.

\begin{lemma}
Let $\mathcal B$ be a finite set satisfying $|\mathcal B|\lesssim\log
n$, $\sup_{b\in\mathcal B}b\log(n)\to0$ as well as $\sup_{b\in
\mathcal B}(\log n)^2/(nb^{2\beta+1})\to0$. Then we obtain for $n\to
\infty$ and $\eta\in(0,1)$
%
\begin{equation}
\label{eqfDiffUniform} \sup_{f\in\mathcal C^\alpha(R,r,\zeta,U_n)}\sup
_{f_\varepsilon\in
\mathcal D^\beta(R,\gamma)} P
\Bigl(\sup_{b\in\mathcal B}\sup_{q^*_\tau\in[q_\tau\wedge\widetilde
q_{\tau,b},q_\tau\vee
\widetilde q_{\tau,b}]}\bigl|\widetilde
f_b\bigl(q^*_\tau\bigr)-f(q_\tau)\bigr|>\eta
f(q_\tau) \Bigr)\to0.
\end{equation}
Moreover, supposing $\min_{b\in\mathcal B}nb^{(2\beta\wedge
1)+2}\gtrsim1$, we obtain uniformly in $f\in\mathcal C^\alpha
(R,r,\zeta)$ and $f_\varepsilon\in\mathcal D^\beta(R,\gamma)$ for
any sequence of critical values $(\delta_b)_{b\in\mathcal B}$
satisfying $\inf_{\mathcal B}\delta_b\to\infty$
%
\begin{equation}
\label{eqqDiffUniform} P \bigl(\exists b\in\mathcal B\dvt |\widetilde
q_{\tau,b}-q_\tau|>
\delta_b \bigl(3Db^{\alpha+1}+n^{-1/2}
\bigl(b^{-\beta+1/2}\vee1\bigr) \bigr) \bigr) \lesssim\sum
_{b\in\mathcal B}\frac{1}{\delta_b}+\mathrm{o}(1).
\end{equation}
\end{lemma}

\begin{pf}
Since $f(q_\tau)\geq r$ and $f\in C^\alpha([q_\tau-\zeta,q_\tau
+\zeta],R)$, decomposition~(\ref{eqUniformCons}) implies
with $\kappa=(\frac{\eta r}{2R})^{1\vee\alpha^{-1}}\wedge\zeta$
%
\begin{eqnarray}\label{eqdecomp1}
&& P \Bigl(\sup_{b\in\mathcal B}\sup_{q^*_\tau\in[q_\tau\wedge
\widetilde q_{\tau,b},q_\tau\vee\widetilde q_{\tau,b}]}\bigl|\widetilde
f_b\bigl(q^*_\tau\bigr)-f(q_\tau)\bigr|>\eta
f(q_\tau) \Bigr)
\nonumber\\[-8pt]\\[-8pt]
&&\quad \leq P \Bigl(\sup_{b\in\mathcal B}
\sup_{x\in[-\kappa,\kappa
]}\bigl|\widetilde f_b(x+q_\tau)-f(x+q_\tau)\bigr|>
\eta r/2 \Bigr)+P\Bigl(\sup_{b\in
\mathcal B}|\widetilde
q_{\tau,b}-q_\tau|>\kappa\Bigr).\nonumber
\end{eqnarray}
Using $b_1=\min\mathcal B$, the first probability can be bounded by
\begin{eqnarray*}
&& \sum_{b\in\mathcal B}P \Bigl( \Bigl\{\sup
_{x\in[-\kappa,\kappa
]}\bigl|\widetilde f_b(x+q_\tau)-f(x+q_\tau)\bigr|>
\eta r/2 \Bigr\}\cap B_\varepsilon(b_1) \Bigr)+P
\bigl(B_\varepsilon(b_1)^c\bigr)
\\
&&\quad \lesssim \log n\sup_{b\in\mathcal B}P \Bigl( \Bigl\{\sup
_{x\in
[-\kappa,\kappa]}\bigl|\widetilde f_b(x+q_\tau)-f(x+q_\tau)\bigr|>
\eta r/2 \Bigr\}\cap B_\varepsilon(b_1) \Bigr)+\mathrm{o}(1) =\mathrm{o}(1),
\end{eqnarray*}
since for all $b$ the probability in the last line converges faster to
zero than $1/\log n$ owing to the concentration inequalities~(\ref
{eqConcentrationDensity}) and (\ref{eqconcBoundEps}) and the
conditions on $b$.
To estimate the second term in~(\ref{eqdecomp1}), we apply Lemma~\ref
{lemConsitency}. Therefore, the conditions $b\log(n)\to0$ and $(\log
n)^2/(nb^{2\beta+1})\to0$ yield the first assertion.

The estimate~(\ref{eqqDiffUniform}) follows from the error
decomposition~(\ref{eqtaylor}), (\ref{eqfDiffUniform}) and
Corollary~\ref{coruniformQ} with $x_n=\mathrm{o}(\inf_{\mathcal B}\delta_b)$:
\begin{eqnarray*}
&&P \bigl(\exists b\in\mathcal B\dvt |\widetilde q_{\tau,b}-q_\tau|>
\delta_b \bigl(3Db^{\alpha+1}+n^{-1/2}
\bigl(b^{-\beta+1/2}\vee1\bigr) \bigr) \bigr)
\\
&&\quad \leq P \biggl(\exists b\in\mathcal B\dvt  \biggl|\int_{-\infty
}^{q_\tau}
\widetilde f_b(x)-f(x)\,\d x \biggr|>\frac{1}{2}f(q_\tau)
\delta_b \bigl(3Db^{\alpha+1}+n^{-1/2}
\bigl(b^{-\beta+1/2}\vee1\bigr) \bigr) \biggr)
\\
&&\qquad{} +P \biggl(\sup_{b\in\mathcal B}\sup_{q^*_\tau\in[q_\tau\wedge
\widetilde q_{\tau,b},q_\tau\vee\widetilde q_{\tau,b}]}\bigl|\widetilde
f_b\bigl(q^*_\tau\bigr)-f(q_\tau)\bigr|>
\frac{1}{2} f(q_\tau) \biggr)
\\
&&\quad \lesssim \sum_{b\in\mathcal B} \biggl(\frac{1}{\delta_b}+
\frac
{1}{\delta_b^2}\frac{x_n}{mb^{1\wedge2\beta+2}\wedge1} \biggr)+\mathrm{o}(1)
\lesssim\sum
_{b\in\mathcal B}\frac{1}{\delta_b}+\mathrm{o}(1).
\end{eqnarray*}\upqed
\end{pf}

The variances $\sigma_{b,X}$ and $\sigma_{b,\varepsilon}$, defined
in (\ref{eqVxVar}) and (\ref{eqsigmaEps}) can be estimated by
$\widetilde\sigma_{b,X}$ and $\widetilde\sigma_{b,\varepsilon}$
from (\ref{eqTildeSigmaX}) and (\ref{eqTildeSigmaEPS}),
respectively. The latter can be decomposed into $\widetilde\sigma
^2_{b,\varepsilon}=\frac{1}{4}\pi^{-2}m^{-1}\widetilde\sigma
^2_{b,\varepsilon,1}\widetilde\sigma^2_{b,\varepsilon,2}$ with
\begin{eqnarray*}
\widetilde\sigma_{b,\varepsilon,1}^2&=&\int_{-1/b}^{1/b}\bigl|
\phi_K(bu)\bigr| \biggl|\frac{\phi_n(u)}{\phi_{\varepsilon,m}(u)} \biggr|^2\,\d u,
\\
\widetilde\sigma^2_{b,\varepsilon,2}&=&\int_{-1/b}^{1/b}\bigl|
\phi_K(bu)\bigr|\frac{|\F a_s(u)|^2}{|\phi_{\varepsilon,m}|^2}\,\d u.
\end{eqnarray*}
The following two lemmas show that these estimators are indeed reasonable.

%
\begin{lemma}\label{lemestSigX}
Let $\mathcal B$ be a finite set satisfying $|\mathcal B|\lesssim\log
n$, $\max_{b\in\mathcal B}b^\alpha\log n\to0$ as well as $\min_{b\in
\mathcal B}nb^{2\beta+2}\to\infty$. Let $\widetilde\sigma
_{b,X}$ and $\sigma_{b,X}$ be given in (\ref{eqTildeSigmaX}) and
(\ref{eqVxVar}), respectively. Then we obtain for all $\eta>0$ as
$n\to\infty$
\[
\sup_{f\in\mathcal C^\alpha(R,r,\zeta)}\sup_{f_\varepsilon\in
\mathcal D^\beta(R,\gamma)} P \bigl(
\exists b\in\mathcal B\dvt |\widetilde\sigma_{b,X}-\sigma_{b,X}|>
\eta m^{-1/2}\bigl(b^{-\beta
+1/2}\vee1\bigr) \bigr)\to0.
\]
\end{lemma}

\begin{pf}
Note that
%
\begin{eqnarray}\label{eqXiSums}
\widetilde\sigma^2_{b,X}&=&\frac{1}{n^2}\sum
_{j=1}^n\xi^2_{j,1}(b)+
\frac{1}{n^2}\sum_{j=1}^n
\xi^2_{j,2}(b)+\frac
{1}{n^2}\sum
_{j=1}^n\xi^2_{j,3}(b)
\nonumber\\[-8pt]\\[-8pt]
&&{}+\frac{2}{n^2}\sum_{j=1}^n
\xi_{j,1}(b)\xi_{j,2}(b)+\frac
{2}{n^2}\sum
_{j=1}^n\xi_{j,1}(b)\xi_{j,3}(b)+
\frac{2}{n^2}\sum_{j=1}^n
\xi_{j,2}(b)\xi_{j,3}(b),\nonumber
\end{eqnarray}
where we have defined
\begin{eqnarray*}
\xi_{j,1}(b)&:=& \int_{-\infty}^0a_s(x)
\F^{-1} \biggl[\phi_K(bu)\mathrm{e}^{\mathrm{i}uY_j} \biggl(
\frac{1}{\phi_{\varepsilon,m}(u)}-\frac{1}{\phi_\varepsilon(u)} \biggr)
\biggr](x+\widetilde
q_{\tau,b})\,\d x,
\\
\xi_{j,2}(b)&:=& \int_{-\infty}^0a_s(x)
\F^{-1} \biggl[\frac{\phi_K(bu)\mathrm{e}^{\mathrm{i}uY_j}}{\phi_\varepsilon(u)} \biggr
](x+q_{\tau})\,\d x,
\\
\xi_{j,3}(b)&:=& \int_{-\infty}^0a_s(x)
\F^{-1} \biggl[\frac{\phi_K(bu)\mathrm{e}^{\mathrm{i}uY_j}(\mathrm{e}^{-\mathrm{i}u\widetilde q_{\tau
,b}}-\mathrm{e}^{-\mathrm{i}u q_{\tau}})}{\phi_\varepsilon(u)} \biggr](x)\,\d x.
\end{eqnarray*}
We will first study these three terms separately.
Applying Plancherel's identity, the Cauchy--Schwarz inequality, the
Neumann type bound (\ref{eqNeumann}) as well as $|\F a_s(u)|\leq
A_s(1+|u|)^{-1}$, the decay of $\phi_\varepsilon$ and the upper bound
on $f$, we obtain
%
\begin{eqnarray}
\E\bigl[\bigl|\xi_{j,1}(b)\bigr|^2\ind_{B_\varepsilon(b)}\bigr]
&\leq&\frac{9}{2\pi^2}\int_{-1/b}^{1/b}
\frac{|\F
a_s(u)|^2}{|\phi_\varepsilon(u)|^2}\,\d u\int_{-1/b}^{1/b}
\frac{|\phi
_K(bu)|^2}{m|\phi_\varepsilon(u)|^2}\,\d u
\nonumber\\[-8pt]\label{eqXiJ1} \\[-8pt]
& \lesssim& \frac{1}{(b^{2\beta
-1}\wedge1)mb^{2\beta+1}},\nonumber
\\
\E\bigl[\bigl|\xi_{j,2}(b)\bigr|^2\bigr] &=&\E\biggl[ \biggl|
\frac{1}{2\pi}\int_{\R}\F a_s(u)\mathrm{e}^{-\mathrm{i}uq_\tau
}
\frac{\phi_K(bu)}{\phi_\varepsilon(u)}\mathrm{e}^{\mathrm{i}uY_j}\,\d u \biggr|^2 \biggr]
\nonumber\\[-8pt]\label{eqXiJ2} \\[-8pt]
&\leq&\frac{\|K\|_{L^1}^2A_s^2R^3}{4\pi^2}\int
_{-1/b}^{1/b}\bigl(1+|u|\bigr)^{2\beta-2}
\,\d u =:S^2_b\nonumber
\end{eqnarray}
as well as the deterministic bound
\begin{eqnarray}
\nonumber
\bigl|\xi_{j,2}(b)\bigr|^2&=& \biggl|\frac{1}{2\pi}\int
_{\R}\F a_s(u)\mathrm{e}^{-\mathrm{i}uq_\tau}
\frac{\phi
_K(bu)}{\phi_\varepsilon(u)}\mathrm{e}^{\mathrm{i}uY_j}\,\d u \biggr|^2 \leq
\frac{\|
K\|^2_{L_1}A_s^2}{4\pi^2}\int_{-1/b}^{1/b}\bigl(1+|u|\bigr)^{2\beta}
\,\d u =:d^2_b.
\end{eqnarray}
Hence, $\Var[\xi_{j,2}(b)^2]\leq\E[\xi_{j,2}(b)^4]\leq
d^2_bS^2_b$ and $|\xi^2_{j,2}(b)-\E[\xi_{j,2}^2(b)]|\leq
2d^2_b$, so that an application of Bernstein's inequality yields for
any $b>0$ and $z>0$
\begin{eqnarray*}
&&P \Biggl( \Biggl|\frac{1}{n}\sum_{j=1}^n
\bigl(\xi_{j,2}^2(b)-\E\bigl[\xi_{j,2}^2(b)
\bigr]\bigr) \Biggr|\geq z \Biggr) \leq2\exp\biggl(-\frac
{z^2n}{2S_b^2d_b^2+({4}/{3})d_b^2z}
\biggr).
\end{eqnarray*}
Setting $z=S^2_b$ and noting $S^2_b\lesssim(b^{-2\beta+1}\vee1),
d^2_b\lesssim b^{-2\beta}$, we see that
%
\begin{eqnarray}
P \Biggl( \Biggl|\frac{1}{n}\sum_{j=1}^n
\bigl(\xi_{j,2}^2(b)-\E\bigl[\xi_{j,2}^2(b)
\bigr]\bigr) \Biggr|\geq S_b^2 \Biggr) &\leq&2\exp
\biggl(-\frac{S_b^2n}{4d_b^2} \biggr) \leq2\exp\bigl(-Cnb^{2\beta
\wedge1}
\bigr)\label{eqXi2Bernstein}
\end{eqnarray}
for some $C>0$. The right-hand side of (\ref{eqXi2Bernstein}) tends
to zero with polynomial rate since $nb^{2\beta\wedge1}\gtrsim\log n$.

We use $\supp a_s\subset[-1,0]$ to write $\xi_{j,3}$ as
\begin{eqnarray*}
\xi_{j,3}(b)&=&\int_{\R} \bigl(a_s(x-
\widetilde q_{\tau,b})-a_s(x-q_\tau) \bigr)
\F^{-1} \biggl[\frac{\phi
_K(bu)\mathrm{e}^{\mathrm{i}uY_j}}{\phi_\varepsilon(u)} \biggr](x)\,\d x
\\
&\leq&\sup_{t\in(-1,0)}\bigl|a_s'(t)\bigr||
\widetilde q_{\tau,b}-q_\tau|\int_{(\widetilde q_{\tau,b}\wedge q_\tau
)-1}^{\widetilde q_{\tau,b}\vee q_\tau}
\biggl|\F^{-1} \biggl[\frac{\phi_K(bu)\mathrm{e}^{\mathrm{i}uY_j}}{\phi
_\varepsilon(u)} \biggr](x) \biggr|\,\d x.
\end{eqnarray*}
The Cauchy--Schwarz inequality and Plancherel's identity yield
\begin{eqnarray*}
\bigl|\xi_{j,3}(b)\bigr|^2 &\leq&\bigl\|a_s'
\ind_{(-1,0)}\bigr\|_\infty^2|\widetilde
q_{\tau,b}-q_\tau|^2\bigl(1+|\widetilde
q_{\tau,b}-q_\tau|\bigr)
\\
&&{}\times \int_{(\widetilde
q_{\tau,b}\wedge q_\tau)-1}^{\widetilde q_{\tau,b}\vee q_\tau}
\biggl|\F^{-1} \biggl[\frac{\phi_K(bu)\mathrm{e}^{\mathrm{i}uY_j}}{\phi_\varepsilon(u)} \biggr
](x) \biggr|^2\,\d x
\\
&\leq&\frac{\|a_s'\ind_{(-1,0)}\|_\infty^2}{2\pi}|\widetilde q_{\tau
,b}-q_\tau|^2\bigl(1+|
\widetilde q_{\tau,b}-q_\tau|\bigr)\int_{\R
} \biggl|
\frac{\phi_K(bu)}{\phi_\varepsilon(u)} \biggr|^2\,\d u
\\
&\lesssim& |\widetilde q_{\tau,b}-q_\tau|^2\bigl(1+|
\widetilde q_{\tau,b}-q_\tau|\bigr)b^{-2\beta-1}.
\end{eqnarray*}
By Lemma~\ref{lemConsitency} $\sup_{b\in\mathcal B}|\widetilde
q_{\tau,b}-q_\tau|=\mathrm{o}_P(1)$. Applying (\ref{eqqDiffUniform}), we
conclude for some constant $C>0$, for $\delta_b=(b^{\alpha+(1/2-\beta
)_+}+n^{-1/2}b^{-\beta-1/2}))^{-1}$ and for any $\eta>0$
%
\begin{eqnarray}\label{eqXiJ3}
&&P \bigl(\exists b\in\mathcal B\dvt \bigl|\xi_{j,3}(b)\bigr|>\eta
\bigl(b^{-\beta
+1/2}\vee1\bigr) \bigr)
\nonumber
\\
&&\quad \leq P \bigl(\exists b\in\mathcal B\dvt |\widetilde q_{\tau,b}-q_\tau|>
\eta Cb^{(\beta\wedge1/2)+1/2} \bigr)+\mathrm{o}(1)
\nonumber\\[-8pt]\\[-8pt]
&&\quad \leq P \bigl(\exists b\in\mathcal B\dvt |\widetilde q_{\tau,b}-q_\tau|>
\eta C\delta_b\bigl(b^{\alpha+1}+n^{-1/2}
\bigl(b^{-\beta
+1/2}\vee1\bigr)\bigr) \bigr)+\mathrm{o}(1)
\nonumber
\\
&&\quad \lesssim \biggl(\sum_{b\in\mathcal B}(\delta_b)^{-1}
\biggr)+\mathrm{o}(1) \lesssim\sup_{b\in\mathcal B}b^\alpha\log n+\sup
_{b\in\mathcal
B}\frac{\log n}{\sqrt nb^{\beta+1/2}}+\mathrm{o}(1)=\mathrm{o}(1).\nonumber
\end{eqnarray}
Combining the variance bounds (\ref{eqXiJ1}), (\ref{eqXiJ2}) and
(\ref{eqXiJ3}), we apply Markov's inequality, the Cauchy--Schwarz
inequality and the concentration result (\ref{eqXi2Bernstein}) on the
decomposition~(\ref{eqXiSums}) to obtain
\begin{eqnarray*}
&& \sup_{b\in\mathcal B} \bigl(n\bigl(b^{2\beta-1}\wedge1\bigr)\bigl|
\widetilde\sigma^2_{b,X}-\sigma^2_{b,X}\bigr|
\bigr)
\\
&&\quad  =\sup_{b\in\mathcal B} \Biggl(\frac{b^{2\beta-1}\wedge1}{n}\sum
_{j=1}^n\bigl(\xi_{j,2}^2(b)-
\E\bigl[\xi_{j,2}^2(b)\bigr]\bigr) \Biggr)+\mathrm{o}_P(1)
=\mathrm{o}_P(1).
\end{eqnarray*}\upqed
\end{pf}

%
\begin{lemma}\label{lemestSigEps}
Let $\mathcal B$ be a finite set satisfying $|\mathcal B|\lesssim\log
n$ as well as $\sup_{b\in\mathcal B}1/(nb^{2\beta+1})\to0$. Let
$\widetilde\sigma_{b,\varepsilon}$ and $\sigma_{b,\varepsilon}$ be
given in (\ref{eqTildeSigmaEPS}) and (\ref{eqsigmaEps}),
respectively. Then we obtain uniformly over $f\in\mathcal C^\alpha
(R,r,\zeta)$ and $f_\varepsilon\in\mathcal D^\beta(R,\gamma)$ for
all $\eta>0$ as $n\to\infty$
\[
P \bigl(\exists b\in\mathcal B\dvt |\widetilde\sigma_{b,\varepsilon
}-
\sigma_{b,\varepsilon}|>\eta(\log n) m^{-1/2}\bigl(b^{-\beta+1/2}\vee1
\bigr) \bigr)\to0.
\]
\end{lemma}

\begin{pf}
We start by showing for $b_1=\min\mathcal B$ that
%
\begin{equation}
\label{equniformConsistency} \sup_{|u|\leq1/b_1} \biggl|\frac{\phi
_\varepsilon(u)}{\phi
_{\varepsilon,m}(u)}
\biggr|=1+\mathrm{o}_P(1).
\end{equation}
To this end, recall $w(u)=(\log(\mathrm{e}+|u|))^{-1/2-\eta}$ for some $\eta
\in(0,1/2)$. Markov's inequality, Lemma~\ref{lemPhiEps} and Theorem
4.1 by Neumann and Rei{\ss} \cite{neumann2009nonparametric} yield for
any $\delta>0$
\begin{eqnarray*}
&&P \biggl(\sup_{|u|\leq1/b_1} \biggl|\frac{\phi_\varepsilon
(u)}{\phi_{\varepsilon,m}(u)}-1 \biggr|\geq\delta
\biggr)
\\
&&\quad \leq P \Bigl(\sup_{|u|\leq1/b_1}m^{1/2}\bigl|
\phi_\varepsilon(u)-\phi_{\varepsilon,m}(u)\bigr|\geq\delta|\log
b_1| \Bigr) +P \Bigl(\inf_{|u|\leq1/b_1}\bigl|
\phi_{\varepsilon,m}(u)\bigr|\leq m^{-1/2}|\log b_1| \Bigr)
\\
&&\quad \leq\bigl(\delta|\log b_1|\bigr)^{-1}\E\Bigl[\sup
_{|u|\leq
1/b_1}m^{1/2}\bigl|\phi_\varepsilon(u)-
\phi_{\varepsilon,m}(u)\bigr| \Bigr]+\mathrm{o}(1)
\\
&&\quad \leq\frac{1}{\delta|\log b_1|w(1/b_1)}\E\Bigl[\sup_{u\in\R
}m^{1/2}w(u)\bigl|
\phi_\varepsilon(u)-\phi_{\varepsilon,m}(u)\bigr| \Bigr]+\mathrm{o}(1) =\mathrm{o}(1),
\end{eqnarray*}
which implies (\ref{equniformConsistency}) holding uniformly in
$\mathcal B$ since $[-1/b_1,1/b_1]$ is the maximal interval for all
$b\in\mathcal B$.

Now, we consider $\widetilde\sigma_{b,\varepsilon,1}$. The uniform
consistency (\ref{equniformConsistency}) implies
\begin{eqnarray*}
\widetilde\sigma_{b,\varepsilon,1}^2&=&\bigl(1+\mathrm{o}_P(1)
\bigr)\int_{\R}\bigl|\phi_K(bu)\bigr| \biggl|
\frac{\phi_n(u)}{\phi_{\varepsilon}(u)} \biggr|^2\,\d u.
\end{eqnarray*}
Chebyshev's inequality yields for all $\eta>0$
\begin{eqnarray*}
&&P \biggl(\sup_{b\in\mathcal B} \biggl| \biggl(\int_{\R}\bigl|
\phi_K(bu)\bigr|\frac{|\phi_n(u)|^2}{|\phi_{\varepsilon}(u)|^2}\,\d u \biggr
)^{1/2}-\E
\biggl[ \biggl(\int_{\R}\bigl|\phi_K(bu)\bigr|
\frac{|\phi
_n(u)|^2}{|\phi_{\varepsilon}(u)|^2}\,\d u \biggr)^{1/2} \biggr] \biggr|>\eta
\log n \biggr)
\\
&&\quad\leq(\eta\log n)^{-2}\sum_{b\in\mathcal B}\E
\biggl[\int_{\R
}\bigl|\phi_K(bu)\bigr|\frac{|\phi_n(u)|^2}{|\phi_{\varepsilon}(u)|^2}
\,\d u \biggr]
\\
&&\quad \lesssim\bigl(\eta^2\log n\bigr)^{-1}\int
_{-1/b_1}^{1/b_1}\frac{\E[|\phi
_n(u)|^2]}{|\phi_{\varepsilon}(u)|^2}\,\d u \lesssim\bigl(
\eta^2\log n\bigr)^{-1},
\end{eqnarray*}
where the last estimate follows from $\E[|\phi_n(u)|^2]\lesssim|\phi
_Y(u)|^2+\E[|\phi_n(u)-\phi_Y(u)|^2]\lesssim|\phi_Y(u)|^2+1/n$,
$f_\varepsilon\in\mathcal D^\beta(R,\gamma), \|f\|_\infty\lesssim
1$ and $nb_1^{2\beta+1}\to\infty$.
Hence, we obtain uniformly in $\mathcal B$
%
\begin{eqnarray}
\widetilde\sigma_{b,\varepsilon,1} &=&\bigl(1+\mathrm{o}_P(1)\bigr) \bigl(
\sigma_{b,\varepsilon,1}+\mathrm{o}_P(\log n) \bigr) =\sigma_{b,\varepsilon,1}+\mathrm{o}_P(
\log n).\label{eqoPsigmaEps1}
\end{eqnarray}

Concerning $\widetilde\sigma_{b,\varepsilon,2}$, we write with use
of (\ref{equniformConsistency})
\begin{eqnarray*}
\widetilde\sigma_{b,\varepsilon,2}^2 &=&\int_{-1/b}^{1/b}\bigl|
\phi_K(bu)\bigr|\frac{|\F a_s(u)|^2}{|\phi
_{\varepsilon,m}(u)|^2}\,\d u =\bigl(1+\mathrm{o}_p(1)
\bigr)\int_{-1/b}^{1/b}\bigl|\phi_K(bu)\bigr|
\frac{|\F a_s(u)|^2}{|\phi
_{\varepsilon}(u)|^2}\,\d u.
\end{eqnarray*}
Moreover, the triangle inequality for the $L^2$-norm and Lemma~\ref
{lemPhiEps}, applied on $B_\varepsilon(b_1)$ yield
\begin{eqnarray*}
\hspace*{-5pt}&& \biggl| \biggl(\int_{-1/b}^{1/b}\bigl|\phi_K(bu)\bigr|
\frac{|\F a_s(u)|^2}{|\phi
_{\varepsilon}(u)|^2}\,\d u \biggr)^{1/2}-\sigma_{b,\varepsilon,2}
\biggr|^2
\\
\hspace*{-5pt}&&\quad \leq 2 \biggl|\E\biggl[ \biggl( \biggl(\int_{-1/b}^{1/b}\bigl|
\phi_K(bu)\bigr|\frac{|\F a_s(u)|^2}{|\phi_{\varepsilon}(u)|^2}\,\d u \biggr
)^{1/2}
\\
\hspace*{-5pt}&&\hspace*{46pt}{}- \biggl(
\int_{-1/b}^{1/b}\bigl|\phi_K(bu)\bigr|
\frac{|\F
a_s(u)|^2}{|\phi_{{\varepsilon,m}}(u)|^2}\,\d u \biggr)^{1/2} \biggr)\ind
_{B_\varepsilon(b_1)}
\biggr] \biggr|^2
\\
\hspace*{-5pt}&&\qquad{} +2P \bigl(\bigl(B_\varepsilon(b_1)\bigr)^c
\bigr)\int_{-1/b}^{1/b}\bigl|\phi_K(bu)\bigr|
\frac{|\F a_s(u)|^2}{|\phi_{\varepsilon}(u)|^2}\,\d u
\\
\hspace*{-5pt}&&\quad \leq 2\E\biggl[ \biggl(\int_{-1/b}^{1/b}\bigl|
\phi_K(bu)\bigr|\bigl|\F a_s(u)\bigr|^2\frac{|\phi_{\varepsilon,m}(u)-\phi_{\varepsilon
}(u)|^2}{|\phi_{\varepsilon}(u)\phi_{\varepsilon,m}(u)|^2}
\,\d u \biggr)\ind_{B_\varepsilon(b_1)} \biggr]
+\mathrm{o}(1)\int_{-1/b}^{1/b}
\frac{|\F
a_s(u)|^2}{|\phi_{\varepsilon}(u)|^2}\,\d u
\\
\hspace*{-5pt}&&\quad \leq\frac{2}{|\log b_1|^{3/2}}\E\biggl[\int_{-1/b}^{1/b}
\frac
{|\F a_s(u)|^2}{|\phi_\varepsilon(u)|^2}m\bigl|\phi_{\varepsilon,m}(u)-\phi
_{\varepsilon}(u)\bigr|^2
\,\d u \biggr]+\mathrm{o}(1) \bigl(b^{-2\beta+1}\vee1\bigr)
\\
\hspace*{-5pt}&&\quad =\mathrm{o}(1) \bigl(b^{-2\beta+1}
\vee1\bigr),
\end{eqnarray*}
where $\mathrm{o}(1)$ is a null sequence which does not depend on $b$. Consequently,
\[
\sup_{b\in\mathcal B} \biggl| \biggl(\int_{-1/b}^{1/b}\bigl|
\phi_K(bu)\bigr|\frac
{|\F a_s(u)|^2}{|\phi_{\varepsilon}(u)|^2}\,\d u \biggr)^{1/2}-\sigma
_{b,\varepsilon,2} \biggr|\bigl(b^{\beta-1/2}\wedge1\bigr)=\mathrm{o}(1).
\]
Using $\sigma_{b,\varepsilon,2}^2\lesssim b^{-2\beta+1}\vee1$ by
the analysis of the convergence rates, we get
%
\begin{eqnarray}
\widetilde\sigma_{b,\varepsilon,2} &=&\bigl(1+\mathrm{o}_p(1)\bigr) \bigl(
\sigma_{b,\varepsilon,2}+\mathrm{o} \bigl(b^{-\beta+1/2}\vee1 \bigr) \bigr) =
\sigma_{b,\varepsilon,2}+\mathrm{o}_P \bigl(b^{-\beta+1/2}\vee1
\bigr).\label{eqoPsigmaEps2}
\end{eqnarray}
Since $\sigma_{b,\varepsilon,1}\lesssim1, \sigma_{b,\varepsilon
,2}\lesssim b^{-\beta+1/2}\vee1$, it remains to combine (\ref
{eqoPsigmaEps1}) and (\ref{eqoPsigmaEps2}) to obtain uniformly in
$\mathcal B$
\begin{eqnarray*}
\widetilde\sigma_{b,\varepsilon}&=&\frac{1}{2\pi}m^{-1/2}\widetilde
\sigma_{b,\varepsilon,1}\widetilde\sigma_{b,\varepsilon,2} =\frac
{1}{2\pi}m^{-1/2}
\bigl(\sigma_{b,\varepsilon,1}+\mathrm{o}_P(\log n) \bigr) \bigl(
\sigma_{b,\varepsilon,2}+\mathrm{o}_P \bigl(b^{-\beta+1/2}\vee1 \bigr)
\bigr)
\\
&=&\sigma_{b,\varepsilon}+\mathrm{o}_P \bigl((\log n)m^{-1/2}
\bigl(b^{-\beta
+1/2}\vee1\bigr) \bigr).
\end{eqnarray*}\upqed
\end{pf}

%
\subsubsection{\texorpdfstring{Proof of Theorem~\protect\ref{thoracleinequality}}{Proof of Theorem 3.2}}
Applying Lemma~\ref{lemPhiEps} and (\ref{eqBminProb}), it suffices
to consider the event
\[
A_0:=\bigl\{b_{j_{0,n}}\leq b_{\widetilde j_n}\leq
n^{-1/(2\alpha
+2(\beta\vee1/2)+1)}\bigr\}\cap B_\varepsilon(b_{j_{0,n}})
\]
with $j_{0,n}$ defined in (\ref{eqj0}). Therefore we can set $\B:=\{
b_{j_{0,n}},\dots,b_{M_n}\}$ in the following.

As seen in error decomposition~(\ref{eqerrordecompadapt}), there
are three stochastic errors $V_{b,X}, V_{b,\varepsilon}$ and $V_{b,c}$
which were treated in Lemma~\ref{lemconcIneq}. This motivates the
following definition. For $\delta_1>0$, let
\begin{eqnarray*}
S_{b,X}:=(1+\delta_1)\sqrt{2\log\log n}\max
_{\mu\in\B\dvtx \mu
\geq b}\sigma_{\mu,X},\qquad S_{b,\varepsilon}:=(
\delta_1\log n)^3\max_{\mu\in\B\dvtx \mu\geq
b}
\sigma_{\mu,\varepsilon}.
\end{eqnarray*}
On the assumption $|\phi_\varepsilon(u)|\gtrsim(1+|u|)^{-\beta}$ we
obtain for $\sigma_{b,\varepsilon}=\frac{1}{2\pi}m^{-1/2}\sigma
_{b,\varepsilon,1}\sigma_{b,\varepsilon,2}$ from (\ref
{eqsigmaEps}) that
\[
\sigma_{b,\varepsilon,2}^2\gtrsim\int_{-1/b}^{1/b}\bigl|
\F a_s(-u)\bigr|^2\bigl(1+|u|\bigr)^{2\beta}\,\d u \gtrsim\int
_{-1/b}^{1/b}\bigl(1+|u|\bigr)^{2\beta-2}\,\d u\sim
b^{-2\beta
+1}\vee1.
\]
Also,\vspace*{1pt} we have $\sigma_{b,\varepsilon,1}=\|\phi_X\|
_{L^2}+\mathrm{o}(1)\geq\|\phi_X\|_{L^2}/2$ for $b$ small
enough and $n$ large enough. Thus, $\sigma_{b,\varepsilon}\gtrsim
m^{-1/2}(b^{-\beta+1/2}\vee1)$. Therefore, Lemma~\ref{lemconcIneq} yields
\begin{eqnarray*}
&&P \bigl(\exists b\in\mathcal B\dvt |V_{b,X}+V_{b,\varepsilon
}+V_{b,c}|
\geq S_{b,X}+S_{b,\varepsilon} \bigr)
\\
&&\quad \leq P \biggl(\exists b\in\mathcal B\dvt |V_{b,X}|\geq
S_{b,X}+\frac{1}{3}S_{b,\varepsilon} \biggr) +P \biggl(\exists b
\in\mathcal B\dvt |V_{b,\varepsilon}|\geq\frac
{S_{b,\varepsilon}}{3} \biggr)
\\
&&\qquad{} +P \biggl(
\exists b\in\mathcal B\dvt |V_{b,c}|\geq\frac{S_{b,\varepsilon}}{3}
\biggr)
\\
&&\quad = \mathrm{o}(1).
\end{eqnarray*}
Hence, the probability of the event
\[
A_1:= \bigl\{\forall b\in\mathcal B\dvt |V_{b,X}+V_{b,\varepsilon
}+V_{b,c}|
\leq S_{b,X}+S_{b,\varepsilon} \bigr\}
\]
converges to one. The variances $S_{b,X}$ and $S_{b,\varepsilon}$ can
be estimated by
\begin{eqnarray*}
\widetilde S_{b,X}:=(1+\delta_1)\sqrt{2\log\log n}\max
_{\mu\in\B\dvtx \mu\geq b}\widetilde\sigma_{\mu,X},\qquad\widetilde
S_{b,\varepsilon}:=(\delta_1\log n)^3\max
_{\mu\in\B\dvtx \mu\geq b}\widetilde\sigma_{\mu,\varepsilon}.
\end{eqnarray*}
Applying Lemmas~\ref{lemestSigX} and \ref{lemestSigEps}, the
triangle inequality of the $\ell^\infty$-norm yields uniformly in
$b\in\mathcal B$
\begin{eqnarray*}
\Bigl|\max_{\mu\geq b}\widetilde\sigma_{\mu,X}-\max
_{\mu
\geq b}\sigma_{\mu,X}\Bigr| &\leq&\max
_{\mu\geq b}|\widetilde\sigma_{\mu,X}-\sigma
_{\mu,X}|=\mathrm{o}_P \biggl(\frac{1}{m^{1/2}(b^{\beta-1/2}\wedge1)} \biggr),
\\
\Bigl|\max_{\mu\geq b}\widetilde\sigma_{\mu,\varepsilon}-\max
_{\mu\geq b}\sigma_{\mu,\varepsilon}\Bigr| &\leq&\max
_{\mu\geq b}|\widetilde\sigma_{\mu,\varepsilon
}-\sigma_{\mu,\varepsilon}|=\mathrm{o}_P
\biggl(\frac{\log n}{m^{1/2}(b^{\beta
-1/2}\wedge1)} \biggr).
\end{eqnarray*}
Using again $\sigma_{b,\varepsilon}\gtrsim m^{-1/2}(b^{-\beta
+1/2}\vee1)$, we thus obtain for all $\eta>0$ that the event
\begin{eqnarray*}
A_2&:=& \bigl\{\forall b\in\mathcal B\dvt  \bigl|(\widetilde
S_{b,X}+\widetilde S_{b,\varepsilon})-(S_{b,X}+S_{b,\varepsilon})
\bigr|\leq\eta(S_{b,X}+S_{b,\varepsilon}) \bigr\}
\end{eqnarray*}
fulfills $P(A_2)\to1$. The same holds true for the events
\begin{eqnarray*}
A_3&:=& \Bigl\{\forall b\in\mathcal B\dvt \sup_{q^*\in[(q_\tau\wedge
\widetilde q_{\tau,b})\vee(q_\tau\wedge\widetilde q_{\tau,b})]}\bigl|
\widetilde f_b\bigl(q^*\bigr)-f(q_\tau)\bigr|\leq\eta
f(q_\tau) \Bigr\},
\\
A_4&:=& \Bigl\{\forall b\in\mathcal B\dvt \sup_{q^*\in[(q_\tau\wedge
\widetilde q_{\tau,b})\vee(q_\tau\wedge\widetilde q_{\tau,b})]}\bigl|
\widetilde f_b\bigl(q^*\bigr)-\widetilde f_b(
\widetilde q_{\tau,b})\bigr|\leq\eta\bigl|\widetilde f_b(
\widetilde q_{\tau,b})\bigr| \Bigr\}
\end{eqnarray*}
by (\ref{eqfDiffUniform}). Therefore, it is sufficient to work in the
following on the event
\begin{eqnarray*}
A&:=&A_0\cap A_1\cap A_2\cap A_3
\cap A_4.
\end{eqnarray*}
We show that the adaptive estimator $\widetilde q_\tau$ mimics the
oracle estimator defined as follows. Recalling the estimate of the bias
$B_b=Db^{\alpha+1}$, let the oracle bandwidth be defined by
%
\begin{eqnarray}
\label{eqoracleband} b_*:=\max\{b\in\B\dvt B_b\leq
S_{b,X}+S_{b,\varepsilon}
\}.
\end{eqnarray}
Note that $b_*$ is well-defined and unique since $B_b$ is monoton
increasing in $b$ while $(S_{b,X}+S_{b,\varepsilon})$ is monton
decreasing. We get the oracle estimator $\widetilde q_{\tau,b_*}$.

Since on $A_4$ for all $b\in\mathcal B$ and $q^*\in[(q_\tau\wedge
\widetilde q_{\tau,b})\vee(q_\tau\wedge\widetilde q_{\tau,b})]$
\[
\bigl|\widetilde f_b\bigl(q^*\bigr)\bigr|\geq\bigl|\widetilde
f_b(\widetilde q_{\tau,b})\bigr|- \bigl|\widetilde f_b
\bigl(q^*\bigr)-\widetilde f_b(\widetilde q_{\tau,b})\bigr|
\geq(1-\eta)\bigl|\widetilde f_b(\widetilde q_{\tau,b})\bigr|,
\]
we have for any $b\in\B$ on the event $A_1\cap A_4$ by (\ref
{eqerrordecompadapt})
\begin{eqnarray*}
|\widetilde q_{\tau,b}-q_\tau| &\leq&\frac
{B_b+|V_{b,X}+V_{b,\varepsilon}+V_{b,c}|
{+|\widetilde M_b(\widetilde q_{\tau,b})|}}{|\widetilde{f}_b(q^*)|}
\leq\frac{B_b+S_{b,X}+S_{b,\varepsilon} {+|\widetilde
M_b(\widetilde q_{\tau,b})|}}{(1-\eta)|\widetilde{f}_b(\widetilde
q_{\tau,b})|}.
\end{eqnarray*}
Furthermore, by the definition of $b_*$ we have on the event $A$ for
any $b\leq b_*$
\begin{eqnarray*}
|\widetilde q_{\tau,b}-q_\tau|&\leq&\frac
{2(S_{b,X}+S_{b,\varepsilon})+|\widetilde M_b(\widetilde q_{\tau
,b})|}{(1-\eta)|\widetilde{f}_b(\widetilde q_{\tau,b})|}.
\end{eqnarray*}
On $A_2$ we estimate $\widetilde S_{b,X}+\widetilde S_{b,\varepsilon
}\geq(1-\eta)(S_{b,X}+S_{b,\varepsilon})$ and thus we have on
$A$ for any $b\leq b_*$
\begin{eqnarray*}
|\widetilde q_{\tau,b}-q_\tau|&\leq&\frac{2(\widetilde
S_{b,X}+\widetilde S_{b,\varepsilon})}{(1-\eta
)^2|\widetilde{f}_b(\widetilde q_{\tau,b})|}+
\frac{|\widetilde
M_b(\widetilde q_{\tau,b})|}{(1-\eta)|\widetilde f_b(\widetilde
q_{\tau,b})|}.
\end{eqnarray*}
Since for any $\delta>0$ we find $\delta_1,\eta>0$ such that $
((1-\eta)^{-2}(2\sqrt2+\delta_1)-2\sqrt2 )\vee(2(1-\eta
)^{-2}\delta_1 )\vee\frac{\eta}{1-\eta}<\delta$, we obtain
$|\widetilde q_{\tau,b}-q_\tau|\leq\widetilde\Sigma_b$ with
$\widetilde\Sigma_b$ as defined in (\ref{eqTildeSigma}). As a
result one has $q_\tau\in\mathcal U_b$ and $q_\tau\in\mathcal U_\mu
$ for all $b\leq b_*$ and $\mu\leq b_*$,\vspace*{1pt} implying $\mathcal
U_\mu\cap\mathcal U_b\ne\varnothing$. By the definition of the
procedure, $\widetilde{b}^*\geq b_*$ and $\mathcal U_{\widetilde
{b}_*}\cap\mathcal U_{b_*}\ne\varnothing$ on the event $A$. This leads to
\begin{eqnarray*}
|\widetilde q_{\tau,{\widetilde{b}^*}}-q_\tau| &\leq&|\widetilde
q_{\tau,b_*}-q_\tau|+ |\widetilde q_{\tau,{\widetilde{b}^*}} -\widetilde
q_{\tau,b_*}| \leq\widetilde\Sigma_{b_*}+(\widetilde
\Sigma_{b_*}+ \widetilde\Sigma_{\widetilde{b}^*}).
\end{eqnarray*}
On $A_2\cap A_3$ we have $\widetilde\Sigma_{b}\lesssim
S_{b,X}+S_{b,\varepsilon}$ since $f(q_\tau)\geq r$ {and
$|\widetilde M_{b}(\widetilde q_{\tau,b})|\leq|\widetilde
M_b(q_\tau)|=|\int_{-\infty}^{q_\tau}(\widetilde f_b-f)|$}. Using
additionally the monotonicity of $(S_{b,X}+S_{b,\varepsilon})$ as well
as $\widetilde{b}^*\geq b_*$, this implies
\[
|\widetilde q_{\tau,{\widetilde{b}^*}}-q_\tau|\lesssim
(S_{b_*,X}+S_{b_*,\varepsilon})
\lesssim\bigl(\sqrt{\log\log n}+ \bigl(\log n^{\delta}
\bigr)^3 \bigr) \bigl(b_*^{-\beta+1/2}\vee1 \bigr)n^{-1/2}.
\]
It remains to note by the definition (\ref{eqoracleband}) of the
oracle $b_*$ and by the assumption $b_{j+1}/b_j\lesssim1$ that
{$b_*\sim((\log n^\delta)^6/n)^{-1/(2\alpha+2(\beta\vee1/2)+1)}$}
as $n\to\infty$.


\section*{Acknowledgements}
This research started when the first author was a Postdoc at EURANDOM,
Eindhoven University of Technology, The Netherlands. The research was
partly supported by the Deutsche Forschungsgemeinschaft through the FOR
1735 ``Structural Inference in Statistics''. Part of the work on this
project was done during a visit of the third author to EURANDOM. We
thank three anonymous referees for helpful comments and suggestions.


%

\printhistory

\begin{thebibliography}{31}
\bibitem{carroll1988optimal}
%
\begin{barticle}[mr]
\bauthor{\bsnm{Carroll},~\bfnm{Raymond~J.}\binits{R.J.}} \AND
\bauthor{\bsnm{Hall},~\bfnm{Peter}\binits{P.}}
(\byear{1988}).
\btitle{Optimal rates of convergence for deconvolving a density}.
\bjournal{J. Amer. Statist. Assoc.}
\bvolume{83}
\bpages{1184--1186}.
\bid{issn={0162-1459}, mr={0997599}}
\end{barticle}
%
\bptok{imsref}%
\endbibitem

\bibitem{carroll2006measurement}
%
\begin{bbook}[mr]
\bauthor{\bsnm{Carroll},~\bfnm{Raymond~J.}\binits{R.J.}},
\bauthor{\bsnm{Ruppert},~\bfnm{David}\binits{D.}},
\bauthor{\bsnm{Stefanski},~\bfnm{Leonard~A.}\binits{L.A.}} \AND
\bauthor{\bsnm{Crainiceanu},~\bfnm{Ciprian~M.}\binits{C.M.}}
(\byear{2006}).
\btitle{Measurement Error in Nonlinear Models: A Modern Perspective},
\bedition{2nd} ed.
\bseries{Monographs on Statistics and Applied Probability}
\bvolume{105}.
\blocation{Boca Raton, FL}:
\bpublisher{Chapman \& Hall/CRC}.
\bid{doi={10.1201/9781420010138}, mr={2243417}}
\end{bbook}
%
\bptok{imsref}%
\endbibitem

\bibitem{comte2011data}
%
\begin{barticle}[mr]
\bauthor{\bsnm{Comte},~\bfnm{F.}\binits{F.}} \AND
\bauthor{\bsnm{Lacour},~\bfnm{C.}\binits{C.}}
(\byear{2011}).
\btitle{Data-driven density estimation in the presence of additive
noise with unknown distribution}.
\bjournal{J. R. Stat. Soc. Ser. B Stat. Methodol.}
\bvolume{73}
\bpages{601--627}.
\bid{doi={10.1111/j.1467-9868.2011.00775.x}, issn={1369-7412}, mr={2853732}}
\end{barticle}
%
\bptok{imsref}%
\endbibitem

\bibitem{dattnerEtAll2011}
%
\begin{barticle}[mr]
\bauthor{\bsnm{Dattner},~\bfnm{I.}\binits{I.}},
\bauthor{\bsnm{Goldenshluger},~\bfnm{A.}\binits{A.}} \AND
\bauthor{\bsnm{Juditsky},~\bfnm{A.}\binits{A.}}
(\byear{2011}).
\btitle{On deconvolution of distribution functions}.
\bjournal{Ann. Statist.}
\bvolume{39}
\bpages{2477--2501}.
\bid{doi={10.1214/11-AOS907}, issn={0090-5364}, mr={2906875}}
\end{barticle}
%
\bptok{imsref}%
\endbibitem

\bibitem{dattner2012estimation}
%
\begin{barticle}[mr]
\bauthor{\bsnm{Dattner},~\bfnm{I.}\binits{I.}} \AND
\bauthor{\bsnm{Reiser},~\bfnm{B.}\binits{B.}}
(\byear{2013}).
\btitle{Estimation of distribution functions in measurement error models}.
\bjournal{J.~Statist. Plann. Inference}
\bvolume{143}
\bpages{479--493}.
\bid{doi={10.1016/j.jspi.2012.09.004}, issn={0378-3758}, mr={2995109}}
\end{barticle}
%
\bptok{imsref}%
\endbibitem

\bibitem{delaigle2008deconvolution}
%
\begin{barticle}[mr]
\bauthor{\bsnm{Delaigle},~\bfnm{Aurore}\binits{A.}},
\bauthor{\bsnm{Hall},~\bfnm{Peter}\binits{P.}} \AND
\bauthor{\bsnm{Meister},~\bfnm{Alexander}\binits{A.}}
(\byear{2008}).
\btitle{On deconvolution with repeated measurements}.
\bjournal{Ann. Statist.}
\bvolume{36}
\bpages{665--685}.
\bid{doi={10.1214/009053607000000884}, issn={0090-5364}, mr={2396811}}
\end{barticle}
%
\bptok{imsref}%
\endbibitem

\bibitem{dudley1992}
%
\begin{barticle}[mr]
\bauthor{\bsnm{Dudley},~\bfnm{R.~M.}\binits{R.M.}}
(\byear{1992}).
\btitle{Fr\'echet differentiability, {$p$}-variation and uniform
{D}onsker classes}.
\bjournal{Ann. Probab.}
\bvolume{20}
\bpages{1968--1982}.
\bid{issn={0091-1798}, mr={1188050}}
\end{barticle}
%
\bptok{imsref}%
\endbibitem

\bibitem{fan1991}
%
\begin{barticle}[mr]
\bauthor{\bsnm{Fan},~\bfnm{Jianqing}\binits{J.}}
(\byear{1991}).
\btitle{On the optimal rates of convergence for nonparametric
deconvolution problems}.
\bjournal{Ann. Statist.}
\bvolume{19}
\bpages{1257--1272}.
\bid{doi={10.1214/aos/1176348248}, issn={0090-5364}, mr={1126324}}
\end{barticle}
%
\bptok{imsref}%
\endbibitem

\bibitem{frese2011blood}
%
\begin{barticle}[pbm]
\bauthor{\bsnm{Frese},~\bfnm{Ethel~M.}\binits{E.M.}},
\bauthor{\bsnm{Fick},~\bfnm{Ann}\binits{A.}} \AND
\bauthor{\bsnm{Sadowsky},~\bfnm{H.~Steven}\binits{H.S.}}
(\byear{2011}).
\btitle{Blood pressure measurement guidelines for physical therapists}.
\bjournal{Cardiopulm. Phys. Ther. J.}
\bvolume{22}
\bpages{5--12}.
\bid{issn={1541-7891}, pmcid={3104931}, pmid={21637392}}
\end{barticle}
%
\bptok{imsref}%
\endbibitem

\bibitem{girardi2003operator}
%
\begin{barticle}[mr]
\bauthor{\bsnm{Girardi},~\bfnm{Maria}\binits{M.}} \AND
\bauthor{\bsnm{Weis},~\bfnm{Lutz}\binits{L.}}
(\byear{2003}).
\btitle{Operator-valued {F}ourier multiplier theorems on {B}esov spaces}.
\bjournal{Math. Nachr.}
\bvolume{251}
\bpages{34--51}.
\bid{doi={10.1002/mana.200310029}, issn={0025-584X}, mr={1960803}}
\end{barticle}
%
\bptok{imsref}%
\endbibitem

\bibitem{goldenshluger1997spatially}
%
\begin{barticle}[mr]
\bauthor{\bsnm{Goldenshluger},~\bfnm{A.}\binits{A.}} \AND
\bauthor{\bsnm{Nemirovski},~\bfnm{A.}\binits{A.}}
(\byear{1997}).
\btitle{On spatially adaptive estimation of nonparametric regression}.
\bjournal{Math. Methods Statist.}
\bvolume{6}
\bpages{135--170}.
\bid{issn={1066-5307}, mr={1466625}}
\end{barticle}
%
\bptok{imsref}%
\endbibitem

\bibitem{hall2008estimation}
%
\begin{barticle}[mr]
\bauthor{\bsnm{Hall},~\bfnm{Peter}\binits{P.}} \AND
\bauthor{\bsnm{Lahiri},~\bfnm{Soumendra~N.}\binits{S.N.}}
(\byear{2008}).
\btitle{Estimation of distributions, moments and quantiles in
deconvolution problems}.
\bjournal{Ann. Statist.}
\bvolume{36}
\bpages{2110--2134}.
\bid{doi={10.1214/07-AOS534}, issn={0090-5364}, mr={2458181}}
\end{barticle}
%
\bptok{imsref}%
\endbibitem

\bibitem{johannes2009deconvolution}
%
\begin{barticle}[mr]
\bauthor{\bsnm{Johannes},~\bfnm{Jan}\binits{J.}}
(\byear{2009}).
\btitle{Deconvolution with unknown error distribution}.
\bjournal{Ann. Statist.}
\bvolume{37}
\bpages{2301--2323}.
\bid{doi={10.1214/08-AOS652}, issn={0090-5364}, mr={2543693}}
\end{barticle}
%
\bptok{imsref}%
\endbibitem

\bibitem{Johannes2010}
%
\begin{barticle}[mr]
\bauthor{\bsnm{Johannes},~\bfnm{Jan}\binits{J.}} \AND
\bauthor{\bsnm{Schwarz},~\bfnm{Maik}\binits{M.}}
(\byear{2013}).
\btitle{Adaptive circular deconvolution by model selection under
unknown error distribution}.
\bjournal{Bernoulli}
\bvolume{19}
\bpages{1576--1611}.
\bid{doi={10.3150/12-BEJ422}, issn={1350-7265}, mr={3129026}}
\end{barticle}
%
\bptok{imsref}%
\endbibitem

\bibitem{kannel1995framingham}
%
\begin{barticle}[author]
\bauthor{\bsnm{Kannel},~\bfnm{W.~B.}\binits{W.B.}}
(\byear{1995}).
\btitle{Framingham study insights into hypertensive risk of
cardiovascular disease.}
\bjournal{Hypertension Research: Official Journal of the Japanese
Society of Hypertension}
\bvolume{18}
\bpages{181--196}.
\end{barticle}
%
\bptok{imsref}%
\endbibitem

\bibitem{kappus2012}
%
\begin{barticle}[mr]
\bauthor{\bsnm{Kappus},~\bfnm{Johanna}\binits{J.}}
(\byear{2014}).
\btitle{Adaptive nonparametric estimation for {L}\'evy processes
observed at low frequency}.
\bjournal{Stochastic Process. Appl.}
\bvolume{124}
\bpages{730--758}.
\bid{doi={10.1016/j.spa.2013.08.010}, issn={0304-4149}, mr={3131312}}
\end{barticle}
%
\bptok{imsref}%
\endbibitem

\bibitem{Lepski1990}
%
\begin{barticle}[mr]
\bauthor{\bsnm{Lepski{\u\i}},~\bfnm{O.~V.}\binits{O.V.}}
(\byear{1990}).
\btitle{A problem of adaptive estimation in {G}aussian white noise}.
\bjournal{Teor. Veroyatn. Primen.}
\bvolume{35}
\bpages{459--470}.
\bid{doi={10.1137/1135065}, issn={0040-361X}, mr={1091202}}
\end{barticle}
%
\bptok{imsref}%
\endbibitem

\bibitem{louniciNickl2011}
%
\begin{barticle}[mr]
\bauthor{\bsnm{Lounici},~\bfnm{Karim}\binits{K.}} \AND
\bauthor{\bsnm{Nickl},~\bfnm{Richard}\binits{R.}}
(\byear{2011}).
\btitle{Global uniform risk bounds for wavelet deconvolution estimators}.
\bjournal{Ann. Statist.}
\bvolume{39}
\bpages{201--231}.
\bid{doi={10.1214/10-AOS836}, issn={0090-5364}, mr={2797844}}
\end{barticle}
%
\bptok{imsref}%
\endbibitem

\bibitem{massart2007}
%
\begin{bbook}[mr]
\bauthor{\bsnm{Massart},~\bfnm{Pascal}\binits{P.}}
(\byear{2007}).
\btitle{Concentration Inequalities and Model Selection}.
\bseries{Lecture Notes in Math.}
\bvolume{1896}.
\blocation{Berlin}:
\bpublisher{Springer}.
\bid{mr={2319879}}
\end{bbook}
%
\bptok{imsref}%
\endbibitem

\bibitem{meister2004effect}
%
\begin{barticle}[mr]
\bauthor{\bsnm{Meister},~\bfnm{Alexander}\binits{A.}}
(\byear{2004}).
\btitle{On the effect of misspecifying the error density in a
deconvolution problem}.
\bjournal{Canad. J. Statist.}
\bvolume{32}
\bpages{439--449}.
\bid{doi={10.2307/3316026}, issn={0319-5724}, mr={2125855}}
\end{barticle}
%
\bptok{imsref}%
\endbibitem

\bibitem{neumann1997effect}
%
\begin{barticle}[mr]
\bauthor{\bsnm{Neumann},~\bfnm{Michael~H.}\binits{M.H.}}
(\byear{1997}).
\btitle{On the effect of estimating the error density in nonparametric
deconvolution}.
\bjournal{J. Nonparametr. Stat.}
\bvolume{7}
\bpages{307--330}.
\bid{doi={10.1080/10485259708832708}, issn={1048-5252}, mr={1460203}}
\end{barticle}
%
\bptok{imsref}%
\endbibitem

\bibitem{neumann2007deconvolution}
%
\begin{barticle}[mr]
\bauthor{\bsnm{Neumann},~\bfnm{Michael~H.}\binits{M.H.}}
(\byear{2007}).
\btitle{Deconvolution from panel data with unknown error distribution}.
\bjournal{J. Multivariate Anal.}
\bvolume{98}
\bpages{1955--1968}.
\bid{doi={10.1016/j.jmva.2006.09.012}, issn={0047-259X}, mr={2396948}}
\end{barticle}
%
\bptok{imsref}%
\endbibitem

\bibitem{neumann2009nonparametric}
%
\begin{barticle}[mr]
\bauthor{\bsnm{Neumann},~\bfnm{Michael~H.}\binits{M.H.}} \AND
\bauthor{\bsnm{Rei{\ss}},~\bfnm{Markus}\binits{M.}}
(\byear{2009}).
\btitle{Nonparametric estimation for {L}\'evy processes from
low-frequency observations}.
\bjournal{Bernoulli}
\bvolume{15}
\bpages{223--248}.
\bid{doi={10.3150/08-BEJ148}, issn={1350-7265}, mr={2546805}}
\end{barticle}
%
\bptok{imsref}%
\endbibitem

\bibitem{NicklReiss2012}
%
\begin{barticle}[mr]
\bauthor{\bsnm{Nickl},~\bfnm{Richard}\binits{R.}} \AND
\bauthor{\bsnm{Rei{\ss}},~\bfnm{Markus}\binits{M.}}
(\byear{2012}).
\btitle{A {D}onsker theorem for {L}\'evy measures}.
\bjournal{J. Funct. Anal.}
\bvolume{263}
\bpages{3306--3332}.
\bid{doi={10.1016/j.jfa.2012.08.012}, issn={0022-1236}, mr={2973342}}
\end{barticle}
%
\bptok{imsref}%
\endbibitem

\bibitem{rosenthal1970}
%
\begin{barticle}[mr]
\bauthor{\bsnm{Rosenthal},~\bfnm{Haskell~P.}\binits{H.P.}}
(\byear{1970}).
\btitle{On the subspaces of {$L\sp{p}$} {$(p>2)$} spanned by sequences
of independent random variables}.
\bjournal{Israel J. Math.}
\bvolume{8}
\bpages{273--303}.
\bid{issn={0021-2172}, mr={0271721}}
\end{barticle}
%
\bptok{imsref}%
\endbibitem

\bibitem{soehlTrabs2012}
%
\begin{barticle}[mr]
\bauthor{\bsnm{S{\"o}hl},~\bfnm{Jakob}\binits{J.}} \AND
\bauthor{\bsnm{Trabs},~\bfnm{Mathias}\binits{M.}}
(\byear{2012}).
\btitle{A uniform central limit theorem and efficiency for
deconvolution estimators}.
\bjournal{Electron. J. Stat.}
\bvolume{6}
\bpages{2486--2518}.
\bid{doi={10.1214/12-EJS757}, issn={1935-7524}, mr={3020273}}
\end{barticle}
%
\bptok{imsref}%
\endbibitem

\bibitem{spokoiny1996}
%
\begin{barticle}[mr]
\bauthor{\bsnm{Spokoiny},~\bfnm{V.~G.}\binits{V.G.}}
(\byear{1996}).
\btitle{Adaptive hypothesis testing using wavelets}.
\bjournal{Ann. Statist.}
\bvolume{24}
\bpages{2477--2498}.
\bid{doi={10.1214/aos/1032181163}, issn={0090-5364}, mr={1425962}}
\end{barticle}
%
\bptok{imsref}%
\endbibitem

\bibitem{stirnemann2012density}
%
\begin{barticle}[mr]
\bauthor{\bsnm{Stirnemann},~\bfnm{J.~J.}\binits{J.J.}},
\bauthor{\bsnm{Comte},~\bfnm{F.}\binits{F.}} \AND
\bauthor{\bsnm{Samson},~\bfnm{A.}\binits{A.}}
(\byear{2012}).
\btitle{Density estimation of a biomedical variable subject to
measurement error using an auxiliary set of replicate observations}.
\bjournal{Stat. Med.}
\bvolume{31}
\bpages{4154--4163}.
\bid{doi={10.1002/sim.5392}, issn={0277-6715}, mr={3040071}}
\end{barticle}
%
\bptok{imsref}%
\endbibitem

\bibitem{triebel2010}
%
\begin{bbook}[mr]
\bauthor{\bsnm{Triebel},~\bfnm{Hans}\binits{H.}}
(\byear{2010}).
\btitle{Theory of Function Spaces}.
\bseries{Modern Birkh\"auser Classics}.
\blocation{Basel}:
\bpublisher{Birkh\"auser}.
\bnote{Reprint of 1983 edition [MR0730762], Also published in 1983 by
Birkh{\"a}user Verlag [MR0781540]}.
\bid{mr={3024598}}
\end{bbook}
%
\bptok{imsref}%
\endbibitem

\bibitem{tsybakov2009}
%
\begin{bbook}[mr]
\bauthor{\bsnm{Tsybakov},~\bfnm{Alexandre~B.}\binits{A.B.}}
(\byear{2009}).
\btitle{Introduction to Nonparametric Estimation}.
\bseries{Springer Series in Statistics}.
\blocation{New York}:
\bpublisher{Springer}.
\bnote{Revised and extended from the 2004 French original, Translated
by Vladimir Zaiats}.
\bid{doi={10.1007/b13794}, mr={2724359}}
\end{bbook}
%
\bptok{imsref}%
\endbibitem

\bibitem{vanderVaart1998}
%
\begin{bbook}[mr]
\bauthor{\bsnm{van~der Vaart},~\bfnm{A.~W.}\binits{A.W.}}
(\byear{1998}).
\btitle{Asymptotic Statistics}.
\bseries{Cambridge Series in Statistical and Probabilistic Mathematics}
\bvolume{3}.
\blocation{Cambridge}:
\bpublisher{Cambridge Univ. Press}.
\bid{doi={10.1017/CBO9780511802256}, mr={1652247}}
\end{bbook}
%
\bptok{imsref}%
\endbibitem
\end{thebibliography}
\end{document}